\documentclass[stand,english,12pt]{article}
\usepackage[centertags]{amsmath}
\usepackage{amsfonts}
\usepackage{amssymb}
\usepackage{amsthm}
\usepackage{newlfont}

\def\ds{\displaystyle}
\def\be{\begin{equation}}
\def\ee{\end{equation}}

\def\Dt{\Delta}
\def\pt{\partial}

\def\up{\stackrel}

\def\lms{\limits}

\begin{document}


 \begin{center}
 \textbf{Yu.~N.~Bratkov}\\
 \textbf{
 The hyperbolic Monge--Ampere equation:\\
 classical solutions on the whole plane
 }
 \end{center}


\begin{abstract}
\noindent
 The Cauchy problem for the hyperbolic Monge--Ampere equation
$$
 \left\{
 \begin{array}{l}
   A + Bz_{xx} + Cz_{xy} + Dz_{yy} + \mbox{hess}\: z = 0 \, ,
   \medskip \\
   z(0,y) = z^o(y),
   \quad
   z_x(0,y) = p^{\: o}(y), \ \ y \in \mathbb{R}
 \end{array}
 \right.
$$
 is considered.
 Here
 ${\mbox{hess}\: z = z_{xx}z_{yy} - z_{xy}^2}$,
 $A,B,C,D$
 depends on
 $x,y,z,z_x,z_y$.
 The equation is hyperbolic when
 ${C^2 - 4BD + 4A > 0}$.
 Sufficient conditions
 on the existence of a (unique)
 $C^3\mbox{-solution}$ on the whole plain are formulated.
\end{abstract}

\tableofcontents

\section{Introduction}

 This well-known problem was posed to the author by E.~V.~Shikin.
 The first publication was made in \cite{Bratkov}
 (it was received by the journal in February 1998). The author
 is grateful to D.~V.~Tunitsky for finding some subtle
 (not essential) error. Here we publish the full text of
 the revised proof.
 Changing is getting another successive approximations.
 The result doesn't depend on this changing.
 Certainly, this paper was a good
 reason for extending and improving the result.

 The Monge--Ampere equation could be reduced to the system of
 five equations in Riemann invariants
  \cite{Tunitsky-dep}.
 The theory of hyperbolic systems is a perfect one, when
 eigenvalues are separable (for example, are separated by
 constants). In our case
 eigenvalues are the solution of the system,
 i.e. unknown functions. It is required to find them.

 The essence of the problem was formulated by J.~Leray in his
 Princeton lectures on hyperbolic equations (1953):
 "It turns out well to prove
 the \textit{local existence theorem} only\dots
 It shows that for hyperbolic equations
 the existence of solutions on the whole depends on
 getting \textit{a priory estimations} for their derivatives\dots
 Except of equations which are linear outside some small area,
 we don't have examples of such known a priory estimations"\
 \cite{Leray}, Chapter IX.
 Comments of N.~H.~Ibragimov (1984):
 "Now we know some nonlinear equations for which the Cauchy problem
 is solvable on the whole. For example, for the Yang-Mills
 equations the theorem of the existence of the solution on the whole
 is proved by two different ways"\ (\cite{Leray},
 Russian edition, p.~162).

 Global classical solvability is a sophisticated kind of sport.
 B. Riemann proved the nonexistence of global classical solutions
 for some system of hyperbolic equations \cite{Yanenko}.
 This system isn't weakly
 nonlinear. Weakly nonlinear
 systems were introduced by N. N. Yanenko in 1955 \cite{Yanenko}.
 The system of two equations in Riemann invariants
$$
      \left(
        \partial_x
        +
        \xi_1
        (u_1, u_2)
        \:
        \partial_y
      \right)
      u_1
      =
      0,
   \qquad
      \left(
        \partial_x
        +
        \xi_2
        (u_1, u_2)
        \:
        \partial_y
      \right)
      u_2
      =
      0
$$
 was considered. Weak nonlinearity
$$
 \pt \xi_i / \pt u_i = 0,
 \quad i=1,2,
$$
 is required
 for global classical solvability.

 The system
$$
 \begin{array}{l}
      \left(
        \partial_x
        +
        \xi_1
        (x,y,u_2)
        \:
        \partial_y
      \right)
      u_1
      =
      f_1(x,y,u_1, u_2),
   \medskip\\
      \left(
        \partial_x
        +
        \xi_2
        (x,y,u_1)
        \:
        \partial_y
      \right)
      u_2
      =
      f_2(x,y,u_1, u_2)
 \end{array}
$$
 was considered by B.~L.~Rozhdestvensky and A.~D.~Sydorenko in 1967.
 This weakly nonlinear system has global classical solutions, when it is
 hyperbolic in the restricted sense. Restricted hyperbolicity
 means separability of eigenvalues of the system.
 Thus by the Rozhdestvensky---Sydorenko theorem
 (\cite{RYa}, Chapter~1, $\S$ 10, Subsection~3)
 the problem
 of a priory estimations for derivatives is reduced to the problem of
 separability of eigenvalues, or to the problem
 $ \xi_1 \neq \xi_2 $.

 Consider an example. Let the coefficients of the Monge--Ampere
 equation depend on ${x,y}$ only. The equation could be reduced to the
 system
$$
 \begin{array}{l}
      \left(
        \partial_x
        +
        u_2
        \:
        \partial_y
      \right)
      u_1
      =
      (u_1 - u_2)\, a_1(x,y,u_1),
   \medskip\\
      \left(
        \partial_x
        +
        u_1
        \:
        \partial_y
      \right)
      u_2
      =
      (u_1 - u_2)\, a_2(x,y,u_2).
 \end{array}
$$
 From this system, subtracting and integrating, we have
$$
 (u_1 - u_2)(x,y)
 =
 (u_1^o - u_2^o)
 \exp
 \left\{
    {\int_0^x}
      (a_1 - a_2 - u_{2y})
    \: d\tau
 \right\}.
$$
 If ${u_1^o(y) \neq u_2^o(y) \, \forall \:
 y\in \mathbb{R},}$ and $u_1, u_2, u_{2y} \neq \infty$ in any finite point,
 then ${u_1 - u_2 \neq 0}$ in any finite point. A priory estimations
 for derivatives in this case are equivalent
 to a priory estimations
 for $u_1, u_2, u_1 - u_2.$

 Author's approach is the next.
 We don't prove
 the existence of
 a priory estimations,
 and we don't search
 a priory estimations.
 We set
 a priory estimations
 for $u_1, u_2, u_1 - u_2$.
 Thinking in this direction was
 blockaded.

 Finally, refer to papers of
 Jia-Xing Hong
 \cite{Hong-93}, \cite{Hong-95}.
 The author found out on the existence of these papers
 in Beijing in 2002.
 The paper \cite{Hong-95}
 on global classical solutions of the equation
$$
 \mbox{hess} \, z = -k^2(x,y)
$$
 is unclear to the author,
 and
 the author couldn't find
 the paper
 \cite{Hong-93}.
 The paper
 \cite{Hong-93}
 could be considered, in particular, as solving the equation
$$
 \mbox{hess} \, z = -k^2(x,y)(1 + z_x^2 + z_y^2)^2.
$$


\section{Systems in Riemann invariants}

\subsection{One model equation}

 Consider the plane ${\mathbb{R}^2 = (x,y)}$.
 Let $u(x,y)$ be an unknown function.
 Let the equation
$$
    \left( \pt_x
           +\xi(x,y) \pt_y
    \right)u(x,y)=f(x,y)
    \eqno(2.1)
$$
 be an example for studying the main concepts,
 such as a characteristic and integration along a characteristic
 (\cite{RYa}, Chapter~1).
 Suppose
 $\xi ,f \in C^1(\mathbb{R}^2).$

{\bf Definition.}
 Consider the plane ${\mathbb{R}^2 = (x,y)}$.
 The curve
$$
  x=\tau, \qquad y=g(\tau ,x,y)
$$
 is a {\it characteristic} of the equation
 (2.1).
 Here the function $g(\tau ,x,y)$
 is a solution of the Cauchy problem
$$
 \left\{
   \begin{array}{l}
     \pt_\tau g=\xi (\tau ,g(\tau ,x,y)),
     \smallskip
     \\
     g(x,x,y)=y.
   \end{array}
 \right.
  \eqno(2.2)
$$

 {\bf Lemma~2.1.}
$$
    \left( \pt_x
           +\xi \pt_y
    \right)g=0.
    \eqno(2.3)
$$

 {\bf Proof.}
 Suppose the existence of functions
 $g_x,\: g_y$.
 Following
 \cite{Pozniak},
 we'll find these functions.

 Differentiating (2.2) with respect to $x$,
 we obtain an ordinary differential
 equation with respect to
 $\pt_x g$:
$$
  \pt_\tau \pt_x g
  =
  \frac{\ds \pt \xi}{\ds \pt g}
  (\tau ,g)
  \:
  \pt_x g.
 \eqno(2.4)
$$
 Differentiating the initial condition in (2.2)
 with respect to
 $x$,
 we obtain
$$
  \pt_x g(x,x,y)=0,
$$
or
$$
  \left.\pt_\tau g(\tau ,x,y)\right| _{\tau =x}+
  \left.\pt_x g(\tau ,x,y)\right| _{\tau =x}=0,
$$
therefore, taking into account (2.2), we have
$$
  \left.\pt_x g(\tau ,x,y)\right| _{\tau =x}=-\xi (x,y).
  \eqno(2.4^0)
$$
 Solving the Cauchy problem (2.4), $(2.4^0)$, we have
$$
  \pt_x g(\tau ,x,y)=-\xi (x,y)\exp
  \left\{
    \int _x^{\tau}
      \frac{\ds \pt \xi}{\ds \pt g}
      (t,g(t,x,y))
    \: dt
  \right\}.
  \eqno(2.5)
$$
 Analogously,
 solving the Cauchy problem
 with respect to the function
 $\pt_y g$
$$
\left\{
  \begin{array}{l}
    \pt_\tau \pt_y g
    =
      \frac{\ds \pt \xi}{\ds \pt g}(\tau ,g)
      \pt_y
      g,
      \\
    \pt_y g(\tau ,x,y)| _{\tau =x}=1,
  \end{array}
\right.
$$
we have
$$
  \pt_y g(\tau ,x,y)=\exp
  \left\{
    \int _x^{\tau}
      \frac{\ds \pt \xi}{\ds \pt g}
      (t,g(t,x,y))
    \: dt
  \right\}.
  \eqno(2.6)
$$
From (2.5), (2.6) we have (2.3). $\Box $

 Consider the Cauchy problem for the equation (2.1).
 Let an initial condition be
$$
  u(0,y)=u^0(y).
  \eqno(2.1^0)
$$
 Here
 $u^0\in C^1({\mathbb{R}}^1).$

 The problem
 (2.1), $(2.1^0)$
 is well-defined \cite{RYa}.
 Wi'll solve it in the half-plane
 ${x\ge 0}$.
 Solving in the half-plane
 ${x\le 0}$
 is analogous.

 {\bf Lemma~2.2.}
 The solution of the problem
 (2.1), $(2.1^0)$
 is
$$
  u(x,y)=u^0(g(0,x,y))+\int _0^x f(\tau ,g(\tau ,x,y))d\tau .
  \eqno(2.7)
$$

{\bf Proof.}
 The condition
 $(2.1^0)$
 is true, because
$$
  u(0,y)=u^o(g(0,0,y))=u^0(y)
$$
 by the initial condition (2.2).
 Differentiating (2.7),
 we use the formula
$$
  \pt_x
  \int _0^x
    \varphi(\tau ,x)
  d\tau =
  \int _0^x
    \pt_x\varphi(\tau ,x)
  d\tau
  +\varphi (x,x).
$$
 Therefore,
$$
    \left( \pt_x
           +\xi \pt_y
    \right)
    u
    =
    \left.
    u_t^0(t)
    \right|_{t=g(0,x,y)}
    \left( \pt_x
           +\xi \pt_y
    \right)
    g(0,x,y)+
$$
$$
  +\int _0^x
    \frac{\ds \pt f}{\ds \pt g}(\tau ,g(\tau ,x,y))
    \left( \pt_x
           +\xi \pt_y
    \right)
    g(\tau ,x,y)
  d\tau
  +f(x,g(x,x,y)).
$$
 By Lemma~2.1
 and the initial condition (2.2),
$$
  \left(
    \pt_x
    +\xi \pt_y
  \right)
  u(x,y)=f(x,y). \quad \Box
$$

 The procedure of solving
   (2.1), $(2.1^0)$
 by
   (2.7)
 is called {\it integrating the equation} (2.1)
 {\it along the characteristic.}


\subsection{Systems in Riemann invariants}

 The developed theory works in more general case.
 Consider vector-functions
$$
\begin{array}{lll}
 u=(u_1,\dots ,u_m),\quad & u_i=u_i(x,y)\in C^1({\mathbb{R}}^2),&
 \smallskip
 \\
 \xi =(\xi _1,\dots ,\xi _m),\quad &\xi _i=\xi_i(x,y,u)\in C^1({\mathbb{R}}^2),&
 \smallskip
 \\
 f=(f_1,\dots ,f_m),\quad & f_i=f_i(x,y,u)\in C^1({\mathbb{R}}^2),&
 i=1,\dots ,m.
\end{array}
$$
{\bf Definition.}
 The system
$$
    \left( \pt_x
           +\xi_i(x,y,u)
           \pt_y
    \right)
    u_i(x,y)=f_i(x,y,u), \quad i=1,\dots ,m,
    \eqno(2.8)
$$
 is called a
 {\it system in Riemann invariants} \cite{RYa}.

 Consider the Cauchy problem for the system (2.8).
 Let an initial condition be
$$
  u(0,y)=u^o(y).
  \eqno(2.8^0)
$$
 Here
 $u^0=(u_1^0,\dots ,u_m^0),\quad
     u_i^0=u_i^0(y)\in C^1(\mathbb{R}),\quad i=1,\dots ,m.$\\
{\bf Definition.}
 Consider the plane ${\mathbb{R}^2 = (x,y)}$.
 The curve
$$
  x=\tau, \qquad y=g_i(\tau ,x,y)
$$
 is an {\it i-th characteristic} of the system
 (2.8).
 Here the function $g_i(\tau ,x,y)$
 is a solution of the Cauchy problem
$$
\left\{
\begin{array}{lr}
    \pt_\tau g_i
    =
    \xi _i
     (\tau,
      g_i(\tau ,x,y),
      u(\tau ,x,y)
     ),
  \smallskip
  \\
    g_i(x,x,y)=y.
\end{array}
\right.
 \eqno (2.9)
$$
 {\bf Lemma~2.3.}
 The $i\mbox{-th}$
 component of the solution of the problem
 (2.8), $(2.8^0)$
 is a result of integrating the
 $i\mbox{-th}$ equation of the system
 along the
 $i\mbox{-th}$
 characteristic,
 ${i=1,\dots ,m:}$
$$
  u_i(x,y)=u_i^0(g_i(0,x,y))+
  {\ds \int _0^x}
    f_i(\tau ,g_i(\tau ,x,y),
      u(\tau ,g_i(\tau ,x,y)))
   \: d\tau .
 \eqno(2.10)
$$
{\bf Proof.}
 Suppose the solution
 $u(x,y)$
 of the problem
 (2.8), $(2.8^0)$
 is known.
 Inserting it in
 (2.8),
 we see that each equation of the system is
 in the form of
 (2.1),
 because it is possible to suppose
$$
  \xi _i(x,y,u(x,y))=\tilde \xi (x,y),\quad
  f_i(x,y,u(x,y))=\tilde f(x,y),\quad i=1,\dots ,m.
$$
 Therefore for
$$
    (
      \pt_x
           +\tilde \xi _i (x,y) \pt_y
    )
    \: u_i(x,y)=\tilde f_i(x,y),\quad i=1,\dots ,m,
$$
 we obtain Lemma 2.1 and Lemma 2.2.
 Further,
 (2.5), (2.6)
 are
$$
  \!\!\!
{\ds
\begin{array}{l}
  \!\!\!
  \pt_x
  g_i(\tau ,x,y)
  \!
  =
  \!\!
  -\xi _i (x,y,u)
  \!
  \exp
  \!
  \left\{ \!
    {\ds\int _x^{\tau}}
    \!\!
    \left[
        \xi _{iy} \! + \!\!\!
        {\ds\sum}_{j=1}^m
        \frac{\ds \pt \xi_i}{\ds \pt u_j}\frac{\ds \pt u_j}{\ds \pt y}
    \right]
    \!\!
    (t,g_i(t,x,y))dt \!
 \right\}
 \!
 ,
 \\
 \\
  \!\!\!\pt_y g_i(\tau ,x,y)=
  \exp
  \left\{
    {\ds\int _x^{\tau}}
    \left[
        \xi _{iy}+
        {\ds\sum}_{j=1}^m
            \frac{\ds \pt \xi_i}{\ds \pt u_j}
            \frac{\ds \pt u_j}{\ds \pt y}
    \right]
    (t,g_i(t,x,y))dt
 \right\}
 \!
 .
\end{array}
}
  \!\!\!\!\!
 \eqno(2.11)
$$
 Finally, (2.7) is
 (2.10). $\Box$

 The formula (2.10) is used for studying properties of solutions.
 It is recursive, so it doesn't give a solution of the problem
 (2.8), $(2.8^0)$
 in an explicit form.
 The solution could be obtained by successive approximations.



\section{Systems in Riemann invariants
         and hyperbolic Monge--Ampere equations}

 Here we follow a chapter of the dissertation
 of Tunitsky \cite{Tunitsky-dis}.
 This chapter was published separately as \cite{Tunitsky-dep}.

 It was well-known that the hyperbolic
 Monge--Ampere equation
 could be reduced to a system of five equations
 of the first order
 (\cite{Cour}, Supplement~1 to Chapter~V, \S~2),
 and it was well-known that this system
 could be reduced to a system in Riemann invariants
 (in \cite{Cour} the system of five equations
 isn't a system in Riemann invariants).
 In Tunitsky's paper the system in invariants
 is written in an explicit form,
 and everyone can use it.


\subsection{The Cauchy problem for the Monge--Ampere equation}

 Consider the plane ${\mathbb{R}^2 = (x,y)}$ \
 and an unknown function
 $z=z(x,y)$
 on this plane.
 Consider
 the Monge--Ampere equation with respect to $z$
   $$
   A+Bz_{xx}+Cz_{xy}+Dz_{yy}+E\,\mbox{\rm hess}\, z=0.
   $$
 Here
 $\mbox{\rm hess}\, z=z_{xx}z_{yy}-z^2_{xy},\quad A,B,C,D,E$
 are functions of
 $x,y,z$,\ $z_x,z_y$; \;\; ${E\ne 0}$.
 Dividing by
 $E$,
 we obtain
 ${E\equiv 1}$,
 so we consider the equation
$$
   A+Bz_{xx}+Cz_{xy}+Dz_{yy}+\mbox{\rm hess}\, z = 0. \eqno (3.1)
$$
 Suppose
 ${A,B,C,D \in C^2(\mathbb{R}^5)}$.

 Let ${z(x,y)}$ be a
 $C^2\mbox{\rm -solution}$
 of the equation (3.1).
 We say that
 the equation (3.1)
 is hyperbolic at the solution
 ${z(x,y)}$
 (\cite{Cour}, Supplement~1 to Chapter~5, ${{\S}~2}$),
 if
 $$
 \Delta^2(x,y,z(x,y),z_x(x,y),z_y(x,y))
 =
 C^2 - 4BD + 4AE > 0.
 \eqno(3.2)
 $$
 Recall that ${E\equiv 1}$.
 By ${\Delta}^2 > 0$, we set
 ${\Delta > 0}$.
 Initial functions for
 the Cauchy problem on the $Oy$ axis are
$$
 z(0,y) = z_0(y), \qquad z_x(0,y)=p_0(y);
 \eqno(3.3)
$$
 here ${z_0 \in C^3(\mathbb{R}^1)}$,
 and
 ${p_0 \in C^2(\mathbb{R}^1)}$.
 Suppose
 $z_0$, $p_0$
 satisfy
 the next two conditions.
 First,
 the axis $Oy$ is free, i.e.
$$
 z_0^{''}(y)
 +
 B
 (0,y,z_0(y),p_0(y),z_0^{'}(y))
 \neq
 0.
 \eqno(3.4)
$$
 Secondly,
 on the axis $Oy$
 the hyperbolic condition (3.2)
$$
 \Delta^2
 (0,y,z_0(y),p_0(y),z_0^{'}(y))
 >
 0
 \eqno(3.5)
$$
 holds.

 The full formulation for the Cauchy problem
 for the Monge--Ampere equation is the next.
 The aim is to find the
 $C^3\mbox{\rm -function}$
 ${z(x,y)}$
 satisfying
 the initial condition (3.3),
 the equation (3.1),
 and (at the solution ${z(x,y)}$)
 the hyperbolic condition (3.2).

 Now we explain
 ${C^3\mbox{-smoothness}}$,
 when
 a classical solution of the equation (3.1)
 is assumed to be
 ${C^2\mbox{-smooth}}$.
 To obtain the system in Riemann invariants,
 the integrability conditions
 ${p_{xy} = p_{yx}}$, ${q_{xy} = q_{yx}}$
 will be used.
 Here $p,\ q$ are assumed to be $z_x,\ z_y$.

\subsection{The system in total differentials}

 Let ${z(x,y)}$ be a $C^3\mbox{-solution}$ of the equation (3.1)
 in some domain $T$. Suppose the equation (3.1)
 is hyperbolic at the solution $z$, and
$$
 z_{yy}(x,y)
 +
 B
 (x,y,z(x,y),z_x(x,y),z_y(x,y))
 \neq
 0,
 \eqno(3.6)
$$
 ${\forall \, (x,y) \in T}$. The inequality (3.6)
 means that vertical lines ${x=\mbox{const}}$ are free.

 By definition, put
$$
 u_1
 =
 \frac
 {\ds C + \Delta - 2z_{xy}}
 {2(z_{yy} + B)}
 ,
 \qquad
 u_2
 =
 \frac
 {\ds C - \Delta - 2z_{xy}}
 {2(z_{yy} + B)}
 .
 \eqno(3.7)
$$
 Functions $u_1$, $u_2$ are tangents of angles of inclinations
 of characteristics of the equation (3.1)
 (\cite{Cour}, Supplement~1 to Chapter~V, \S~2).
 By the hyperbolic condition
 (3.2), we obtain ${u_1 \neq u_2}$,
 so it is possible to solve
 (3.7) (uniquely)
 with respect to $z_{xy}$,
 $z_{yy}$:
$$
 z_{yy}
 =
 \frac
 {\ds \Delta}
 {u_1 - u_2}
 -
 B
 ,
 \qquad
 z_{xy}
 =
 \frac{\ds \Delta}{\ds 2}
 \frac
 {u_1 + u_2}
 {u_2 - u_1}
 +
 \frac{\ds C}{\ds 2}
 .
 \eqno(3.8)
$$
 Substututing (3.8) in (3.1),
 we obtain a linear equation with respect to $z_{xx}$.
 Solving it, we have
$$
 z_{xx}
 =
 \Delta
 \frac
 {u_1 u_2}
 {u_1 - u_2}
 -
 D
 .
 \eqno(3.9)
$$
 Following Monge, by definition, put
$$
 z_x = p, \qquad z_y = q.
 \eqno(3.10)
$$
 From (3.8), (3.9),
 we obtain that
 $z$, $p$, $q$
 as functions of variables
 $x$, $y$ in the domain $T$
 are a solution of the system in total differentials
$$
\begin{array}{l}
 p_x
 =
 \Delta
 (x,y,z,p,q)
 \frac
 {\ds u_1 u_2}
 {\ds u_1 - u_2}
 (x,y)
 -
 D
 (x,y,z,p,q)
 ,
 \medskip
 \\
 p_y
 =
 \frac{\ds \Delta}{\ds 2}
 (x,y,z,p,q)
 \frac
 {\ds u_1 + u_2}
 {\ds u_2 - u_1}
 (x,y)
 +
 \frac{\ds C}{\ds 2}
 (x,y,z,p,q)
 ,
 \medskip
 \\
 q_x
 =
 \frac{\ds \Delta}{\ds 2}
 (x,y,z,p,q)
 \frac
 {\ds u_1 + u_2}
 {\ds u_2 - u_1}
 (x,y)
 +
 \frac{\ds C}{\ds 2}
 (x,y,z,p,q)
 ,
 \medskip
 \\
 q_y
 =
 \frac
 {\ds \Delta (x,y,z,p,q)}
 {\ds (u_1 - u_2)(x,y)}
 -
 B
 (x,y,z,p,q)
 .
\end{array}
 \eqno(3.11)
$$
 Right sides of
 (3.11)
 are continuously differentiable functions of
 $x,y,z,p,q$.

 Thus, we proved the next.
 If
 ${z(x,y)}$ is a $C^3\mbox{-solution}$ of the equation (3.1)
 in the domain $T$, and inequalities (3.2) и (3.6) hold,
 then three functions
 $z,p,q$
 are defined, and they are a solution of the system
 (3.10)--(3.11).

 Conversely,
 suppose there exist two
 continuously differentiable functions
 ${(u_1,\ u_2)}$
 in the domain $T$,
 such as
 ${u_1 \neq u_2}$,
 and suppose there exists a
 $C^1\mbox{-solution}$
 of the system (3.10)--(3.11)
 in the domain $T$,
 such as
$$
 \Delta
 (x,y,z(x,y),z_x(x,y),z_y(x,y))
 >
 0
 .
 \eqno(3.12)
$$
 Clearly, in this case
 we have
 ${z_{xx} = p_x}$, ${z_{xy} = p_y = q_x}$, ${z_{yy} = q_y}$.
 Therefore, ${z \in C^3(T)}$.
 Substituting second derivatives of $z$ in the left part of (3.1),
 we see that ${z(x,y)}$ is a solution of the equation (3.1).
 By (3.12), the equation (3.1) is hyperbolic at $z$, and
 the inequality (3.6) holds.

 Summarize the obtained results as a lemma.

 {\bf Lemma~3.1.}
 There exists in the domain $T$ a function ${z(x,y) \in C^3(T)}$
 satisfying (3.1), (3.2), (3.6)
 iff
 there exists in the domain $T$
 a $C^1\mbox{-solution}$ ${(z,p,q)}$
 of the system (3.10)--(3.11),
 satisfying (3.12). $\Box$

\subsection{The system in Riemann invariants}

 The system of differential equations in total differentials is
 overdetermined, so, generally speaking, it doesn't have a
 solution. The solution exists iff integrability conditions
 hold. Integrability conditions for the equations (3.10)--(3.11)
 could be reduced to differential relations for functions
 $u_1, u_2, z, p, q$.

 Indeed, right parts of the system (3.10)--(3.11)
 are continuously differentiable, so
 functions $z,\ p,\ q$ are twice continuously differentiable.
 Therefore, on $T$ we have equalities
$$
 z_{xy} = z_{yx},
 \qquad
 p_{xy} = p_{yx},
 \qquad
 q_{xy} = q_{yx}.
 \eqno(3.13)
$$
 By (3.10)--(3.11), the first equality (3.13) is an identity.
 The second equality (3.13) is
$$
\begin{array}{ll}
 &
 \Delta_y
 \frac
 {\ds u_1 u_2}
 {\ds u_1 - u_2}
 +
 \Delta
 \frac
 {\ds u_1^2 u_{2y} - u_2^2 u_{1y}}
 {\ds (u_1 - u_2)^2}
 -
 D_y
 +
 \left(
    \frac
    {\ds u_1 u_2}
    {\ds u_1 - u_2}
    \Delta_z
    -
    D_z
 \right)
 q
 +
 \medskip
 \\
 &
 +
 \left(
    \frac
    {\ds u_1 u_2}
    {\ds u_1 - u_2}
    \Delta_p
    -
    D_p
 \right)
 \left(
    \frac{\ds \Delta}{\ds 2}
    \frac
    {\ds u_1 + u_2}
    {\ds u_2 - u_1}
    +
    \frac{\ds C}{\ds 2}
 \right)
 +
 \medskip
 \\
 &
 +
 \left(
    \frac
    {\ds u_1 u_2}
    {\ds u_1 - u_2}
    \Delta_q
    -
    D_q
 \right)
 \left(
    \frac
    {\ds \Delta}
    {\ds u_1 - u_2}
    -
    B
 \right)
 =
 \medskip
 \\
 =
 &
 \frac{\ds \Delta_x}{\ds 2}
 \frac
 {\ds u_1 + u_2}
 {\ds u_2 - u_1}
 +
 \Delta
 \frac
 {\ds u_2 u_{1x} - u_1 u_{2x}}
 {\ds (u_2 - u_1)^2}
 +
 \frac{\ds C_x}{\ds 2}
 +
 \left(
    \frac
    {\ds u_1 + u_2}
    {\ds u_2 - u_1}
    \frac{\ds \Delta_z}{\ds 2}
    +
    \frac{\ds C_z}{\ds 2}
 \right)
 p
 +
 \medskip
 \\
 &
 +
 \left(
    \frac
    {\ds u_1 + u_2}
    {\ds u_2 - u_1}
    \frac{\ds \Delta_p}{\ds 2}
    +
    \frac{\ds C_p}{\ds 2}
 \right)
 \left(
    \Delta
    \frac
    {\ds u_1 u_2}
    {\ds u_1 - u_2}
    -
    D
 \right)
 +
 \medskip
 \\
 &
 +
 \left(
    \frac
    {\ds u_1 + u_2}
    {\ds u_2 - u_1}
    \frac{\ds \Delta_q}{\ds 2}
    +
    \frac{\ds C_q}{\ds 2}
 \right)
 \left(
    \frac{\ds \Delta}{\ds 2}
    \frac
    {\ds u_1 + u_2}
    {\ds u_2 - u_1}
    +
    \frac{\ds C}{\ds 2}
 \right)
 .
\end{array}
 \eqno(3.14)
$$
 Analogously, the third equality (3.13) is
$$
\begin{array}{ll}
 &
 \frac{\ds \Delta_y}{\ds 2}
 \frac
 {\ds u_1 + u_2}
 {\ds u_2 - u_1}
 +
 \Delta
 \frac
 {\ds u_2 u_{1y} - u_1 u_{2y}}
 {\ds (u_2 - u_1)^2}
 +
 \frac{\ds C_y}{\ds 2}
 +
 \left(
    \frac
    {\ds u_1 + u_2}
    {\ds u_2 - u_1}
    \frac{\ds \Delta_z}{\ds 2}
    +
    \frac{\ds C_z}{\ds 2}
 \right)
 q
 +
 \medskip
 \\
 &
 +
 \left(
    \frac
    {\ds u_1 + u_2}
    {\ds u_2 - u_1}
    \frac{\ds \Delta_p}{\ds 2}
    +
    \frac{\ds C_p}{\ds 2}
 \right)
 \left(
    \frac{\ds \Delta}{\ds 2}
    \frac
    {\ds u_1 + u_2}
    {\ds u_2 - u_1}
    +
    \frac{\ds C}{\ds 2}
 \right)
 +
 \medskip
 \\
 &
 +
 \left(
    \frac
    {\ds u_1 + u_2}
    {\ds u_2 - u_1}
    \frac{\ds \Delta_q}{\ds 2}
    +
    \frac{\ds C_q}{\ds 2}
 \right)
 \left(
    \frac
    {\ds \Delta}
    {\ds u_1 - u_2}
    -
    B
 \right)
 =
 \medskip
 \\
 =
 &
 \frac
 {\ds \Delta_x}
 {\ds u_1 - u_2}
 +
 \Delta
 \frac
 {\ds u_{2x} - u_{1x}}
 {\ds (u_1 - u_2)^2}
 -
 B_x
 +
 \left(
    \frac
    {\ds \Delta_z}
    {\ds u_1 - u_2}
    -
    B_z
 \right)
 p
 +
 \medskip
 \\
 &
 +
 \left(
    \frac
    {\ds \Delta_p}
    {\ds u_1 - u_2}
    -
    B_p
 \right)
 \left(
    \Delta
    \frac
    {\ds u_1 u_2}
    {\ds u_1 - u_2}
    -
    D
 \right)
 +
 \medskip
 \\
 &
 +
 \left(
    \frac
    {\ds \Delta_q}
    {\ds u_1 - u_2}
    -
    B_q
 \right)
 \left(
    \frac{\ds \Delta}{\ds 2}
    \frac
    {\ds u_1 + u_2}
    {\ds u_2 - u_1}
    +
    \frac{\ds C}{\ds 2}
 \right)
 .
\end{array}
 \eqno(3.15)
$$

 (3.14)--(3.15) is a linear system of two algebraic
 equations with respect to
 ${u_{1x} + u_2 u_{1y}}$,
 ${u_{2x} + u_1 u_{2y}}$.
 The determinant is
 ${\Delta^2/(u_1 - u_2)^3}$.
 By ${u_1 \neq u_2}$,
 this system has a unique solution.
 The solution is
$$
\left\{
\begin{array}{l}
    u_{1x} + u_2 u_{1y}
    =
    E_o + E_1 u_1 + E_2 u_2 + E_3 u_1^2 +
    E_4 u_1 u_2 + E_5 u_1^2 u_2
    ,
 \medskip
 \\
    u_{2x} + u_1 u_{2y}
    =
    I_o + I_1 u_1 + I_2 u_2 + I_3 u_2^2 +
    I_4 u_1 u_2 + I_5 u_1 u_2^2
    .
\end{array}
\right.
 \eqno(3.16)
$$
 Coefficients $E_j,\ I_j,\ {0 \le j \le 5}$,
 depend on functions $B,C,D,\Delta$
 and on their first derivatives:
$$
\begin{array}{l}
 E_0 = I_0 = D_q,
 \qquad
 E_5 = I_5 = -B_p,
 \medskip
 \\
 E_1
 =
 \phantom{-}
 \alpha_1 + \alpha_2
 +
 \frac{\ds 1}{\ds 4\Delta}
 (
    \phantom{-}
    CC_q - 3\Delta C_q + C\Delta_q + \Delta\Delta_q - 2\Delta D_p
 )
 ,
 \medskip
 \\
 E_2
 =
 -\alpha_1 - \alpha_2
 +
 \frac{\ds 1}{\ds 4\Delta}
 (
   -CC_q - \phantom{3}\Delta C_q - C\Delta_q - \Delta\Delta_q - 2\Delta D_p
 )
 ,
 \medskip
 \\
 E_3
 =
 -\beta_1 + \beta_2
 +
 \frac{\ds 1}{\ds 4\Delta}
 (
   -CC_p + \phantom{3}\Delta C_p + C\Delta_p - \Delta\Delta_p + 2\Delta B_q
 )
 ,
 \medskip
 \\
 E_4
 =
 \phantom{-}
 \beta_1 - \beta_2
 +
 \frac{\ds 1}{\ds 4\Delta}
 (
   \phantom{-}
   CC_p + 3\Delta C_p - C\Delta_p + \Delta\Delta_p + 2\Delta B_q
 )
 ,
 \medskip
 \\
 I_1
 =
 \phantom{-}
 \alpha_1 - \alpha_2
 +
 \frac{\ds 1}{\ds 4\Delta}
 (
    \phantom{-}
    CC_q - \phantom{3}\Delta C_q - C\Delta_q + \Delta\Delta_q - 2\Delta D_p
 )
 ,
 \medskip
 \\
 I_2
 =
 -\alpha_1 + \alpha_2
 +
 \frac{\ds 1}{\ds 4\Delta}
 (
   -CC_q - 3\Delta C_q + C\Delta_q - \Delta\Delta_q - 2\Delta D_p
 )
 ,
 \medskip
 \\
 I_3
 =
 \phantom{-}
 \beta_1 + \beta_2
 +
 \frac{\ds 1}{\ds 4\Delta}
 (
   \phantom{-}
   CC_p + \phantom{3}\Delta C_p + C\Delta_p + \Delta\Delta_p + 2\Delta B_q
 )
 ,
 \medskip
 \\
 I_4
 =
 -\beta_1 - \beta_2
 +
 \frac{\ds 1}{\ds 4\Delta}
 (
   -CC_p + 3\Delta C_p - C\Delta_p - \Delta\Delta_p + 2\Delta B_q
 )
 ,
\end{array}
 \eqno(3.17)
$$
 where
$$
\begin{array}{l}
 \alpha_1
 =
 \frac{\ds 1}{\ds 2\Delta}
 (
    2D_y + C_x + 2D_z q + C_z p + C D_p - D C_p - 2B D_q
 )
 ,
 \medskip
 \\
 \alpha_2
 =
 \frac{\ds 1}{\ds 2\Delta}
 (
    \Delta_x + \Delta_z p - D \Delta_p
 )
 ,
 \medskip
 \\
 \beta_1
 =
 \frac{\ds 1}{\ds 2\Delta}
 (
    2B_x + C_y + 2B_z p + C_z q + C B_q - B C_q - 2D B_p
 )
 ,
 \medskip
 \\
 \beta_2
 =
 \frac{\ds 1}{\ds 2\Delta}
 (
    \Delta_y + \Delta_z q - B \Delta_q
 )
 .
\end{array}
$$

 Let ${z(x,y)}$ be a $C^3\mbox{-solution}$ of (3.1)
 in $T$.
 Suppose inequalities (3.2), (3.6) hold.
 By Lemma~3.1, so that $z,\ p=z_x,\ q=z_y$ satisfy (3.10)--(3.11),
 where $u_1$, $u_2$ are defined by (3.7).

 Multiplying the second equation (3.10) by $u_1$, and then
 adding with the first equation (3.10), we get
$$
 z_x + u_1 z_y = p + u_1 q.
 \eqno(3.18)
$$
 Analogously, consider (3.11).
 Multiplying the second equation (3.11) by $u_2$, and then
 adding with the first equation (3.11), we get
$$
 p_x + u_2 p_y = \frac{\ds C + \Delta}{\ds 2} u_2 - D.
 \eqno(3.19)
$$
 Multiplying the fourth equation (3.11) by $u_1$, and then
 adding with the third equation (3.11), we get
$$
 q_x + u_1 q_y = \frac{\ds C + \Delta}{\ds 2} - B u_1.
 \eqno(3.20)
$$

 Thus, the next fact is proved.
 \\
 \textbf{Theorem~3.1.}
 Let ${z(x,y)}$ be a $C^3\mbox{-solution}$ of (3.1)
 in the domain $T$.
 Suppose
 (3.2), (3.6) are satisfied by $z$. Then the set of functions
 ${(u_1,\ u_2,\ z,\ p,\ q)}$, where $u_1,\ u_2$ we get from
 (3.7), ${p=z_x}$, ${q=z_y}$, is a $C^1\mbox{-solution}$
 of the system of five equations (3.16)--(3.20) in $T$. $\Box$

\subsection{Reducing the Monge--Ampere equation
 to a system in Riemann invariants}

 Suppose the domain $T$ has a nonempty intersection with the axis
 $Oy$.
 The next statement in some sense is inverse with respect to
 Theorem~3.1.

 \textbf{Theorem~3.2.}
 Let
 $(u_1,\ u_2,\ z,\ p,\ q)$
 be a $C^1\mbox{-solution}$ of (3.16)--(3.20)
 in the domain $T$.
 Let initial values for this solution be
$$
\begin{array}{l}
 z(0,y) = z_0(y),
 \medskip
 \\
 p(0,y) = p_0(y),
 \medskip
 \\
 q(0,y) = z_0^{'}(y),
 \medskip
 \\
 u_1(0,y)
 =
 {\ds
    \frac
    {(C + \Delta) ((0,y,z_0(y),p_0(y),z_0^{'}(y))) - 2p_0^{'}(y)}
    {2(z_0^{''} + B(0,y,z_0(y),p_0(y),z_0^{'}(y)))}
 }
 ,
 \medskip
 \\
 u_2(0,y)
 =
 {\ds
    \frac
    {(C - \Delta) ((0,y,z_0(y),p_0(y),z_0^{'}(y))) - 2p_0^{'}(y)}
    {2(z_0^{''} + B(0,y,z_0(y),p_0(y),z_0^{'}(y)))}
 }
 .
\end{array}
 \eqno(3.21)
$$
 Let $T$ be a domain in which this solution is defined.
 Suppose (3.12) holds.
 Then
 $z$ is a $C^3\mbox{-solution}$ of (3.1)--(3.3) in the domain $T$,
 and, further,
 ${z_x = p},\ {z_y = q},$ and (3.6) holds.
 \\
 \textbf{Proof.}
 First let us prove that
 ${u_1 \neq u_2}$ in the domain $T$.
 By (3.21) and (3.5), we get
$$
 u_1(0,y) - u_2(0,y)
 =
 {\ds
    \frac
    {\Delta
        (0,y,z_0(y),p_0(y),z_0^{'}(y))
    }
    {z_0^{''} (y) + B
        (0,y,z_0(y),p_0(y),z_0^{'}(y))
    }
 }
 \neq 0
 .
 \eqno(3.22)
$$
 Consider (3.16). Subtracting the second equation
 from the first, and taking into account (3.17), we obtain
$$
\begin{array}{l}
  (u_1 - u_2)_x + u_2 (u_1 - u_2)_y
  =
 \medskip
 \\
 \phantom{12345}
  =
  (u_1 - u_2)
  \left(
    2\alpha_2 +
    {\ds
      \frac
      {C \Delta_q - \Delta C_q}
      {2\Delta}
    }
   \right.
    +
 \medskip
 \\
 \phantom{12345}
    +
    (u_1 + u_2)
    \left(
      \beta_2
      +
      \frac{\ds B_q}{\ds 2}
      +
      {\ds
        \frac
        {C \Delta_p + \Delta C_p}
        {4 \Delta}
      }
    \right)
      +
 \medskip
 \\
 \phantom{12345}
    +
  \left.
      (u_2 - u_1)
      \left(
        \beta_1
        +
        {\ds
          \frac
          {C C_p + \Delta \Delta_p}
          {4 \Delta}
        }
      \right)
      -
      B_p u_1 u_2 + u_{2y}
  \right)
  .
\end{array}
 \eqno(3.23)
$$
 A solution of (3.16)--(3.20)
 is defined in
 the domain $T$,
 so
 on the segment ${[0,x]}$
 there exist a solution of the Cauchy problem
$$
\left\{
\begin{array}{l}
 {\ds
   \frac{d g}{d\tau}
 }
   =
   u_2(\tau, g)
   ,
 \medskip
 \\
   g(x,x,y) = y
   .
\end{array}
\right.
$$
 If ${(x,y) \in T}$, then
 ${(\tau,\ g(\tau,x,y)) \in T}$
 for ${0 \le \tau \le x}$.
 Integrating the equation (3.23)
 along the characteristic
 ${\eta = g(\tau,x,y)}$ over $[0,x]$, we get
$$
\begin{array}{l}
  (u_1 - u_2)(x,y)
  =
 \medskip
 \\
 \phantom{12345}
  =
  (u_1 - u_2)
  (0, g(0,x,y))
  \times
 \medskip
 \\
 \phantom{12345}
  \times
  \exp
  \{
  {\ds \int_0^x}
  \left(
    2\alpha_2 +
    {\ds
      \frac
      {C \Delta_q - \Delta C_q}
      {2\Delta}
    }
   \right.
    +
 \medskip
 \\
 \phantom{12345}
    +
    (u_1 + u_2)
    \left(
      \beta_2
      +
      \frac{\ds B_q}{\ds 2}
      +
      {\ds
        \frac
        {C \Delta_p + \Delta C_p}
        {4 \Delta}
      }
    \right)
      +
 \medskip
 \\
 \phantom{12345}
    +
  \left.
      (u_2 - u_1)
      \left(
        \beta_1
        +
        {\ds
          \frac
          {C C_p + \Delta \Delta_p}
          {4 \Delta}
        }
      \right)
      -
      B_p u_1 u_2 + u_{2y}
  \right)
 \medskip
 \\
 \phantom{12345}
  (\tau, g(\tau, x,y))
  \: d\tau
  \}
  .
\end{array}
 \eqno(3.24)
$$
 From (3.24), taking into account (3.22), we get
 ${u_1 \neq u_2}$ in $T$.

 Further, by definition, put
$$
\begin{array}{l}
  r(x,y,z,p,q)
  =
  \Delta
  {\ds
    \frac
    {u_1 u_2}
    {u_1 - u_2}
  }
  -
  D
  ,
 \medskip
 \\
  s(x,y,z,p,q)
  =
  {\ds
    \frac
    {\Delta}
    {2}
  }
  {\ds
    \frac
    {u_1 + u_2}
    {u_2 - u_1}
  }
  +
  {\ds
    \frac
    {C}
    {2}
  }
  ,
 \medskip
 \\
  t(x,y,z,p,q)
  =
  {\ds
    \frac
    {\Delta}
    {u_1 - u_2}
  }
  -
  B
  .
\end{array}
 \eqno(3.25)
$$
 Therefore, equations (3.19)--(3.20)
 are
$$
  p_x + u_2 p_y = r + u_2 s
  ,
  \eqno(3.19')
$$
$$
  q_x + u_1 q_y = s + u_1 t
  .
  \eqno(3.20')
$$
 Taking into account ${u_1 \neq u_2}$ in $T$, we obtain equivalence
 of (3.16)
 to (3.14)--(3.15).
 Equations (3.14)--(3.15)
 are
$$
  r_y + r_z q + r_p s + r_q t = s_x + s_z p + s_p r + s_q s
  ,
  \eqno(3.14')
$$
$$
  s_y + s_z q + s_p s + s_q t = t_x + t_z p + t_p r + t_q s
  .
  \eqno(3.15')
$$
 We obtain characteristics of the system
 (3.18), $(3.19')$, $(3.20')$
 from the Cauchy problem
$$
 \left\{
 \begin{array}{l}
   {\ds
    \frac{d g_i}{d \tau}
   }
   =
   u_{3-i} (\tau, g_i)
   ,
 \medskip
 \\
 g_i (x,x,y) = y \qquad (i = 1,2).
 \end{array}
 \right.
 \eqno(3.26)
$$
 Integrating
 (3.18), $(3.19')$, $(3.20')$
 along relevant characteristics
 over $[0,x]$, we get
$$
 \begin{array}{l}
   z(x,y)
   =
   z_0(g_2(0,x,y))
   +
   {\ds \int_0^x}
   \{
     p + u_1 q
   \}
   (\tau, g_2(\tau, x,y))
   \: d\tau
   ,
 \medskip
 \\
   p(x,y)
   =
   p_0(g_1(0;x,y))
   +
   {\ds \int_0^x}
   \{
     r + u_2 s
   \}
   (\tau, g_1(\tau, x,y))
   \: d\tau
   ,
 \medskip
 \\
   q(x,y)
   =
   z_0^{'}(g_2(0,x,y))
   +
   {\ds \int_0^x}
   \{
     s + u_1 t
   \}
   (\tau, g_2(\tau, x,y))
   \: d\tau
   .
 \end{array}
 \eqno(3.27)
$$
 Right parts of
  (3.26)
 are continuously differentiable.
 Therefore, functions
 ${g_i(\tau, x,y)}\ {(i=1,2)}$
 are continuously differentiable,
 and they have continuous secondary derivatives
 with respect to
 $\tau,\ x$ and
 with respect to
 $\tau,\ y$.
 Also for derivatives of
 $g_i$
 with respect to $x$
 and with respect to
 $y$
 we have
$$
   g_{ix}
 (\tau,x,y)
 +
 u_{3-i}
 (x,y)
   g_{iy}
 (\tau,x,y)
 =
 0
 .
 \eqno(3.28)
$$

 Now we find first derivatives of functions
 $z,\ p$ и $q$.
 Differentiating
 (3.27), taking into account
 $(3.14')-(3.15')$, (3.28),
 integrating by parts,
 taking into account
 the initial data
 (3.21),
 we obtain
$$
 \begin{array}{l}
   z_y(x,y)
   =
   q(x,y)
   +
   {\ds \int_0^x}
   \{
     (p_y - s) + (s - q_x)
   \}
   (\tau, g_2(\tau, x,y))
    g_{2y}
    (\tau,x,y)
   \: d\tau
   ,
 \medskip
 \\
   z_x(x,y)
   =
   p(x,y)
   +
   {\ds \int_0^x}
   \{
     (p_y - s) + (s - q_x)
   \}
   (\tau, g_2(\tau, x,y))
    g_{2x}
    (\tau,x,y)
   \: d\tau
   ,
 \medskip
 \\
   p_y(x,y)
   =
   s(x,y,z(x,y),p(x,y),q(x,y))
   +
 \medskip
 \\
   \phantom{123}
   +
   {\ds \int_0^x}
   \{
     s_z(p - z_x) + r_z(z_y - q) + s_p(r - p_x) + r_p(p_y - s)
     +
 \medskip
 \\
   \phantom{123}
   +
     s_q(s - q_x) + r_q(q_y - t)
   \}
   (\tau, g_1(\tau, x,y))
     g_{1y}
    (\tau,x,y)
   \: d\tau
   ,
 \medskip
 \\
   p_x(x,y)
   =
   r(x,y,z(x,y),p(x,y),q(x,y))
   +
 \medskip
 \\
   \phantom{123}
   +
   {\ds \int_0^x}
   \{
     s_z(p - z_x) + r_z(z_y - q) + s_p(r - p_x) + r_p(p_y - s)
     +
 \medskip
 \\
   \phantom{123}
   +
     s_q(s - q_x) + r_q(q_y - t)
   \}
   (\tau, g_1(\tau, x,y))
     g_{1x}
    (\tau,x,y)
   \: d\tau
   ,
 \medskip
 \\
   q_y(x,y)
   =
   t(x,y,z(x,y),p(x,y),q(x,y))
   +
 \medskip
 \\
   \phantom{123}
   +
   {\ds \int_0^x}
   \{
     t_z(p - z_x) + s_z(z_y - q) + t_p(r - p_x) + s_p(p_y - s)
     +
 \medskip
 \\
   \phantom{123}
   +
     t_q(s - q_x) + s_q(q_y - t)
   \}
   (\tau, g_2(\tau, x,y))
    g_{2y}
    (\tau,x,y)
   \: d\tau
   ,
 \medskip
 \\
   q_x(x,y)
   =
   s(x,y,z(x,y),p(x,y),q(x,y))
   +
 \medskip
 \\
   \phantom{123}
   +
   {\ds \int_0^x}
   \{
     t_z(p - z_x) + s_z(z_y - q) + t_p(r - p_x) + s_p(p_y - s)
     +
 \medskip
 \\
   \phantom{123}
   +
     t_q(s - q_x) + s_q(q_y - t)
   \}
   (\tau, g_2(\tau, x,y))
    g_{2x}
    (\tau,x,y)
   \: d\tau
   .
 \end{array}
 \eqno(3.29)
$$
 By (3.29), we get
$$
 \begin{array}{l}
   z_y(x,y)
   =
   q(x,y)
   ,
 \medskip
 \\
   z_x(x,y)
   =
   p(x,y)
   ,
 \medskip
 \\
   p_y(x,y)
   =
   q_x(x,y)
   =
   s(x,y,z(x,y),p(x,y),q(x,y))
   ,
 \medskip
 \\
   p_x(x,y)
   =
   r(x,y,z(x,y),p(x,y),q(x,y))
   ,
 \medskip
 \\
   q_y(x,y)
   =
   t(x,y,z(x,y),p(x,y),q(x,y))
   .
 \end{array}
$$
 By (3.25), so that
 three functions
  $z,\
 p,\ q$ are a $C^1\mbox{-solution}$ of the system (3.10)--(3.11).
 By Lemma~3.1, so that $z$ is a $C^3\mbox{-solution}$
 of the problem
 (3.1)--(3.3)
 in the domain $T$.
 $\Box$

 \textbf{Remark.}
 The system (3.16)--(3.20)
 consists of five equations
 with respect to five unknown functions
 $u_1,\ u_2,\ z,\ p,\ q$. If coefficients
 $A,\ B,\ C,\ D$ of the equation (3.1)
 don't depend on
 $z$,
 then equations (3.16), (3.19), (3.20)
 are a closed system of four equations
 with respect to four unknown functions
 $u_1,\ u_2,\ p,\ q$.
 If we know
 $u_1,\ u_2,\ p,\ q$,
 then
 we can find
 $z$ from (3.18) or from the first equation (3.27).

 Let equations
 (3.16)
 be a closed system with respect to
 $u_1,\ u_2$.
 Obviously, this situation holds
 iff
$$
 {\ds
   \frac{\pt E_j}{\pt z}
   =
   \frac{\pt E_j}{\pt p}
   =
   \frac{\pt E_j}{\pt q}
   =
   \frac{\pt I_j}{\pt z}
   =
   \frac{\pt I_j}{\pt p}
   =
   \frac{\pt I_j}{\pt q}
   =
   0
 }
 \eqno(3.30)
$$
 ${(j = 0,\dots, 5)}$. Here we get $E_j,\ I_j$ from
 (3.17). In this case, after finding
 ${u_1,\ u_2}$, one could get functions $z,\ p,\ q$
 from
 (3.18)--(3.20)
 or from (3.11).



\subsection{Final form of the system in Riemann invariants}

 The system (3.16) was obtained by D.~V.~Tunitsky
 \cite{Tunitsky-dep}.
 Right parts of this system are polynomial,
 generators are
 $u_1,\ u_2,$
 coefficients contain unknown functions
 $p,\ q.$
 So it isn't easy to formulate conditions on the coefficients.
 We propose another form for the system
 (3.16). Let generators be unknown functions $u_1, u_2, p, q$,
 and let coefficients be known functions.

 By definition, put
$$
   {r = u_1, \quad s = u_2.}
 \eqno(3.31)
$$
 Here $r,\ s$ are characteristic variables
 (\cite{Cour}, Supplement~1 to Chapter~V, \S~2).
 Also $r,\ s$ are
 Riemann invariants and eigenvalues of the system (3.16)
 (see \cite{RYa}, Chapter~1). We underline that
 in this context
 $r,\ s$
 aren't
 the Monge notations for second derivatives
 from
 (3.25).
 Both notations (3.31) and (3.25) are traditional.

 The final system is
$$
    \left(
      \partial_x
           +\xi(\omega) \partial_y
    \right)
    \omega = f_\omega(x,y,r,s,p,q,z),
 \eqno(3.32)
$$
 here $\omega$ is an index, ${\omega \in \{ r,s,p,q,z \} }$;
 the function ${\xi(\omega)}$ is
$$
 \xi(r)=s, \quad
 \xi(s)=r, \quad
 \xi(p)=s, \quad
 \xi(q)=r, \quad
 \xi(z)=r \, ;
$$
 functions ${f_\omega}$ are
$$
\begin{array}{l}
f_r=\rho_0+\rho_1r+\rho_2s+\rho_3pr+\rho_4qr+\rho_5ps+\rho_6qs+\rho_7r^2
+\rho_8rs+\rho_9pr^2+
\smallskip \\
\phantom{f^r}+\rho_{10}qr^2+\rho_{11}prs+\rho_{12}qrs+\rho_{13}r^2s
=f_r(\rho ,r,s,p,q),
 \smallskip \\
 f_s(\sigma ,r,s,p,q)=f_r(\sigma ,s,r,p,q),
 \smallskip \\
 f_p=\pi_0+\pi_1 s,\quad
 f_q=\kappa_0+\kappa_1 r,\quad
 f_z=p+qr;
\end{array}
$$
 vector-functions $\rho ,\sigma ,\pi ,\kappa$ depend on $x,y,z,p,q$,
$$
\begin{array}{ll}
\rho_{0\phantom{0}}=D_q\: , & \sigma_{0\phantom{0}}=D_q\: ,
\smallskip \\
 \rho_{3\phantom{0}}=\phantom{-}\frac{1}{2\Dt}(C+\Dt)_z\:
,&
            \sigma_{3\phantom{0}}=-\frac{1}{2\Dt}(C-\Dt)_z\: ,
\smallskip \\
\rho_{4\phantom{0}}=\phantom{-}\frac{1}{\Dt}D_z\: ,&
            \sigma_{4\phantom{0}}=-\frac{1}{\Dt}D_z\: ,
\smallskip \\
\rho_{5\phantom{0}}=-\frac{1}{2\Dt}(C+\Dt)_z\: ,&
            \sigma_{5\phantom{0}}=\phantom{-}\frac{1}{2\Dt}(C-\Dt)_z\: ,
\smallskip \\
\rho_{6\phantom{0}}=-\frac{1}{\Dt}D_z\: ,&
            \sigma_{6\phantom{0}}=\phantom{-}\frac{1}{\Dt}D_z\: ,
\smallskip \\
\rho_{9\phantom{0}}=-\frac{1}{\Dt}B_z\: ,&
            \sigma_{9\phantom{0}}=\phantom{-}\frac{1}{\Dt}B_z\: ,
\smallskip \\
\rho_{10}=-\frac{1}{2\Dt}(C-\Dt)_z\: ,\phantom{1234567}&
                      \sigma_{10}=\phantom{-}\frac{1}{2\Dt}(C+\Dt)_z\: ,
\smallskip \\
\rho_{11}=\phantom{-}\frac{1}{\Dt}B_z\: ,&
                      \sigma_{11}=-\frac{1}{\Dt}B_z\: ,
\smallskip \\
\rho_{12}=\phantom{-}\frac{1}{2\Dt}(C-\Dt)_z\: ,&
                      \sigma_{12}=-\frac{1}{2\Dt}(C+\Dt)_z\: ,
\smallskip \\
\rho_{13}=-B_p\: ,&\sigma_{13}=-B_p\: ,
\end{array}
$$
$$
\begin{array}{l}
\rho_1=\phantom{-}\frac{1}{2\Dt} ( \phantom{-}
\Dt_x+C_x+2D_y+CD_p-D\Dt_p-DC_p-2BD_q+\frac{1}{2}CC_q-\frac{3}{2}\Dt
C_q+
\smallskip \\
\phantom{\rho_1=} +\frac{1}{2}C\Dt_q+\frac{1}{2}\Dt\Dt_q-\Dt D_p
)\: ,
\smallskip \\
\rho_2=-\frac{1}{2\Dt} ( \phantom{-}
\Dt_x+C_x+2D_y+CD_p-D\Dt_p-DC_p-2BD_q+\frac{1}{2}CC_q+\frac{1}{2}\Dt
C_q+
\smallskip \\
\phantom{\rho_1=} +\frac{1}{2}C\Dt_q+\frac{1}{2}\Dt\Dt_q+\Dt D_p
)\: ,
\smallskip \\
\rho_7=-\frac{1}{2\Dt} (
-\Dt_y+C_y+2B_x+CB_q+B\Dt_q-BC_q-2DB_p+\frac{1}{2}CC_p-\frac{1}{2}\Dt
C_p-
\smallskip \\
\phantom{\rho_1=} -\frac{1}{2}C\Dt_p+\frac{1}{2}\Dt\Dt_p-\Dt B_q
)\: ,
\smallskip \\
\rho_8=\phantom{-}\frac{1}{2\Dt} (
-\Dt_y+C_y+2B_x+CB_q+B\Dt_q-BC_q-2DB_p+\frac{1}{2}CC_p+\frac{3}{2}\Dt
C_p-
\smallskip \\
\phantom{\rho_1=} -\frac{1}{2}C\Dt_p+\frac{1}{2}\Dt\Dt_p+\Dt B_q
)\: ,\smallskip \\
\sigma_1=-\frac{1}{2\Dt} (
-\Dt_x+C_x+2D_y+CD_p+D\Dt_p-DC_p-2BD_q+\frac{1}{2}CC_q+\frac{3}{2}\Dt
C_q-
\smallskip \\
\phantom{\rho_1=} -\frac{1}{2}C\Dt_q+\frac{1}{2}\Dt\Dt_q+\Dt D_p
)\: ,
\smallskip \\
\sigma_2=\phantom{-}\frac{1}{2\Dt} (
-\Dt_x+C_x+2D_y+CD_p+D\Dt_p-DC_p-2BD_q+\frac{1}{2}CC_q-\frac{1}{2}\Dt
C_q-
\smallskip \\
\phantom{\rho_1=} -\frac{1}{2}C\Dt_q+\frac{1}{2}\Dt\Dt_q-\Dt D_p
)\: ,
\end{array}
$$
$$
\begin{array}{l}
\sigma_7=\phantom{-}\frac{1}{2\Dt} ( \phantom{-}
\Dt_y+C_y+2B_x+CB_q-B\Dt_q-BC_q-2DB_p+\frac{1}{2}CC_p+\frac{1}{2}\Dt
C_p+
\smallskip \\
\phantom{\rho_1=} +\frac{1}{2}C\Dt_p+\frac{1}{2}\Dt\Dt_p+\Dt B_q
)\: ,
\smallskip \\
\sigma_8=-\frac{1}{2\Dt} ( \phantom{-}
\Dt_y+C_y+2B_x+CB_q-B\Dt_q-BC_q-2DB_p+\frac{1}{2}CC_p-\frac{3}{2}\Dt
C_p+
\smallskip \\
\phantom{\rho_1=} +\frac{1}{2}C\Dt_p+\frac{1}{2}\Dt\Dt_p-\Dt B_q
)\: ,
\end{array}
$$
$$
\begin{array}{ll}
\pi_0=-D\: ,                  & \kappa_0=\frac{1}{2}(C+\Dt) \: ,
\\
\pi_1=\frac{1}{2}(C+\Dt)\: ,\phantom{1234567} & \kappa_1=-B \: .
\end{array}
$$

 Taking into account (3.31), we get
 initial conditions (3.21) for the system (3.32):
$$
\begin{array}{l}
  r(0,y)=r^0(y)=\frac{\ds (C+\Dt)(0,y,z^0(y),p^0(y),z^0_y(y))-2p^0_y(y)}
           {\ds 2(z^0_{yy}(y)+B(0,y,z^0(y),p^0(y),z^0_y(y)))}\; ,
  \medskip \\
  s(0,y)=s^0(y)=\frac{\ds (C-\Dt)(0,y,z^0(y),p^0(y),z^0_y(y))-2p^0_y(y)}
           {\ds 2(z^0_{yy}(y)+B(0,y,z^0(y),p^0(y),z^0_y(y)))}\; ,\medskip \\
  p(0,y)=p^0(y)\; ,\quad
  q(0,y)=q^0(y)=z^0_y(y)\; ,\quad
  z(0,y)=z^0(y)\; .
\end{array}
 \eqno(3.32^0)
$$

 Let the system $(3.32^0)$ be
$$
 \left\{
 \begin{array}{l}
 p_y^0 = - \frac{\ds \Delta}{\ds 2}
 \frac{\ds r^0 + s^0}{\ds r^0 - s^0} +
 \frac{\ds C}{\ds 2} \, ,
 \medskip
 \\
 q_{y}^0 =
 \frac{\ds \Delta}{\ds r^0 - s^0} - B \, ,
 \medskip
 \\
 z_y^0 = q^0 \, .
 \end{array}
 \right.
 \eqno (3.32^{00})
$$
 Let $r^0,\ s^0$
 be initial functions.
 Therefore we get
 $z^0,\ p^0$
 from $(3.32^{00})$.
 If coefficients $A,B,C,D$ depend on $x,\ y$ only,
 then the system $(3.32^{00})$ is linear.
 In general case, the system $(3.32^{00})$ is nonlinear.

 If coefficients $A,B,C,D$ of the Monge--Ampere equation
 depend on $x,\ y$ only,
 then the equation could be reduced to the system
 with respect to $r,\ s$
$$
\begin{array}{l}
    ( \partial_x + s \, \partial_y )
    r
    =
    (r - s)(a_{r 1}+a_{r 2}\, r)\, ,
  \medskip \\
    ( \partial_x + r \, \partial_y )
    s
    =
    (r - s)(a_{s 1}+a_{s 2}\, s)\, ,
\end{array}
 \eqno (3.33)
$$
where
$$
\begin{array}{ll}
    a_{r1}=\frac{\ds 1}{\ds 2\Dt}(2D_y+C_x+{\Dt}_x)\, , &
    a_{r2}=\frac{\ds 1}{\ds 2\Dt}(-2B_x-C_y+{\Dt}_y)\, ,
    \medskip
    \\
    a_{s1}=\frac{\ds 1}{\ds 2\Dt}(2D_y+C_x-{\Dt}_x)\, , &
    a_{s2}=\frac{\ds 1}{\ds 2\Dt}(-2B_x-C_y-{\Dt}_y) \, .
\end{array}
$$
 The system (3.33) is a system of the first and the second
 equations from the system (3.32).
 After getting
 $r,\ s$ (from (3.33))
 we solve the third, the fourth, and the fifth linear equations
 (3.32) with respect to $p, \ q, \ z$.
 By $(3.32^{00})$,
 we get
 initial functions $p^0,\ q^0,\ z^0$ from
 $r^0,\ s^0$.



\section{Successful approximations}
\subsection{Iterative loop}

   Let ${{\up{n}{\omega}}(x,y)},\;
 {\omega=r,s,p,q,z,}$
 be known functions.
 Let
 ${{\up{n+1}{\omega}}(x,y)},\
 {\omega=r,s,p,q,z,}$
 be a solution of the nonlinear Cauchy problem
$$
\left\{
  \begin{array}{l}
    (
        \pt_x
        +
        {\up{n+1}{\xi}}(\omega)
        \pt_y
    )
    \up{n+1}{\omega}
    =
    f_\omega(x,y,\up{n}{r},\up{n}{s},\up{n}{p},\up{n}{q},\up{n}{z}),
    \smallskip \\
    {\up{n+1}{\omega}}
    (0,y)=\omega^0(y),
    \quad \omega=r,s,p,q,z.
  \end{array}
\right.
 \eqno(4.1)
$$

 The Cauchy problem (4.1) falls into four
 independent problems.
 The first is a problem for a nonlinear system
 with respect to
 ${\up{n+1}{r}}$, ${\up{n+1}{s}}$.
 Another are three problems for independent linear equations
 with respect to
 $p,\ q,\ z.$
 Before solving the nonlinear system we'll get a priori bounds.

 By definition, put
$$
  {{\up{0}{\omega}}(x,y)}
  =
 \omega^0(y), \quad \omega=r,s,p,q,z.
 \eqno(4.2)
$$

 By (3.32),
 it follows that
 the vector-function $\omega_x$
 is defined by $\omega,\ \omega_y$, ${\omega=r,s,p,q,z}$.
 Suppose the existence of continuous functions
 $\omega,\ \omega_y,$
 then there exists a $C^1\mbox{-solution}$ of (3.32), $(3.32^0)$.
 Proof of the existence of continuous vector-functions
 $\omega,\ \omega_y$
 is a proof of uniform convergence of
 $\{{\up{n}{\omega}}\},\ \{{\up{n}{\omega}_y}\}, \ {\omega=r,s,p,q,z.}$
 The main part of
 this proof
 is uniform boundedness of
 $\{{\up{n}{\omega}}\},\ \{{\up{n}{\omega}_y}\}.$


\subsection{Uniform boundedness}

   Suppose vector-functions
 $\rho, \sigma, \pi, \kappa$
 and functions
  $r^0,\: s^0,\: p^0,\: q^0$
 are
 ${C^1\mbox{-smooth}}$
 and bounded,
 ${z^0 \in C^1(\mathbb{R})}$.

 Consider
 functions ${f_r,\, f_s}$
 as right parts of (3.32).
 There are two kinds of monomials.
 Some monomials contain generators ${p,\ q}$,
 another monomials don't contain generators ${p,\ q}$.
 So we distinguish coefficients
 ${\rho, \, \sigma}$
 for the two cases
 by introducing the next two sets of indexes.
 By definition, put
 ${J_{rs} = \{ 0,1,2,7,8,13 \} }$
 for the first case (no ${p,\ q}$),
 and put
 ${J_{pq} = \{ 3,4,5,6,9,10,11,12 \} }$
 for the second case (${p,\ q}$ are).

 By definition, put
$$
\begin{array}{rcl}
    U_0&=&\max\lms_{\omega=r,s}\:\sup\lms_{y\in \mathbb{R}}
    |\omega^0(y)|
    = \hbox{const},
 \medskip\\
    \alpha_1(x)&=&
      \sup\lms_
      {
        {
            {
            \atop\scriptstyle (y,z,p,q)\in \mathbb{R}^4
            }
            \atop
            {
            \atop\scriptstyle j \, \in \, J_{rs}
            }
        }
        \atop
        {
          \atop\scriptstyle a \, \in \, \{ \rho, \sigma \}
        }
      }
      |a_j(x,y,z,p,q)|,
 \medskip\\
    \alpha_2(x)&=&
      \sup\lms_
      {
        {
            {
            \atop\scriptstyle (y,z,p,q)\in \mathbb{R}^4
            }
            \atop
            {
            \atop\scriptstyle j \, \in \, J_{pq}
            }
        }
        \atop
        {
          \atop\scriptstyle a \, \in \, \{ \rho, \sigma \}
        }
      }
      |a_j(x,y,z,p,q)|,
 \medskip\\
    \alpha_3 &=&
      \sup\lms_
      {
        {
            {
            \atop\scriptstyle (x,y,z,p,q)\in \mathbb{R}^5
            }
            \atop
            {
            \atop\scriptstyle j=0,1
            }
        }
        \atop
        {
          \atop\scriptstyle a \, \in \, \{ \pi, \kappa \}
        }
      }
      |a_j(x,y,z,p,q)|
      = \hbox{const}.
\end{array}
 \eqno(4.3)
$$

{\bf Lemma~4.1.}
 Suppose
$$
\begin{array}{l}
  |p^0|\le 1
   \; , \quad
   |q^0|\le 1
   \; ,
  \bigskip \\
  U_0
  +
  6{\ds\int_{-\infty}^{+\infty}}
    \alpha_1(x)
  \: dx
  +
  8{\ds\int_{-\infty}^{+\infty}}
    (1 + 2\alpha_3 |x|)\alpha_2(x)
  \: dx
  \le 1
  \; .
\end{array}
  \eqno(4.4)
$$
 Suppose there exists ${n \ge 0}$ such that
$$
\begin{array}{ll}
 |\up{n}{\omega}(x,y)|
 \le
 1
  \: ,
 &
 \omega=r,s
 \; ,
 \smallskip
 \\
 |\up{n}{\omega}(x,y)|
 \le
 1 + 2\alpha_3 |x|
 \: ,
 &
 \omega=p,q
 \; ,
 \smallskip
 \\
 |\up{n}{z}(x,y)|
 \le
   \max\lms_{t \in [y - |x|,\ y + |x|]}
   |z^0(t)|
   + 2|x| + 2\alpha_3 x^2,
\end{array}
  \eqno(4.5)
$$
 for ${\forall \: (x,y)\in [0,+\infty )\times \mathbb{R}}$.
 Suppose the existence of functions
 ${\up{n+1}{\omega}},\ \omega = r,s,p,q,z$.
 Then we have the same bounds (4.5)
 for the number ${n+1}$.

{\bf Proof.} By (4.2), (4.4), we obtain (4.5) for ${n=0}$.
 Suppose (4.5) holds for some number $n$;
 then we shall prove (4.5) for the number ${n+1}$.

 Consider right parts $f_\omega$ of the system (3.32).
 By (2.10), we obtain the next estimations.
 For $\omega=r,s$ we have
$$
\begin{array}{rl}
  |{\up{n+1}{\omega}}(x,y)|
  &
  \le
  |\omega^0|
  +
  {\ds \int_0^x}
  |f_\omega|
  \: d\tau
  \le
 \medskip
 \\
  &
 \le
 U_0
 +
 6
 {\ds \int_0^{+\infty}}
 \!\!\!
   \alpha_1(\tau)
 \: d\tau
 +
 4
 {\ds \int_0^{+\infty}}
 \!\!\!
 \alpha_2(\tau)
 \: |{\up{n}{p}}|
 \: d\tau
 +
 4
 {\ds \int_0^{+\infty}}
 \!\!\!
 \alpha_2(\tau)
 \: |{\up{n}{q}}|
 \: d\tau
 \le
 \medskip
 \\
  &
 \le
 U_0
 +
 6
 {\ds \int_0^{+\infty}}
 \!\!\!
   \alpha_1(\tau)
 \: d\tau
 +
 8
 {\ds \int_0^{+\infty}}
 \!\!\!
 \alpha_2 (\tau)
 (1 + 2\alpha_3 \tau)
 \: d\tau
 \le
 1
 \: .
\end{array}
$$
 For $\omega=p,q,\; a=\pi,\kappa$
 we have
$$
\begin{array}{rl}
  |{\up{n+1}{\omega}}(x,y)|
  &
  \le
  |\omega^0|
  +
  {\ds \int_0^x}
  |f_\omega|
  \, d\tau
  \le
  \phantom{12}
    \langle\mbox{see (4.3) for }\alpha_3\rangle
 \medskip\\
  &
 \le
 1
 +
 {\ds \int_0^x}
 (|a_0|+|a_1|)
 \, d\tau
 \le
 1
 +
 2\alpha_3
 {\ds \int_0^x}
 d\tau
 \le
 1
 +
 2\alpha_3
 \:
 |x|
 \: .
\end{array}
$$
 For ${\omega=z,}$ taking into account ${|{\up{n}{r}}|\le 1,}$ that is
 ${{\up{n+1}{g_z}}(\tau,x,y)\in [y - x + \tau, \ \ y + x -
 \tau]}$,
 we have
$$
\begin{array}{l}
  |{\up{n+1}{z}}(x,y)|
  \le
  |z^0|
  +
  {\ds \int_0^x}
  |\up{n}{p}+\up{n}{q}\up{n}{r}|
  \, d\tau
  \le
 \max\lms_{t\in [y - |x|,\ y + |x|]} |z^0(t)|
 +
 \medskip\\
 +
 {\ds \int_0^x}
     2(1 + 2\alpha_3 \tau)
 \, d\tau
 =
 \max\lms_{t\in [y - |x|,\ y + |x|]} |z^0(t)|
 +
 2 |x|
 +
 2\alpha_3 x^2
 \: .
 \ \Box
\end{array}
$$

{\bf Corollary.} Suppose conditions (4.4) are satisfied; then
${\{\up{n}{\omega}\}}\: ,\; {\omega=r,s,p,q,z}$ are uniformly
bounded on the compactum
$$
  G(\bar x,\bar y)=\{(x,y)|\; x\in[0,\bar x],\;
   y\in[\bar y - \bar x + x,\ \ \bar y + \bar x - x]\}
   \eqno (4.6)
$$
 for
 ${\forall \: (\bar x,\bar y)\in [0,+\infty)\times\mathbb{R}}$.

{\bf Proof.}
 Consider any point ${(x,y)\in G(\bar x,\bar y)}$.
 By ${|{\up{n}{r}}|\le 1,}\; {|{\up{n}{s}}|\le 1}$,
 so that
 outgoing (from ${(x,y)}$) characteristics
 are being inside the compactum
 ${G(\bar x,\bar y)}$.
 Therefore, functions
 $\up{n}{\omega},\; {\omega=r,s,p,q,z}$,
 are defined on the compactum
 ${G(\bar x,\bar y)}$.
 By (4.5),
 we obtain uniform boundedness.
 $\Box$


\subsection{Hyperbolicity in the restricted sence}

 Let the initial data $r^0, s^0$
 be separated by
 some constant.
 More exactly,
 $\exists \, \delta > 0$
 such that
$$
  \inf \lms_{\scriptstyle y\in \mathbb{R}} r^0(y)
  -
  \sup \lms_{\scriptstyle y\in \mathbb{R}} s^0(y)
  \geq \delta > 0 \, .
  \eqno (4.7)
$$

{\bf Lemma~4.2.} Suppose
 $\exists \, \varepsilon \in (0,\delta]$
 such that
$$
  6
  {\ds\int_{-\infty}^{+\infty}}
      \alpha_1(x)
  \: dx
  +
  8
  {\ds\int_{-\infty}^{+\infty}}
    (1 + 2\alpha_3 |x|)
    \alpha_2(x)
  \: dx
  \leq
  \frac{\delta - \varepsilon}{2}
  \, .
  \eqno (4.8)
$$
 Then for $n=0,1,2,\dots$
$$
\inf \lms_{\scriptstyle (x,y)\in \mathbb{R}^2} \up{n}{r}(x,y) -
  \sup \lms_{\scriptstyle (x,y)\in \mathbb{R}^2} \up{n}{s}(x,y)
  \geq \varepsilon > 0 .
  \eqno (4.9)
$$

{\bf Proof.} By (2.10), taking into account (4.3), (4.5), (4.8),
 we have
$$
\begin{array}{rl}
 {\up{n+1}{r}}(x,y)
 &
 \ge
 \inf \lms_{\scriptstyle y\in \mathbb{R}} r^0(y)
 -
 \left(
  6
  {\ds\int_{-\infty}^{+\infty}}
    \alpha_1(x)
  \: dx
  +
  8
  {\ds\int_{-\infty}^{+\infty}}
    (1 + 2\alpha_3 |x|)
    \alpha_2(x)
  \: dx
 \right)
 \ge
 \medskip
 \\
 &
 \ge
 \inf \lms_{\scriptstyle y\in \mathbb{R}} r^0(y)
 - \frac{\ds \delta - \varepsilon}{\ds 2} \, ,
 \medskip
  \\
 {\up{n+1}{s}}(x,y)
 &
 \le
 \sup \lms_{\scriptstyle y\in \mathbb{R}} s^0(y)
 +
 \left(
  6{\ds\int_{-\infty}^{+\infty}}\alpha_1(x)dx
  +
  8
  {\ds\int_{-\infty}^{+\infty}}
  (1 + 2\alpha_3 |x|)
  \alpha_2(x)
  \: dx
 \right)
 \leq
 \medskip
 \\
 &
 \leq
 \sup \lms_{\scriptstyle y\in \mathbb{R}}
 s^0(y)
 +
 \frac{\ds \delta - \varepsilon}{\ds 2}
 \, .
\end{array}
$$
 Then, by (4.7),
$$
 \inf \lms_{\scriptstyle (x,y)\in \mathbb{R}^2}
 {\up{n+1}{r}}(x,y)
 -
 \sup \lms_{\scriptstyle (x,y)\in \mathbb{R}^2}
 {\up{n+1}{s}}(x,y)
 \geq
 \inf \lms_{\scriptstyle y\in \mathbb{R}} r^0(y)
 -
 \sup \lms_{\scriptstyle y\in \mathbb{R}} s^0(y)
 - (\delta - \varepsilon) \geq \varepsilon
 \; .
 \;\square
$$


\subsection{Solvability of the iterative system}

 \textbf{Lemma~4.3.} There exists a ${C^1\mbox{-solution}}$
 of (4.1)
 on the whole plain.
 \\
 \textbf{Proof.} Consider the subsystem of two equations
 with respect to
 $\up{n+1}{r}$, $\up{n+1}{s}$.
 It is weakly nonlinear
 (\cite{RYa}, Chapter~1, ${\S}~10,$ Subsection~3).
 (For (2.8) weak nonlinearity is
 ${\ds {\pt \xi_i}/{\pt u_i} = 0}$ for all $i$.)
 By (4.5), the solution of this system is bounded
 on the whole half-plane.
 By (4.9), the system is hyperbolic in the restricted sense.
 By the Rozhdestvensky---Sydorenko theorem
 (\cite{RYa}, Chapter~1, ${\S}~10,$ Subsection~3),
 first derivatives of functions
 $\up{n+1}{r}$, $\up{n+1}{s}$
 aren't infinite at finite $x$.
 Therefore, by the corollary of this theorem
 (\cite{RYa}, Chapter~1, ${\S}~10,$ Subsection~3),
 the subsystem of the first and the second equations of (4.1)
 is solvable on the whole (half-)plane,
 that is, it has a global $C^1\mbox{-smooth}$ solution.

 Indeed, the problem (4.1) is solvable locally,
 i.e. in some neighborhood of the axis
 ${x = 0}$. We get this well-known fact,
 for example,
 from the existence theorem
 (\cite{RYa}, Chapter~1, ${\S}~8,$ Subsection~2),
 taking into account Corollary of
 Lemma~4.1.
 Consider prolongation of the local solution.
 By the Rozhdestvensky---Sydorenko theorem,
 a strong break is impossible.
 (A strong break is an infinite first derivative at finite $x$.)
 In the proof of this theorem some majorant
 of the (modulus of the) classical solution
 is constructed.
 So,
 before going to infinity,
 the solution (it is classical yet)
 must go outside the majorant,
 but it is impossible.

 Further, from studying weak breaks
 (finite jumps of first derivatives),
 \textit{
 let a solution and its first derivatives be bounded,
 and let a weak break of a hyperbolic quasilinear system
 be
 propagating along a
 characteristic;
 then the weak break
 can't arise or disappear}
 (\cite{RYa}, Chapter~1, ${\S}~10,$ Subsection~1).
 We have smooth initial functions,
 so weak breaks aren't.
 Thus,
 there exists
 a $C^1\mbox{-smooth}$
 solution on the whole half-plane.

 Another three equations of the system
 (4.1)
 are linear.
 Therefore,
 their Cauchy problems
 (4.1)
 are $C^1\mbox{-solvable}$ on the whole half-plane.
 $\Box$

 \textbf{Lemma~4.4.}
 Let
 ${{\up{n+1}{\omega}}(x,y)}$
 be a solution of
 (4.1).
 Then
 for $\forall$
 $\tau \in [0,x]$
 there exist
 characteristics of the problem
 (4.1).

 \textbf{Proof.}
 By definition,
 characteristics are solutions of
 (2.9).
 For systems
$$
 \frac{d}{d \tau}
 \ u_i
 =
 f_i (\tau , u_1, \dots , u_n),
 \quad i=1, \dots , n
$$
 we have the next classical Cauchy theorem
 of the existence and uniqueness.

 {\it
 Consider a closed domain
$$
 \bar G =
 \{
   (\tau , u_1, \dots , u_n)
   | \quad
   |\tau - \tau^0| \leq a, \quad
   |u_i - u_i^0| \leq b, \quad
   i=1, \dots , n.
 \}
$$
 Suppose the next conditions hold in this domain:\\
 1) functions $f_i$ are continuous;\\
 2) ${|f_i| \leq A}$;\\
 3) functions $f_i$ are Lipschitz
    with respect to ${u_1, \dots , u_n}$.\\
 Then
 for
 $\forall\ \tau$ such that
 ${|\tau - \tau^o| \leq \min (a, b/A)}$
 there exists a unique solution
 of the system
 with initial functions
 ${u_i(\tau^0) = u^0}$.
 }

 In our case, by uniform estimation
 (4.5),
 we have
 ${A = 1}$.
 Further, ${a = x}$,
 and any big value of $b$
 could be chosen.
 Partial derivatives with respect to
 ${u_1, \dots , u_n}$
 are bounded in
 $\bar G$,
 because,
 by Lemma~4.3,
 right parts
 are a classical solution of (4.1),
 so they belong to
 ${C^1([0, +\infty) \times \mathbb{R})}$.
 Therefore,
 there exist both characteristics
 for
 ${\forall \: \tau \in [0, x]}$.
 Thus we can integrate along characteristics.
 $\ \Box$



\section{Derivatives of successful approximations}

 The Rozhdestvensky---Sydorenko theorem
 (\cite{RYa}, Chapter~1, ${\S}~10,$ Subsection~3)
 was proved for an exact solution of a system of two equations.
 We expand it to successful approximations for
 a system of five equations with two different characteristics.

 When we write
 ${{(\up{n}{u})}_y}$,
 we mean that the first operation is getting the
 ${n\mbox{-th}}$ successful approximation for $u$,
 and the second operation is differentiation with respect to $y$.
 Not conversely.

 Let ${G(\bar x, \bar y)}$
 be a compactum
 (4.6),
 ${\forall (\bar x,\bar y)\in [0,+\infty )\times \mathbb{R}}$.

 {\bf Lemma~5.1.}
 There exists a function
 $\Phi (x) \in C^0(\mathbb{R})$
 such that
$$
|{(\up{n}{\omega})_y} (x,y)| \leq \Phi (x)
 \eqno (5.1)
$$
 for $\omega = r,s,p,q,z$, $\forall (\bar x,\bar y)\in G(\bar x,
 \bar y)$, $n=0,1,2,\dots$

 {\bf Proof.}
  By definition, put
  ${\up{n}{g}_{\omega}}
 (x,y_0) = {\up{n}{g}_{\omega}}
 (x,0,y_0)$.
 Then ${\up{n}{g}_{\omega}} (x,y_0)$
 is a solution of the Cauchy problem
 $$
\left\{
 \begin{array}{l}
    \pt_x
    {\up{n}{g}}_{\omega}
    =
    \xi ({\up{n}{\omega}})
    (x, \up{n}{g}_{\omega}),
 \smallskip \\
    {\up{n}{g}_{\omega}}(0,y_0) = y_0, \quad \omega = r,s,p,q,z.
 \end{array}
 \right.
 \eqno (5.2)
 $$
 It means that the curve $(x, {\up{n}{g}_{\omega}}(x,y_0))$
 is a characteristic that goes through
 the point
 $(0, y_0)$.
 The formula (2.10) is
$$
{\up{n+1}{\omega}} (x, {\up{n+1}{g_\omega}} (x,y_0)) =
 \omega^0 (y_0) + \int_0^x f_{\omega} (\tau, {\up{n+1}{g_\omega}} (\tau,y_0),
  {\up{n}{\Omega}} (\tau, {\up{n+1}{g_\omega}} (\tau,y_0)) )\; d\tau
  \, ,
$$
where $\Omega = (\omega) = (r,s,p,q,z)$.
 By definition, put
$$
{\up{n+1}{\bar\omega}} (x,y_0) = {\up{n+1}{\omega}} (x,
{\up{n+1}{g_\omega}} (x,y_0)) \, .
 \eqno (5.3)
$$
 Then for ${\omega = r,s,p,q,z}$
 we obtain
$$
 {\up{n+1}{\bar\omega}} (x, {\up{n+1}{g_\omega}} (x,y_0)) =
 \omega^0 (y_0) + {\ds \int}_0^x f_{\omega} (\tau, {\up{n+1}{g_\omega}} (\tau,y_0),
  {\up{n}{\Omega}} (\tau, {\up{n+1}{g_\omega}} (\tau,y_0)) )\; d\tau
  \, .
 \eqno (5.4)
$$
 Therefore,
 ${\up{n+1}{\bar\omega}} (x,y_0)$
 is a solution of the Cauchy problem
$$
\left\{
\begin{array}{l}
 \pt_x {\up{n+1}{\bar\omega}} =
  f_{\omega} (x, {\up{n+1}{g_\omega}} (x,y_0),
  {\up{n}{\Omega}} (x, {\up{n+1}{g_\omega}} (x,y_0)) )\, ,
  \smallskip \\
 {\up{n+1}{\bar\omega}}(0,y_0) = \omega^0 (y_0), \quad
 \omega = r,s,p,q,z.
\end{array}
\right.
 \eqno (5.5)
$$
 Differentiating (5.3) with respect to $y_0$,
 we get
$$
 \pt_{y_0} {\up{n+1}{\bar\omega}} (x, y_0) =
 \left.
 \pt_y {\up{n+1}{\bar\omega}} (x, y)
 \right|_{y = {\up{n+1}{g_\omega}} (x,y_0)}
 \pt_{y_0} {\up{n+1}{g_\omega}} (x,y_0) \, ,
$$
 therefore,
$$
 \pt_y {\up{n+1}{\omega}} (x, y) =
 \pt_{y_0}{\up{n+1}{\bar\omega}} (x, y_0)
 \pt_{y_0}
 {\up{n+1}{g_\omega}} (x,y_0) \, .
 \eqno (5.6)
$$
 Differentiating (5.2) with respect to $y_0$,
 we get
$$
\pt_x
 \left(\pt_{y_0}
 {\up{n}{g}_{\omega}} (x,y_0)
 \right) =
 \left.
 \pt_y \xi ({\up{n}{\omega}}) (x,y)
 \right|_{y = {\up{n}{g}_{\omega}} (x,y_0)}
 \pt_{y_0}
 {\up{n}{g}_{\omega}} (x,y_0) \, ,
$$
 that is
$$
 \pt_x
 \ln
 \left(
 \pt_{y_0}
 {\up{n}{g}_{\omega}} (x,y_0)
 \right)
 =
 \left.
 \pt_y \xi ({\up{n}{\omega}}) (x,y)
 \right|_{y = {\up{n}{g}_{\omega}} (x,y_0)}
 \, ,
\eqno (5.7)
$$
 and the initial condition
$$
 (
 \pt_{y_0}
 {\up{n}{g}_{\omega}}
 )
  (0,y_0) = 1 \, .
 \eqno (5.7^0)
$$
 Let $v(x,y)\in C^1$
 be an arbitrary function.
 By definition, put
$$
 \left(
     \frac{d}{dx}
     \: v
 \right)_{\omega}
 =
 (
    \pt_x
    +
    {\up{n+1}{\xi}} (\omega)
    \pt_y
 )
 \: v
 \, .
$$
 Subtracting from
$$
 (
    \pt_x
    +
    {\up{n+1}{r}} \pt_y
 )
 {\up{n+1}{s}}
 =
 f_s(x,y,{\up{n}{\Omega}})
$$
 the equality
$$
 (
    \pt_x
    +
    {\up{n+1}{s}} \pt_y
 )
 {\up{n+1}{s}}
 =
 \left(
     \frac{d}{dx} {\up{n+1}{s}}
 \right)_r \, ,
$$
 we get
$$
 \pt_y
     {\up{n+1}{s}} =
 \frac{f_s(x,y,{\up{n}{\Omega}})
         -
         (\frac{\ds d}{\ds dx}{\up{n+1}{s}})_r
      }
      {\ds {\up{n+1}{r}} - {\up{n+1}{s}}} \, .
$$
 Transforming, we get
$$
\begin{array}{rl}
 \pt_y
 {\ds\up{n+1}{s}}
 &
 =
 \
 \frac{\ds f_s(x,y,{\up{n}{\Omega}})
           -
           f_r(x,y,{\up{n}{\Omega}})
      }
      {\ds {\up{n+1}{r}} - {\up{n+1}{s}}}
 +
 \frac{{\ds f_r(x,y,{\up{n}{\Omega}})} \, - \,
 {\ds (}
 \frac{\ds d}{\ds dx}{\ds\up{n+1}{s}}
        {\ds )}_r
      }
      {\ds {\up{n+1}{r}} - {\up{n+1}{s}}}
 =
  \smallskip \\
 &
 =
 \
 \frac{\ds f_s(x,y,{\up{n}{\Omega}})
           -
           f_r(x,y,{\up{n}{\Omega}})
      }{\ds {\up{n+1}{r}} - {\up{n+1}{s}}}
 +
 \frac{
        {\ds (}
           \frac{\ds d}{\ds dx}
           {\ds \up{n+1}{r}}
        {\ds )}_r
   \, - \,
 {\ds (}
      \frac{\ds d}{\ds dx}
      {\ds \up{n+1}{s}}
 {\ds )}_r
 }
 {\ds {\up{n+1}{r}} - {\up{n+1}{s}}}
 =
  \smallskip \\
 &
 =
 \
 \frac{\ds f_s(x,y,{\up{n}{\Omega}}) -
           f_r(x,y,{\up{n}{\Omega}})
      }
      {\ds {\up{n+1}{r}} - {\up{n+1}{s}}}
 +
 \frac{
        \left(
            \frac{\ds d}{\ds dx}
            \left(
                {\ds {\up{n+1}{r}} - {\up{n+1}{s}} }
            \right)
        \right)_r
      }
      {\ds {\up{n+1}{r}} - {\up{n+1}{s}}}
 =
  \smallskip \\
 &
 =
 \
 \frac{\ds f_s(x,y,{\up{n}{\Omega}}) -
           f_r(x,y,{\up{n}{\Omega}})
      }
      {\ds {\up{n+1}{r}} - {\up{n+1}{s}}}
 +
   \left(
     \frac{\ds d}{\ds dx}
     \ln
     \left(
        {\ds {\up{n+1}{r}} - {\up{n+1}{s}} }
     \right)
   \right)_r
   \, .
 \end{array}
 \eqno (5.8)
$$
 By (5.7), we have
$$
 \left.
    \left(
 \!\!
        \left(
            \pt_x
            +
            {\up{n+1}{s}}
            \pt_y
        \right)
        \ln
        \frac
        {
            \pt_{y_0}
            {\up{n+1}{g_r}} (x,y_0)
        }
        {
            {\ds {\up{n+1}{r}} - {\up{n+1}{s}} }
        }
    \right)
 \right|
        _{y = {\up{n+1}{g_r}} (x,y_0)}
 \!\!\!\!\!\!\!\!
 \!\!\!\!\!\!\!\!
 \!\!\!
 =
 \
    \left.
        \frac{\ds f_s(x,y,{\up{n}{\Omega}}) -
                f_r(x,y,{\up{n}{\Omega}})
            }
            {\ds {\up{n+1}{r}} - {\up{n+1}{s}}}
    \:
    \right|
           _{y = {\up{n+1}{g_r}} (x,y_0)}
   \, .
$$
 By (2.10), integrating along a characteristic,
 and taking into account
 $(5.7^0)$,
 we get
$$
            \pt_{y_0}
            {\up{n+1}{g_r}} (x,y_0)
 =
        \frac
        {
            {\ds {\up{n+1}{r}} - {\up{n+1}{s}} }
        }
        {\ds r^0 - s^0}
    \exp
    \left\{
        {\ds \int_0^x}
        \frac
             {\ds {\up{n}{f_s}} - {\up{n}{f_r}}}
             {\ds {\up{n+1}{r}} - {\up{n+1}{s}}}
        \: d\tau
    \right\}
   \, .
   \eqno(5.9)
$$
 Functions $f_{\omega},\ \omega = r,s,p,q,z$,
 are continuously differentiable with respect to
 $r,s,p,q,z,x,y$,
 and
 $\{ {\up{n}{\Omega}} \}$
 is uniformly bounded.
 Therefore,
 there exist a constant
 $a$
 such that for
 ${\omega = r,s,p,q,z,}$
 ${\mu = r,s,p,q,z,x,y}$
 we have
$$
|f_{\omega}|\leq a,
 \quad
 \left|
   \frac{\ds \pt f_\omega}{\ds \pt \mu}
 \right|
  \leq a.
\eqno (5.10)
$$
 By (4.5), (4.9), (5.10),
 we get from (5.9)
 the following estimate:
$$
\frac{1}{\psi (x)} \leq
   \pt_{y_0} {\ds\up{n+1}{g_r}} (x,y_)
   \leq \psi (x) \, ,
   \eqno (5.11)
$$
 where
$$
\psi (x) =
 (2/\varepsilon)
 \exp
 \left\{
     2ax/\varepsilon
 \right\}
 \, .
$$

 Estimates for $\pt_{y_0}{\up{n+1}{g_s}}(x,y_0)$
 we get analogously.
 Recall that there are two characteristics only,
 therefore, for
 ${\omega = r,s,p,q,z}$
 ${{\ds\up{n+1}{g_\omega}} (x,y_0)}$
 is
 ${{\ds\up{n+1}{g_r}} (x,y_0)}$
 or
 ${{\ds\up{n+1}{g_s}} (x,y_0)}$.
 Thus, for
 $\omega = r,s,p,q,z$, $n=1,2,3,\dots$
 we obtain estimates
$$
\frac{1}{\psi (x)} \leq
   \pt_{y_0} {\ds\up{n+1}{g_\omega}} (x,y_0)
   \leq \psi (x) \, .
   \eqno (5.12)
$$
 We don't need an estimation for ${n = 0}$.

 Differentiating (5.5) with respect to $y_0$,
 we get
$$
\left\{
 \begin{array}{l}
   \pt_x
   (
     \pt_{y_0}
     {\up{n+1}{\bar \omega}}
   )
     =
     \sum \lms_{\mu = r,s,p,q,z}
        \frac{\ds\pt f_\omega}
          {\ds\pt {\up{n}{\mu}}}
        \frac{\ds\pt {\up{n}{\mu}}}
          {\ds\pt {\up{n+1}{g_\omega}}}
        \frac{\ds\pt {\up{n+1}{g_\omega}}}
          {\ds\pt y_0}
        +
        \frac{\ds\pt f_\omega}
          {\ds\pt {\up{n+1}{g_\omega}}}
        \frac{\ds\pt {\up{n+1}{g_\omega}}}
          {\ds\pt y_0}
   \, ,
   \medskip \\
       \pt_{y_0}
     {\up{n+1}{\bar \omega}}
     (0, y_0)
     =
     \omega_y^0(y_0)
     \, , \qquad \omega = r,s,p,q,z
     \, .
 \end{array}
 \right.
 \eqno (5.13)
$$
 By definition, put
$$
V_0 = \max\lms_{\omega = r,s,p,q,z}
 \sup\lms_{y_0 \in
                  \left.
                      {G(\bar x, \bar y)}
                  \right|_{x = 0}
          }
 |\omega_y^0 (y_0)|
 \, .
$$
 Taking into account (5.13),
 consider the majorant problem
$$
        \frac{\ds d}{\ds dx}
        V
        =
        5a\psi^2 (x)\; V + a\psi (x) \, ,
  \qquad
        V(0)=V_0 \, .
 \eqno (5.14)
$$
 By linearity, the Cauchy problem (5.14)
 is solvable on the whole compactum
 ${{G(\bar x, \bar y)}}$.

 Recall that ${\psi \geq 1}$
 and
 right parts of (5.14)
 are nonnegative.
 Therefore,
 the initial approximation $\up{0}{\omega}$
 satisfy to
$$
|\pt_y {\up{0}{\omega}}|
 \leq
 \psi (x)\; V(x) \, , \quad
 \omega = r,s,p,q,z.
$$
 Suppose
$$
|\pt_y {\up{n}{\omega}}|
 \leq
 \psi (x)\; V(x) \, , \quad
 \omega = r,s,p,q,z.
 \eqno (5.15)
$$
 By (5.10), (5.12), (5.15),
 we get
$$
 \begin{array}{l}
 \left|
 \pt_{y_0} {\up{n+1}{\bar\omega}}
 \right|
 \up{(5.13),\ (2.10)}{\le}
 |\omega^0_y| +
 \medskip \\
 + {\ds\int_0^x}
 \left(
   \sum\lms_{\mu = r,s,p,q,z}
     \left|
        \frac{\ds\pt f_\omega}
          {\ds\pt {\up{n}{\mu}}}
     \right|
     \left|
        \frac{\ds\pt {\up{n}{\mu}}}
          {\ds\pt {\up{n+1}{g_\omega}}}
     \right|
     \left|
        \frac{\ds\pt {\up{n+1}{g_\omega}}}
          {\ds\pt y_0}
     \right|
     +
     \left|
        \frac{\ds\pt f_\omega}
          {\ds\pt {\up{n+1}{g_\omega}}}
     \right|
     \left|
        \frac{\ds\pt {\up{n+1}{g_\omega}}}
          {\ds\pt y_0}
     \right|
 \right)
 \; d\tau
 \le
 \medskip \\
 \le
 V_0
 +
 {\ds\int_0^x}
 (5a\psi(\tau)\psi(\tau)\; V(\tau) + a\psi(\tau))
 \; d\tau
 \up{(5.14)}{=}
 V(x)
 \, .
 \end{array}
\eqno(5.16)
$$
 From (5.6), by estimates (5.12), (5.16),
 we get
$$
 \left|
   \pt_y
   {\up{n+1}{\omega}} (x,y)
 \right|
 \leq
 \psi (x) \; V(x) \, , \quad \omega = r,s,p,q,z.
$$
 Thus,
$$
\Phi (x) = \psi (x) \; V(x) \, . \ \square
$$



\section{Existence and uniqueness of a solution}

 The uniqueness of a solution of the Cauchy problem
 (3.32), $(3.32^0)$
 follows the uniqueness theorem
 (\cite{RYa}, Chapter~1, {\S}~8, Subsection~2).

 In this section, as far as in the Section~2, we consider
 the Cauchy problem for the general system
 (2.8), $(2.8^0)$
 with respect to the unknown vector-function $u=(u_1, \dots , u_m)$.
 Suppose successive approximations
$\{{\up{n}{u}}(x,y)\}$
 and their derivatives
$\{{(\up{n}{u})}_y(x,y)\}$
 are uniformly bounded on the compactum
 (4.6)
 for an arbitrary point
 ${(\bar x, \bar y)\in[0,+\infty)\times \mathbb{R}}$.
 All considerations are being made over this compactum.
 We follow here the standard scheme of the proof
 from \cite{RYa}, Chapter~1, {\S}~8, Subsection~2.
 Also \cite{Tunitsky-dis} was used.


   Suppose $\varphi (u)\in C^1$.
   By definition, put
$\tilde \varphi (\lambda )=\varphi (\bar u+\lambda (u-\bar u)),
\lambda \in \mathbb{R}$.
 By the Newton---Leibniz formula,
 we have
 \cite{Tunitsky-dis}
$$
  \tilde \varphi (1)-\tilde \varphi (0)
  \! = \!
  \varphi (u)-\varphi (\bar u)
  \! =
  \!\!
  {\ds \int_0^1}
  \!\!\!
  \tilde \varphi _{\lambda}d\lambda =
  \!\!
  {\ds \int_0^1
  \!\!
  \left(
     \sum\lms_{j=1}^{m}\frac{\pt \varphi}{\pt u_j}
     (\bar u+\lambda (u-\bar u))(u_j-\bar u_j)
     \!\!
  \right)
  \!\!
   d\lambda.
}
$$
 Finally, we obtain
$$
  \varphi (u)-\varphi (\bar u)=
{\ds
   \sum\lms_{j=1}^m(u_j-\bar u_j)
   \int_0^1 \frac{\pt \varphi}{\pt u_j} (\bar u+\lambda (u-\bar u))
   \; d\lambda .
}
 \eqno(6.1)
$$


\subsection{Continuity of the solution}

{\bf Lemma~6.1.}\ The vector-function
$\lim\lms_{n\to\infty}{\up{n}{u}}(x,y)$
 is continuous.
\\
 {\bf Proof.}\
 A sufficient condition on convergence of
 a functional sequence to a continuous function
 is given by
 the classical theorem:
 \textit{elements of a sequence must be
 continuous functions,
 and the sequence must be uniformly convergent.}

 Uniform convergence of the sequence
 $\{\up{n}{u}\}$
 follows
 uniform convergence of the series
 ${\ds\sum_{n=0}^{\infty}}(\up{n+1}{u}-\up{n}{u})(x,y)$.

 Let
 ${\up{n+1}{u}}$,
 ${\up{n}{u}}$
 be successful approximations.
 By (4.1), we have
$$
\begin{array}{l}
    ( \pt_x
           +
           \xi_i(x,y,{\up{n+1}{u}})
           \pt_y
    )
    {\up{n+1}{u_i}}
    =
    f_i(x,y,{\up{n}{u}})
     \; ,
    \smallskip \\
    ( \pt_x
           +
           \xi_i(x,y,{\up{n+1}{u}})
           \pt_y
    )
    {\up{n}{u}}_i
    +
    \xi_i(x,y,{\up{n}{u}})
    {\up{n}{u}}_{iy}
    =
    f_i(x,y,{\up{n-1}{u}})
    +
    \xi_i(x,y,{\up{n+1}{u}})
    {\up{n}{u}}_{iy}
      \; ,
    \smallskip \\
    {\up{n+1}{u}}
    (0,y)
    =
    {\up{n}{u}}
    (0,y)
    =
    u^0(y)
    \, .
\end{array}
$$
 By subtracting,
 with respect to
 ${{\up{n+1}{u}} - {\up{n}{u}}}$
 we get
 the Cauchy problem
$$
\left\{
\begin{array}{ll}
    (
           \pt_x
           +
           \xi_i(x,y,{\up{n+1}{u}})
           \pt_y
    )
    (
      {{\up{n+1}{u_i}} - {\up{n}{u}}_i}
    )
    \;
    =
    &
    f_i(x,y, {\up{n}{u}})
    -
    f_i(x,y, {\up{n-1}{u}})
    -
       \smallskip \\
     &
    -
    {\up{n}{u}}_{iy}
    (
        \xi_i(x,y,{\up{n+1}{u}})
        -
        \xi_i(x,y,{\up{n}{u}})
    )
    \; ,
    \smallskip \\
    ({\up{n+1}{u}} - {\up{n}{u}})
    (0,y)
    =0
    \, .
\end{array}
\right.
$$
 By (6.1), with respect to
 ${\up{n+1}{r}} = {\up{n+1}{u}} - {\up{n}{u}}$
 we get
$$
\left\{
\begin{array}{l}
    ( \pt_x
           +
           \xi_i(x,y,\up{n+1}{u})
           \pt_y
    )
    \up{n+1}{r_i}
    =
    \smallskip
    \\
{\ds
   =
   \sum\lms_{j=1}^m
   {\up{n}{r}}_j
   \int_0^1
\left(
   \frac{\pt f_i}{\pt u_j}
   (
     x,y,
     \up{n-1}{u}
     +
     \lambda
     \up{n}{r}
   ) -
   {\up{n}{u}}_{iy}
   \frac{\pt \xi_i}{\pt u_j}
   (
       x,y,
       {\up{n}{u}}
       +
       \lambda
       {\up{n+1}{r}}
   )
\right)
   d\lambda
}
     \; ,
     \smallskip \\
     {\up{n+1}{r}}
     (0,y)
     =
     0
     \, .
\end{array}
\right.
$$
 By (2.10), we get
$$
    |\up{n+1}{r_i}|\le
{\ds
   \int_0^x
   \max\lms_k|{\up{n}{r}}_k|
   \sum\lms_{j=1}^m
   \int_0^1
\left(
   \left|\frac{\pt f_i}{\pt u_j}\right| +
   |{\up{n}{u}}_{iy}|
   \left|\frac{\pt \xi_i}{\pt u_j}\right|
\right)
   d\lambda \; d\tau.
}
 \eqno(6.2)
$$
 From (3.32), taking into account
 {$C^1\mbox{-boundedness}$}
 of vector-functions $\rho,\sigma,\pi,\kappa$
 and uniform boundedness of
 $\{\up{n}{u}\}$,
 we obtain
 uniform boundedness of
 ${\ds \left|\frac{\pt f_i}{\pt u_j}\right|}$.
 ${\ds \left|\frac{\pt \xi_i}{\pt u_j}\right|}$
 is $0$ or $1$,
 and
 $\{{\up{n}{u}}_{yi}\}$
 is uniformly bounded,
 therefore,
 we have
$$
{\ds
   \sum\lms_{j=1}^m
   \int_0^1
\left(
   \left|\frac{\pt f_i}{\pt u_j}\right| +
   |{\up{n}{u}}_{iy}|
   \left|\frac{\pt \xi_i}{\pt u_j}\right|
\right)
   d\lambda
   \le C\; ,
}
$$
 where $C$ is some constant.
 So the left part is uniformly bounded.
 By definition, put
$$
  R_n(x)=\max\lms_i\sup\lms_{(\tau,\, y)\in [0,x]\times \mathbb{R}}
  |{\up{n}{r}}_i(\tau ,y)|.
$$
 Then (6.2) is
$$
  R_{n+1}(x)\le C\int_0^xR_n(\tau )\; d\tau .
   \eqno(6.3)
$$
 By recursive applying (6.3),
 we have
$$
  R_{n+1}(x)
  \le
  \max\lms_{\tau\in [0,x]}R_1(\tau )\frac{(Cx)^n}{n!}
  \le
  \max\lms_{\tau\in [0,\bar x]}R_1(\tau )\frac{(C \bar x)^n}{n!}\; ,
$$
 therefore, by the Weierstrass criterion,
 we obtain uniform convergence of
 ${\ds\sum_{n=0}^{\infty}}({\up{n+1}{u}}-{\up{n}{u}})(x,y)$
 in
 ${[0,\bar x]\times \mathbb{R}}$.
 Thus
 $\lim\lms_{n\to\infty}{\up{n}{u}}(x,y)$
 is continuous
 in
 ${[0,\bar x]\times \mathbb{R}.}\ \Box$


\subsection{Continuous differentiability of the solution}

{\bf Lemma~6.2.}\ The vector-function
 $\lim\lms_{n\to\infty}{\up{n}{u}}(x,y)$
 is continuously differentiable.
\\
{\bf Proof.}\ To prove continuity of
 $\lim\lms_{n\to\infty}{{(\up{n}{u})}_y}$,
 we get the Arzela theorem.
 Uniform boundedness of
 $\{{\up{n}{u}}_y\}$
 is proved already.
 Now we shall prove
 equicontinuity of
 $\{{\up{n}{u}}_y\}$.

 First,
 we shall prove
 equicontinuity of
 $\{\up{n}{u}\},\{\up{n}{g}\}$.
 By
(6.1),
 we have
$$
 \!
 \begin{array}{ll}
  {\up{n}{u}}(x_1,y_1)-{\up{n}{u}}(x_2,y_2)&
  \!\!\!\!
  =(x_1-x_2){\ds\int_0^1}
  {\up{n}{u}}_x(x_2+\lambda (x_1-x_2),y_2+\lambda (y_1-y_2))
  \; d\lambda +
 \medskip   \\
  &
  \!\!
  +(y_1-y_2){\ds\int_0^1}
  {\up{n}{u}}_y(x_2+\lambda (x_1-x_2),y_2+\lambda (y_1-y_2))
  \; d\lambda.
\end{array}
$$
 From (4.1) and from uniform boundedness of
 $\{\up{n}{u}\},\, \{{\up{n}{u}}_y\}$,
 we get
 uniform boundedness of
 $\{{\up{n}{u}}_x\}$.
 Finally, we obtain
 equicontinuity of
 $\{\up{n}{u}\}$.

 By (6.1),
 with respect to the vector-function
 ${\up{n}{g}}(\tau ,x,y)$,
 we have
$$
\begin{array}{l}
  {\up{n}{g}}(\tau ,x_1,y_1)-{\up{n}{g}}(\tau ,x_2,y_2)=
                                \medskip   \\
  =(x_1-x_2){\ds\int_0^1}{\up{n}{g}}_x(
  \tau ,x_2+\lambda (x_1-x_2),y_2+\lambda (y_1-y_2))
  \; d\lambda +
                                \medskip   \\
  +\:(y_1-y_2){\ds\int_0^1}{\up{n}{g}}_y(
  \tau ,x_2+\lambda (x_1-x_2),y_2+\lambda (y_1-y_2))
  \; d\lambda.
\end{array}
$$
 Uniform boundedness of
 $\{{\up{n}{g}}_x\},\{{\up{n}{g}}_y\}$
 follows (2.11)
 and
 uniform boundedness of
 $\{\up{n}{u}\},\{{\up{n}{u}}_y\}$
 (here
  ${\xi_{yi}=0,}$
 ${\pt\xi_i/\pt{\up{n}{u}}_j}$
 is 0 or 1).
 Therefore,
$$
  |{\up{n}{g}}_i(\tau ,x_1,y_1)-{\up{n}{g}}_i(\tau ,x_2,y_2)|
  \le \mbox{\rm const}\,(|x_1-x_2|+|y_1-y_2|),
$$
 i.e.
 $\{\up{n}{g}\}$
 is equicontinuous.

 Consider the function
 $\up{n+1}{u_{iy}}$
 on the compactum
 $G(\bar x, \bar y)$:
$$
\begin{array}{l}
  {\up{n+1}{u_{iy}}}(x,y)
  =
  u^0_{iy}
  ({\up{n+1}{g_i}}(0,x,y))
  +
 \medskip   \\
  \phantom{123}
  +
  {\ds\int_0^x}
  \left\{
    -
    {\up{n+1}{u_{iy}}}
    \left(
        \xi_{iy}
        +
        \sum\lms_{j}
            \frac{\ds \pt \xi_{i}}
                 {\ds \pt {\up{n+1}{u_j}}}
            \cdot
            {\up{n+1}{u_{jy}}}
    \right)
    +
    \left(
        f_{iy}
        +
        \sum\lms_{j}
            \frac{\ds \pt f_{i}}
                 {\ds \pt {\up{n}{u}_j}}
            \cdot
            {\up{n}{u}_{jy}}
    \right)
  \right\}
 \medskip
 \\
  \phantom{123}
  \left(
     \tau ,
     {\up{n+1}{g_i}}
     (\tau ,x,y),
     {\up{n}{u}}
     (\tau ,{\up{n+1}{g_i}}(\tau ,x,y))
   \right)
   d\tau
   \, .
\end{array}
 \eqno(6.4)
$$
 Here
 ${\xi_{iy} = 0}$,
 ${
            \pt \xi_i/
                 \pt {\up{n+1}{u_j}}
 }$
 is $0$ or $1$.
 $\{\up{n}{u}\},\{\up{n}{g}\}$
 are equicontinuous,
 therefore
 we have equicontinuity of
$$
\begin{array}{l}
\{u^0_{iy}({\up{n+1}{g_i}}(0,x,y))\}
                          \: ,     \medskip   \\
\left\{
    f_{iy}
  \left(
     \tau ,{\up{n+1}{g_i}}(\tau ,x,y),
     {\up{n}{u}}(\tau ,{\up{n+1}{g_i}}(\tau ,x,y))
   \right)
\right\}
                          \: ,     \medskip   \\
\left\{
            \frac{\ds \pt f_{i}}
                 {\ds \pt {\up{n}{u}_j}}
  \left(
     \tau ,{\up{n+1}{g_i}}(\tau ,x,y),
     {\up{n}{u}}(\tau ,{\up{n+1}{g_i}}(\tau ,x,y))
   \right)
\right\}
 \: .
\end{array}
$$
 We shall use the Cantor theorem: \textit{if
 a function is continuous on a compactum from
 $\mathbb{R}^s$,
 then
 the function is equicontinuous on this compactum}.
 Taking into account
 uniform boundedness of
 $\{{\up{n}{u}}_y\}$
 and the Cantor theorem,
 we obtain
 equicontinuity of right parts of (6.4),
 i.e. of
 $\{{\up{n}{u}}_y\}$,
 on the compactum
 $G(\bar x,\bar y)$.

 We shall use the Arzela theorem: \textit{if a functional sequense
 is uniformly bounded and equicontinuous on a compactum,
 then there exists a uniformly convergent subsequence on the compactum.}
 Therefore,
 there exists
 a uniformly convergent subsequence
 $\{{\up{n_k}{u_y}}\}$.
 By the theorem on termwise differentiating of a functional sequence,
 the vector-function
 $\lim\lms_{k\to\infty}{\up{n_k}{u}}(x,y)$
 is continuously differentiable with respect to $y$,
 and
$$
  \pt_y \lim\lms_{k\to\infty}{\up{n_k}{u}}=
  \lim\lms_{k\to\infty}{\up{n_k}{u_y}}.                     \eqno(6.5)
$$
 Taking into account
 the previously proved uniform convergence
 of
 $\{\up{n}{u}\}$,
 we have
$$
  \lim\lms_{k\to\infty}{\up{n_k}{u}}=
  \lim\lms_{n\to\infty}{\up{n}{u}}.                          \eqno(6.6)
$$
 By (6.6),
 we get from (6.5) the next rule
 for differentiating
 the vector-function
 $\lim\lms_{n\to\infty}{\up{n}{u}}$
 with respect to $y$:
$$
  \pt_y \lim\lms_{n\to\infty}{\up{n}{u}}=
  \lim\lms_{k\to\infty}{\up{n_k}{u_y}}.                      \eqno(6.7)
$$
 Vector-functions
 $\up{n_k}{u_y}$
 are continuous,
 therefore
 the vector-function
 $\pt_y \lim\lms_{n\to\infty}{\up{n}{u}}$
 is continuous.

 By passage to the limit,
 we get from (4.1) the next rule
 for differentiating
 the vector-function
 $\lim\lms_{n\to\infty}{\up{n}{u}}(x,y)$
 with respect to $x$:
$$
  \pt_x \lim\lms_{n\to\infty}{\up{n}{u}}_i(x,y)=
  f_i(x,y,\lim\lms_{n\to\infty}{\up{n}{u}})-
  \xi_i(x,y,\lim\lms_{n\to\infty}{\up{n}{u}})
  (\pt_y \lim\lms_{n\to\infty}{\up{n}{u}}_i).      \eqno(6.8)
$$
 By continuity of vector-functions
 $f,\; \xi ,\;\lim\lms_{n\to\infty}{\up{n}{u}},\;
 \pt_y \lim\lms_{n\to\infty}{\up{n}{u}}$,
 we get continuity of the derivative
 with respect to $x$.

 Now we proved that
 ${\lim\lms_{n\to\infty}{\up{n}{u}}(x,y)
 ~\in~C^1(G(\bar x,\bar y))}$.

 Put the limit function into the Cauchy problem
 (2.8), $(2.8^0)$.
 Taking into account (6.8),
 we see that the limit function is a solution of the Cauchy problem.

 It remains to note that for
 ${\forall\ (x,y)\in [0,+\infty)\times \mathbb{R}}$
 we can get some compactum
 ${G(\bar x,\bar y) \ni (x,y)}$.
 Consider an intersection of two such compact sets.
 Initial approximations are coincide
 in the intersection,
 and they are equal to the initial approximation
 for the half-plane
 ${[0,+\infty)\times \mathbb{R}}$.
 Therefore,
 limit functions are coincide too,
 so they are restrictions of the limit function
 for the half-plane
 ${[0,+\infty)\times \mathbb{R}}$.

 The Arzela theorem is a pure existence theorem.
 Note that we use the Arzela theorem
 in the proof,
 but we don't use this theorem for constructing the solution.
 $\Box$



\section{The main result}
\subsection{Theorem of the existence and the uniqueness}

 If
 (4.4), (4.7), (4.8)
 are satisfied,
 then there exists a unique
 ${C^1\mbox{-solution}}$
 of
 (3.32), $(3.32^0)$
 in
 the half-plane ${x \geq 0}$.
 It was proved above.
 Now we formulate conditions on
 coefficients of the equation (3.1)
 and on initial functions.
 Under these conditions
 the inequalities
 (4.4), (4.7), (4.8)
 are satisfied.

 Let $M_1, M_2$ be arbitrary positive constants,
 and let
 $\varepsilon , \delta$
 be constants under conditions
$$
0 < \varepsilon \leq \delta ,
 \quad (\delta - \varepsilon)/2 < 1 \, ;
 \eqno (7.1)
$$
 here $\eta (x) \in C^0(\mathbb{R})$
 is an arbitrary nonnegative function.

 By definition, put
$$
\begin{array}{l}
 N_1 = \max \{
               M_1, \frac{1}{2} M_2 (4M_1 + 9 M_1^2)
            \},
 \smallskip
 \\
    N_2 = M_1 M_2 \, ,
 \medskip
 \\
    \tilde{\eta}(x) = \frac{\ds 1}{\ds 1 + 2M_1 |x|} \; \eta (x) \, .
\end{array}
\eqno (7.2)
$$
 Suppose
$$
 \begin{array}{lll}
   1) & \mbox{коэффициенты}\ A,B,C,D \in C^2(\mathbb{R}^5)\, ;
      & \phantom{1234567890}
   \smallskip
   \\
   2)
      & \mbox{начальные функции} \ z^0 \in C^3(\mathbb{R}),
        \ p^0 \in C^2(\mathbb{R});
      &
   \medskip
   \\
   3) & |a|\leq M_1, \ \ 1/\Delta \leq M_2 \, ,
      &
   \medskip
   \\
      & \left|
          \frac{\ds \pt a}{\ds \pt \omega}
        \right|
        \leq M_1 \eta (x) \, ,
        \qquad
        \left|
          \frac{\ds \pt a}{\ds \pt z}
        \right|
        \leq M_1 \tilde{\eta} (x) \, ,
      &
   \smallskip
   \\
      & \mbox{где} \ \ a=B,C,D,\Delta, \quad \omega = x,y,p,q;
      &
   \smallskip
   \\
   4) & (6 N_1 + 8 N_2)
        {\ds\int_{-\infty}^{+\infty}}
          \eta (x)
        \; dx
        \leq
        (\delta - \varepsilon)/2
        \, ;
      &
   \medskip
   \\
   5) & |r^0|, |s^0|
        \le
        1 -
        (\delta - \varepsilon)/2
        \, , \quad
        |z_y^0|, |p^0| \le 1
        \, ;
      &
   \medskip
   \\
   6) & \inf\lms_{y \in \mathbb{R}} r^0 (y) -
        \sup\lms_{y \in \mathbb{R}} s^0 (y) \geq \delta > 0 \, .
\end{array}
\eqno (7.3)
$$

{\bf Theorem~7.1.}
 Under conditions (7.3)
 there exists a unique
 $C^3\mbox{-smooth}$
 solution of the Cauchy problem
 (3.1), (3.3).

 {\bf Proof.}
 First, we shall prove that
 conditions
 (4.4), (4.7), (4.8)
 (and conditions before them
 in the beginning of Subsection 4.2)
 follow conditions (7.3).

 In conditions (4.4)
 we suppose
 ${C^1\mbox{-smoothness}}$
 and boundedness
 of
 $r^0, s^0$.
 Now we shall prove that
 the condition
 ${r^0, s^0 \in C^1}$
 follows (7.3).
 From the first and the second equations of
 $(3.32^0)$,
 we have
$$
 {\ds
    r^0 (y) - s^0 (y)
    =
    \frac{\Delta
        (0,y,z^0(y),p^0(y),z_y^0(y))
    }
    {z_{yy}^0 (y) + B
        (0,y,z^0(y),p^0(y),z_y^0(y))
    }
 }
 \, ,
$$
 hence,
$$
 z_{yy}^0 (y) + B
 (0,y,z^0(y),p^0(y),z_y^0(y))
  =
 \frac{\Delta
     (0,y,z^0(y),p^0(y),z_y^0(y))
 }
 {r^0 (y) - s^0 (y)}
 \, .
$$
 Therefore, by (7.3), point 6),
 and by the condition
  ${1/\Delta \leq M_2}$
  from (7.3), point 3),
  we get
$$
 z_{yy}^0 (y) + B
 (0,y,z^0(y),p^0(y),z_y^0(y))
 \neq
 0.
$$

 Then, by $(3.32^0)$
 and by
 (7.3), points 1), 2),
 we get
 ${r^0, s^0 \in C^1}$.
 By
 (7.3), point 5),
 we have
 boundedness of
 $r^0, s^0$.

 The conditions (4.4)
 are formulated under
 assumptions of
 ${C^1\mbox{-smoothness}}$
 and boundedness of
 $p^0, q^0$.
 By $(3.32^0)$,
 we have
 ${q^0 = z_y^0}$.
 By (7.3), point 2),
 we get
 $C^1\mbox{-smoothness}$ of $p^0, q^0$.
 By (7.3), point 5),
 we have boundedness of
 $p^0, q^0$.

 Further, the conditions (4.4)
 are formulated under
 assumptions of
 ${C^1\mbox{-smoothness}}$
 and boundedness of
 vector-functions
 $\rho, \sigma, \pi, \kappa$.
 By (3.32),
 these vector-functions depend on
 ${1/\Delta}$
 and on first derivatives of
 ${B, C, D, \Delta}$.
 By (7.3), point 3),
 $1/\Delta$
 is bounded away from zero
 by some positive constant.
 Then, by (7.3), point 1),
 from ${C^2\mbox{-smoothness}}$ of
 ${A, B, C, D}$,
 we have
 ${C^1\mbox{-smoothness}}$
 of
 vector-functions
 $\rho, \sigma, \pi, \kappa$.
 By (7.3), point 3),
 we have boundedness of
 vector-functions
 $\rho, \sigma, \pi, \kappa$.

 By (7.3), point 3), we obtain
$$
 \begin{array}{l}
   |\rho_i|, |\sigma_i|\leq M_1 \eta \, , \quad i=0,13;
 \smallskip
 \\
   |\rho_i|, |\sigma_i|\leq \frac{1}{2} M_2 (4M_1 + 9 M_1^2) \eta \, , \quad i=1,2,7,8;
 \smallskip
 \\
   |\rho_i|, |\sigma_i|\leq M_1 M_2 \tilde{\eta} \, , \quad i=3,4,5,6,9,10,11,12;
 \smallskip
 \\
   |\pi_i|, |\kappa_i|\leq M_1 \, , \quad i=0,1;
 \end{array}
$$
 Hence, by (4.3), (7.2), we have
$$
 \begin{array}{l}
   \alpha_1 (x) \leq \max
     \{
       M_1, \frac{1}{2} M_2 (4M_1 + 9 M_1^2)
     \}
     \eta (x) = N_1 \eta (x) \, ,
 \smallskip
 \\
   \alpha_2 (x) \leq M_1 M_2 \tilde{\eta} (x) =
      N_2 \tilde{\eta} (x) \, ,
 \smallskip
 \\
   \alpha_3 \leq M_1 \, .
 \end{array}
 \eqno (7.4)
$$
 Then, by points 4), 5) of (7.3),
 we satisfy the second inequality (4.4):
$$
 \begin{array}{l}
   U_0 +
        6
        {\ds\int_{-\infty}^{+\infty}}
            \alpha_1 (x)
        \; dx
        +
        8
        {\ds\int_{-\infty}^{+\infty}}
        (1 + 2\alpha_3 |x|)
        \; \alpha_2 (x)
        \; dx
        \le
 \medskip
 \\
   \le
        U_0 +
        (6 N_1 + 8 N_2)
        {\ds\int_{-\infty}^{+\infty}}
            \eta (x)
        \; dx
        \le
        (
            1
            -
            (\delta - \varepsilon)/2
        )
        +
            (\delta - \varepsilon)/2
        =
        1
        \, .
 \end{array}
$$

 The first inequality (4.4) follows (7.3), point 5).
 (4.7) is (7.3), point 6).
 Taking into account (7.2) and (7.4),
 we see that
 (4.8) is (7.3), point 4).

 Thus,
 (4.4), (4.7), (4.8)
 follows (7.3).
 By
 (4.4), (4.7), (4.8),
 there exists a unique
 $C^1\mbox{-solution}$
 of (3.32), $(3.32^0)$
 in the half-plane
 ${x \ge 0}$.

 In the half-plane
 ${x \le 0}$
 a solution is constructed analogously.
 By the same initial functions,
 we obtain
 $C^1\mbox{-smoothness}$ of a solution
 on the whole plain.

 By
 ${1/\Delta \leq M_2}$
 (it is (7.3), point 3)),
 we have
    $
    {\Delta > 0}.
    $
 By Theorem~3.2,
 it is sufficient for corresponding
 a $C^1\mbox{-solution}$ of (3.32), $(3.32^0)$
 to
 a {$C^3\mbox{-solution}$} of (3.1), (3.3).
 $\square$

 {\bf Example~7.1.}
 Conditions (7.3) define nonempty set of
 coefficients and initial data.
 Indeed, we take
$
 {A=1/16,} \quad
 {B=1/2,} \quad
 {C=0,} \quad {D=0,} \quad
 {z^0 = 1,} \quad {p^0 = 1.}$ Then, by (3.2), we have
 ${\Delta^2 = 4A,}$
 so
 ${\Delta = 1/2.}$
 $(3.32^0)$ is
$
 {r^0, s^0 = (C \pm \Delta)/2B,}
$
 therefore,
$
 {r^0, s^0 = \pm 1/8.}
$
 Thus, constants
$
 {M_1 = 1/2,} \quad
 {M_2 = 2}\, , \quad
 {\delta = 1/4,} \quad
 {\varepsilon = 1/4}
$
 and functions
 ${\eta(x) = 0}, \quad {\tilde \eta(x) = 0}$
 satisfy (7.3).
 Constants
 $N_1,\ N_2$ are defined by (7.2). $\Box$

{\bf Example~7.2.}
 \textit{Equations with constant coefficients.}
 (3.32) is reduced to its independent subsystem (3.33):
$$
 \begin{array}{l}
 ( \partial_x
   +
   s \, \partial_y
 )
 \: r = 0 \, ,
 \medskip
 \\
 ( \partial_x +
    r \, \partial_y
 )
 \: s = 0 \, .
 \end{array}
$$
 For this system
 conditions (7.3)
 are conditions on initial functions
 $r^0,\ s^0$.
 We get $r,\ s$
 by solving this Cauchy problem,
 and then we solve linear equations
 (3.32)
 with respect to
 $p, \ q, \ z$.
 By $(3.32^{00})$,
 initial functions $p^0,\ q^0,\ z^0$
 are functions of
 $r^0,\ s^0$.
 $\Box$


\subsection{Initial conditions}

 We set some sufficient conditions such that
 conditions on initial functions in (7.3) are satisfied.

Let
$$
 {m_1 \geq 0}, \quad {m_2 \geq 0}, \quad
 {L_1 > 0}, \quad {L_2 >0}, \quad {L_3 > 0}
$$
 be arbitrary constants.
 By definition, put constants
$$
\begin{array}{l}
 m_3 = \frac{\ds m_1}{\ds L_1 - m_1} \: (1 + L_3 + \frac{\ds m_2}{\ds 2L_1}) +
       \frac{\ds m_2}{\ds 2L_1} \, ,
 \medskip
 \\
 m_4 = \frac{\ds m_1 (2 + L_3)}{\ds L_1 - m_1} +
       \frac{\ds m_2}{\ds 2(L_1 - m_1)} \, ,
 \medskip
 \\
 L_4 = \frac{\ds 1}{\ds 2 M_2 L_2} \, ,
\end{array}
\eqno(7.5)
$$
 where $M_2$ is the constant from (7.3), point 3).

 Suppose
$$
 \begin{array}{rlrl}
   1)
      & z^0 \in C^3(\mathbb{R}),
        \ p^0 \in C^2(\mathbb{R}) \, ;
        \phantom{1234567}
 &  6) & |C| \leq m_2 \, ;
   \smallskip
   \\
   2) & |z_y^0| \leq 1, \ \ |p^0| \leq 1 \, ;
 &   7) & \frac{\ds |\Delta|}{\ds 2L_1} \le L_3 \, ;
   \medskip
   \\
   3) & |z_{yy}^0| \leq m_1, \quad |p_y^0| \leq m_1 \, ;
 &   8) & L_3 + m_3 \le 1 - \frac{\ds \delta - \varepsilon}{\ds 2} \, ;
   \medskip
   \\
   4) & 0 < L_1 \le B \le L_2 \, ;
 &   9) & L_4 - m_4 \ge \delta/2 > 0 \, .
   \medskip
   \\
   5) & L_1 - m_1 > 0 \, ;
\end{array}
\eqno (7.6)
$$

{\bf Theorem~7.2.} Suppose coefficients $A, B, C, D$
 of the equation
 (3.1) satisfy (7.3),
 and, moreover,
 conditions (7.6) are satisfied.
 Then $z^0,\ p^0$, and
 functions $r^0,\ s^0$
 as functions of $z^0,\ p^0$
 (see $(3.32^0)$)
 satisfy (7.3).

{\bf Proof.}
 (7.3), point 2) follows (7.6), point 1).
 (7.3), point 5) for functions $z_y^0, \ p^0$
 is (7.6), point 2).

 Now we shall prove that for $r^0, \ s^0$
 the condition
 (7.3), point 5)
 is satisfied.
 We have
$$
\begin{array}{l}
\ds
 \frac{C \pm \Delta - 2 p_y^0}{2(z_{yy}^0 + B)} -
 \frac{C \pm \Delta}{2 B}
 =
 \frac{ (C \pm \Delta)B - 2p_y^0 B - (C \pm \Delta)(z_{yy}^0 + B)}
 {2(z_{yy}^0 + B)B}
 =
 \bigskip
 \\
\ds
 =
 \frac{-2 p_y^0 B - (C \pm \Delta)z_{yy}^0}
 {2(z_{yy}^0 + B)B}
 =
 - \frac{p_y^0}{z_{yy}^0 + B} -
 \frac{(C \pm \Delta)z_{yy}^0}{2(z_{yy}^0 + B)B} \, ,
\end{array}
$$
 therefore,
$$
 \frac{C \pm \Delta - 2 p_y^0}{2(z_{yy}^0 + B)}
 =
 \frac{C \pm \Delta}{2 B} -
 \frac{p_y^0}{z_{yy}^0 + B} -
 \frac{(C \pm \Delta)z_{yy}^0}{2(z_{yy}^0 + B)B} \, .
\eqno(7.7)
$$
 By definition, put the index ${\omega = r,s.}$
 By $(3.32^0)$, (7.7),
 we get
$$
\begin{array}{l}
\ds
 |\omega^0|
 =
 \left|
 \frac{C \pm \Delta - 2 p_y^0}{2(z_{yy}^0 + B)}
 \right|
 \le
 \frac{|C| + |\Delta|}{2|B|} +
 \frac{|p_y^0|}{|z_{yy}^0 + B|} +
 \frac{(|C| + |\Delta|)|z_{yy}^0|}{2|z_{yy}^0 + B|\cdot |B|}
 \le
 \bigskip
 \\
\ds
 \up{(7.6),\ \mbox{{\scriptsize p.}}\, 3)-6)}{\le}
 \frac{|\Delta| + m_2}{2L_1} +
 \frac{1}{L_1 - m_1} \:
 \left(
 m_1 +
 \frac{m_1(|\Delta| + m_2)}{2L_1}
 \right)
 =
 \bigskip
 \\
\ds
 =
 \frac{|\Delta|}{2L_1} +
 \frac{m_2}{2L_1} +
 \frac{m_1}{L_1 - m_1} \:
 \left(
 1 +
 \frac{|\Delta|}{2L_1} +
 \frac{m_2}{2L_1}
 \right)
 \up{(7.6),\ \mbox{{\scriptsize p.}}\, 7)}{\le}
 \bigskip
 \\
\ds
 \le
 L_3 +
 \frac{m_2}{2L_1} +
 \frac{m_1}{L_1 - m_1} \:
 \left(
 1 +
 L_3 +
 \frac{m_2}{2L_1}
 \right)
 \up{(7.5)}{=}
\ds
 L_3 + m_3
 \up{(7.6),\ \mbox{{\scriptsize p.}}\, 8)}{\le}
 1 - \frac{\delta - \varepsilon}{2} \, .
\end{array}
$$
 Now we shall prove that (7.3), point 6) is satisfied.
 We have
$$
\begin{array}{l}
\ds
 \frac{C \pm \Delta - 2 p_y^0}{2(z_{yy}^0 + B)} -
 \frac{\pm \Delta}{2 B}
 =
 \frac{(C \pm \Delta)B - 2p_y^0 B - (\pm \Delta)(z_{yy}^0 + B)}
 {2(z_{yy}^0 + B)B}
 =
 \bigskip
 \\
\ds
 =
 \frac{CB - 2 p_y^0 B - (\pm \Delta)z_{yy}^0}
 {2(z_{yy}^0 + B)B}
 =
 \frac{C - 2 p_y^0}{2(z_{yy}^0 + B)} -
 \frac{(\pm \Delta)}{2B} \cdot
 \frac{z_{yy}^0}{z_{yy}^0 + B}
 \, ,
\end{array}
$$
 therefore,
$$
 \frac{C \pm \Delta - 2 p_y^0}{2(z_{yy}^0 + B)}
 =
 \frac{\pm \Delta}{2 B}
 -
 \frac{(\pm \Delta)}{2B} \cdot
 \frac{z_{yy}^0}{z_{yy}^0 + B}
 +
 \frac{C - 2 p_y^0}{2(z_{yy}^0 + B)}
 \, .
\eqno(7.8)
$$
 By definition, put the index ${\omega = r,s}$.
 By $(3.32^0)$, (7.8),
 we get
$$
\begin{array}{l}
\ds
 |\omega^0|
 =
 \left|
 \frac{C \pm \Delta - 2 p_y^0}{2(z_{yy}^0 + B)}
 \right|
 \ge
 \bigskip
 \\
 \ds
 \ge
 \inf
 \left|
   \frac{\pm \Delta}{2B}
 \right|
 -
 \sup
 \left|
   \frac{\pm \Delta}{2B} \cdot
   \frac{z_{yy}^0}{z_{yy}^0 + B}
 \right|
 -
 \sup
 \left|
   \frac{C - 2 p_y^0}{2(z_{yy}^0 + B)}
 \right|
 \up{(7.6),\ \mbox{{\scriptsize p.}}\, 4)}{\ge}
 \bigskip
 \\
 \ds
 \up{(7.3),\ \mbox{{\scriptsize p.}}\, 3)}{\ge}
 \frac{1}{M_2} \cdot \frac{1}{2L_2}
 -
 \sup
 \left|
   \frac{\pm \Delta}{2B} \cdot
   \frac{z_{yy}^0}{z_{yy}^0 + B}
 \right|
 -
 \sup
 \left|
   \frac{C - 2 p_y^0}{2(z_{yy}^0 + B)}
 \right|
 \ge
 \bigskip
 \\
 \ds
 \up{(7.6),\ \mbox{{\scriptsize p.}}\, 3),5),7)}{\ge}
 \frac{1}{M_2} \cdot \frac{1}{2L_2}
 -
 L_3 \cdot \frac{m_1}{L_1 - m_1}
 -
 \frac{m_2 + 2m_1}{2(L_1 - m_1)}
 \up{(7.5)}{\ge}
 \bigskip
 \\
 \ds
 \ge
 L_4 - m_4
 \up{(7.6),\ \mbox{{\scriptsize p.}}\, 9)}{\ge}
 \delta/2 > 0 \, .
\end{array}
$$
 The sign of $\omega^0$ is defined by its principal part
 ${\pm \Delta/2B,}$
 hence,
 $r^0$ and $s^0$ have different signs.
 Moreover, from ${B > 0}$ ((7.6), point 4)) and from
 ${\Delta > 0}$
 ((7.3), point 3)),
 we get
 ${r^0 > 0}, \quad {s^0 < 0}$.
 Then,
 from the proved estimate
 ${|\omega^0| \ge \delta/2}$,
 we obtain
 $
 {\inf\lms_{y \in \mathbb{R}} r^0 (y)
 -
 \sup\lms_{y \in \mathbb{R}} s^0 (y) \geq \delta > 0}
 $.
 $\Box$

{\bf Example~7.3.} Conditions (7.6)
 are satisfied by coefficients
 and initial functions
 of Example~7.1.
 In this case, we have
 ${M_2 = 2}$,
 ${m_i = 0}$, ${L_i = 1/2}$, ${i = 1,2,3,4}$.
 Setting small perturbations of
 coefficients and of initial functions from
 Example~7.1,
 we obtain an example of an equation
 with nonconstant coefficients, dependent on
 $x,y,z,z_x,z_y$,
 and with nonconstant initial functions. $\Box$

\section{Supplement. Contact approach}

\subsection{A contact transformation can transform
  a classical solution into a solution which is singular
  at each point}

 Consider the next example.
 The author is grateful to V.~V.~Lychagin
 \cite{Ly75}, \cite{Ly79}
 and to L.~V.~Zilbergleit
 for acquainting foundations on this theme.

 {\bf Example~8.1.}
 Consider the
 $\mbox{J}^1 (\mathbb{R}^2)$
 space
 with coordinates
 $x,y,z,p,q$.
 Here
  $x,y$ mean independent variables,
  $z$ means an unknown function $z(x,y)$,
  and
 $p,q$ mean first derivatives
 $z_x, z_y$ respectively.
 There exists the
 $\Lambda^2 \mbox{J}^1 (\mathbb{R}^2)$
 space
 over
 $\mbox{J}^1 (\mathbb{R}^2)$.

 Assuming
 ${p = f_x}$,
 ${q = f_y}$,
 where $f(x,y)$
 is some function over
 $\mathbb{R}^2$,
 we have
$$
\begin{array}{l}
 dx \wedge dq + dy \wedge dp =
 dx \wedge d (f_y) +
 dy \wedge d (f_x) =
 dx \wedge
 \left(
 f_{xy} \: dx +
 f_{yy} \: dy
 \right) +
 \medskip\\
 +
 dy \wedge
 \left(
 f_{xx} \: dx +
 f_{xy} \: dy
 \right) =
 f_{xy} \:
 dx \wedge dx +
  f_{yy} \:
 dx \wedge dy +
  f_{xx} \:
 dy \wedge dx +
 \medskip\\
 +
  f_{xy} \:
 dy \wedge dy =
  f_{yy} \:
 dx \wedge dy -
  f_{xx} \:
 dx \wedge dy =
 \left(
  f_{yy} -
  f_{xx}
 \right)
 dx \wedge dy
\end{array}
$$
 and
$$
\begin{array}{l}
 dp \wedge dq + dx \wedge dy =
 d
  (f_x)
  \wedge d
  (f_y)
  + dx \wedge dy =
 \left(
  f_{xx}
  \: dx +
  f_{xy}
  \: dy
 \right)
 \wedge
 (
 f_{xy}
  \: dx +
 \medskip\\
 +
  f_{yy}
  \: dy
 )
 +
 dx \wedge dy
 =
 f_{xx}
 f_{xy}
 \: dx \wedge dx
 +
  f_{xx}
  f_{yy}
 \: dx \wedge dy
 +
 f_{xy}
 f_{xy}
 \: dy \wedge dx
 +
 \medskip\\
 +
 f_{xy}
 f_{yy}
 \: dy \wedge dy
 +
 dx \wedge dy
 =
 \left(
  f_{xx}
  f_{yy}
  -
    \left(
      f_{xy}
    \right)
    ^2
 \right)
 \: dx \wedge dy
 +
 dx \wedge dy
 =
 \medskip\\
 =
 \left(
 \mbox{hess}\: f + 1
 \right)
 dx \wedge dy \, .
\end{array}
$$
 The next formulae for exterior forms
 \cite{Bourbaki}
 were used:
$$
\begin{array}{l}
  \omega_1 \wedge (\omega_2 + \omega_3) =
  \omega_1 \wedge \omega_2 +
  \omega_1 \wedge \omega_3 \, ,
  \smallskip \\
  a \omega_1 \wedge \omega_2 = \omega_1 \wedge a \omega_2
  = a (\omega_1 \wedge \omega_2) \, ,
  \smallskip \\
  \omega \wedge \omega = 0 \, ,
\end{array}
$$
 where $\omega,\, \omega_1,\, \omega_2,\, \omega_3$
 are exterior forms,
 $a$ is a constant.

 Therefore, to each form
$$
 dx \wedge dq + dy \wedge dp
$$
 from
 $\Lambda^2 \mbox{J}^1 (\mathbb{R}^2)$
 we assign the form
$$
 \left(
  f_{yy} -
  f_{xx}
 \right)
 dx \wedge dy
$$
 from
 $\Lambda^2 (\mathbb{R}^2)$,
 or a linear wave equation,
 and to each form
$$
dp \wedge dq + dx \wedge dy
$$
 from
 $\Lambda^2 \mbox{J}^1 (\mathbb{R}^2)$
 we assign the form
$$
 \left(
 \mbox{hess}\: f + 1
 \right)
 dx \wedge dy
$$
 from
 $\Lambda^2 (\mathbb{R}^2)$,
 or a simple Monge--Ampere equation.

 Consider the Ampere transformation
$$
 \bar x = - p , \quad
 \bar y = y , \quad
 \bar z = z - p \: x , \quad
 \bar p = x , \quad
 \bar q = q .
$$
 It is a contact transformation,
 i.e. it conserves the form
 ${dz - p\: dx - q \: dy}$.
 Namely,
$$
\begin{array}{l}
   d \bar z - \bar p \: d \bar x - \bar q \: d \bar y =
   d(z - p \: x) - x \: d(-p) - q \: dy =
  \smallskip \\
   =
   dz - dp \: x - p \: dx + x \: dp - q \: dy =
   dz - p \: dx - q \: dy \, .
\end{array}
$$
 The Ampere transformation
 takes the Monge--Ampere equation
${\mbox{hess} \: z = -1}$
 to the linear wave equation
 ${z_{xx} - z_{yy} = 0}$,
 because
$$
 d \bar p \wedge d \bar q
 +
 d \bar x \wedge d \bar y
 =
 dx \wedge dq
 +
 d(-p) \wedge dy
 =
 dx \wedge dq
 -
 dp \wedge dy
 =
 dx \wedge dq
 +
 dy \wedge dp
 \, .
$$
 The Ampere transformation
 takes the classical solution
 ${z=xy}$
 of the equation
 ${\mbox{hess} \: z = -1}$,
 which is a 2-dimensional integral variety
 ${(u,v,uv,v,u)}$,
 to the integral variety
 ${(-v,v,0,u,u)}$,
 which is a multivalued solution of the wave equation
 ${z_{xx} - z_{yy} = 0}$.
 The projecting map from
 $\mathbb{R}^5$
 to
 ${\mathbb{R}^2 = (x,y)}$
 takes
 the 2-dimensional integral variety
 ${(-v,v,0,u,u)}$
 to the 1-dimensional line
 ${(-v,v)}$,
 which isn't a 2-dimensional domain.
 So
 the 2-dimensional integral variety
 ${(-v,v,0,u,u)}$
 at any its point
 couldn't be used as a classical solution
 of the wave equation.
 $\square$

\end{document}


УДК 517.957

 \begin{center}
 \textbf{Ю.~Н.~Братков}\\
 \textbf{
 Гиперболическое уравнение Монжа-Ампера: классические решения на
 всей плоскости
 }
 \end{center}


 \begin{abstract}
 \noindent
 На плоскости ${\mathbb{R}^2 = (x,y)}$ рассматривается
 задача Коши для гиперболического уравнения Монжа-Ампера
$$
 \left\{
 \begin{array}{l}
   A + Bz_{xx} + Cz_{xy} + Dz_{yy} + \mbox{hess}\: z = 0 \, ,
   \medskip \\
   z(0,y) = z^0(y),
   \quad
   z_x(0,y) = p^{\: 0}(y), \ \ y \in \mathbb{R}
   \, .
 \end{array}
 \right.
$$
Здесь ${\mbox{hess}\: z = z_{xx}z_{yy} - z_{xy}^2}$, коэффициенты
$A,B,C,D$ зависят от $x,y,z,z_x,z_y$. Уравнение гиперболично, если
${C^2 - 4BD + 4A > 0}$. Формулируются достаточные условия
существования (единственного) $C^3\mbox{-решения}$ на всей
плоскости.
 \end{abstract}


\tableofcontents




\newcounter{lemma_counter}
\newcounter{theorem_counter}
\newcounter{example_counter}






\section{Введение}

Эта известная задача была поставлена автору Е.~В.~Шикиным. Первая
публикация сделана в \cite{Братков-ФПМ} (получено редакцией в
феврале 1998). Автор благодарен Д.~В.~Туницкому, обнаружившему
весьма тонкую (непринципиальную) ошибку. Ниже публикуется полный
текст исправленного доказательства: заменен метод последовательных
приближений. Эта замена на результат не влияет. Пользуясь случаем,
автор расширил и дополнил результат.

Используется сведение к системе пяти уравнений в римановых
инвариантах \cite{Туницкий-деп}. Теория гиперболических систем
весьма популярна для случая, когда собственные значения системы
отделимы (например, разделены константами). В случае уравнения
Монжа-Ампера собственные значения совпадают с решением системы,
т.е. они неизвестны и их требуется найти.

 Суть проблемы сформулировал в принстонских лекциях
 (1953) Ж.~Лере (\cite{Лере}, гл.~IX): "Удается доказать
только \textit{локальную теорему существования}\dots\ Она
показывает, что для гиперболических уравнений существование
решений в целом зависит от получения \textit{априорных оценок} для
их производных".

Построение классических решений на всей плоскости --- особый вид
спорта. Если в заданной области решение имеет особенности, то их
заменой переменных можно вывести за пределы области. К этому
сводятся многие работы. Но если заданная область
--- вся плоскость, то выводить сингулярности некуда. Мы получаем
качественно иную задачу.

Еще Б.~Риман показал, что некая конкретная гиперболическая система
\cite{Ян} регулярных решений на всей плоскости не имеет. Класс
слабо нелинейных систем
$$
 \begin{array}{l}
    {\ds
      \left(
        \pt_x
        +
        \xi_1
        (u_{2})
        \pt_y
      \right)
      u_1
      =
      0,
    }
   \medskip
   \\
    {\ds
      \left(
        \pt_x
        +
        \xi_2
        (u_{1})
        \pt_y
      \right)
      u_2
      =
      0
    }
 \end{array}
$$
(имеется в виду ${\pt \xi_i/\pt u_i = 0},$ ${i=1,2}$) как систем,
имеющих регулярные решения на всей плоскости, ввел в 1955~г.
Н.~Н.~Яненко \cite{Ян}. Рассмотрение систем двух уравнений в
инвариантах с ненулевой правой частью было проведено в 1967~г.
Б.~Л.~Рождественским и А.~Д.~Сидоренко (\cite{РЯ}, гл.~1, $\S$ 10,
п.~3). Подходящим оказался класс (слабо нелинейных)
 систем, гиперболических в узком смысле, т.е. систем с отделимыми
собственными значениями. Таким образом, теорема
Рождественского-Сидоренко переводит проблему априорной оценки
производных в проблему априорной оценки отделимости собственных
значений (${\xi_1 \neq \xi_2}$ на всей плоскости).

Приведем пример. Гиперболическое уравнение Монжа-Ампера с
коэффициентами, зависящими только от ${x,y}$, сводится к системе
$$
 \begin{array}{l}
    {\ds
      \left(
        \pt_x
        +
        u_2
        \pt_y
      \right)
      u_1
      =
      (u_1 - u_2)\, a_1(x,y,u_1),
    }
   \medskip
   \\
    {\ds
      \left(
        \pt_x
        +
        u_1
        \pt_y
      \right)
      u_2
      =
      (u_1 - u_2)\, a_2(x,y,u_2).
    }
 \end{array}
$$
 Вычтем из первого уравнения второе, разделим на
 $u_1 - u_2$ и,
интегрируя вдоль характеристики, получим равенство
$$
 (u_1 - u_2)(x,y)
 =
 (u_1^0(y) - u_2^0(y))
 \exp
 \left\{
    {\int_0^x}
      (a_1 - a_2 - u_{2y})
    \: d\tau
 \right\}.
$$
 Из него следует, что если ${u_1^0(y) \neq u_2^0(y)}$
 для
 ${\forall\ y\in\mathbb{R},}$ а
 $u_1,\ u_2,\ u_{2y}$ не уходят в бесконечность
 ни в какой конечной точке, то ${u_1 - u_2 \neq 0}$ в любой
 конечной точке.
 Попросту говоря, здесь мы имеем эквивалентность
 априорных оценок для производных и для $u_1, \ u_2, \ {u_1 - u_2}$,
 где условие ${u_1 - u_2 \neq 0}$ --- это условие
 ${\xi_1 \neq \xi_2}$.

 Подход автора: автор волевым образом задает априорную оценку для
 $u_1, \ u_2, \ {u_1 - u_2}$.
 Мышление в этом направлении было блокировано.

 Упомянем работы Чан-Шен
 Хона (Jia-Xing Hong)
 \cite{Хон-93}, \cite{Хон-95},
 о которых автор узнал в 2002~г. в Пекине.
 Статью \cite{Хон-95}
 о глобальных классических решениях уравнения
$$
 \mbox{hess} \, z = -k^2(x,y)
$$
 автор не смог понять, а
 статью \cite{Хон-93},
 которую можно рассматривать, в частности, как решение уравнения
$$
 \mbox{hess} \, z = -k^2(x,y)(1 + z_x^2 + z_y^2)^2,
$$
 не смог достать.


\section{Системы в римановых инвариантах}
\subsection{Одно модельное уравнение}

 Пусть имеется плоскость ${\mathbb{R}^2 = (x,y)}$,
и пусть $u(x,y)$ --- неизвестная функция. Рассмотрим на примере
уравнения
$$
    \left( \pt_x
           +\xi(x,y) \pt_y
    \right)u(x,y)=f(x,y)
    \eqno(2.1)
$$
основные понятия характеристики и интегрирования вдоль
характеристики (\cite{РЯ}, гл.~1). Пусть $\xi ,f \in
C^1(\mathbb{R}^2).$

{\bf Определение.} Кривая на плоскости $(x,y),$ задаваемая
уравнениями
$$
  x=\tau, \qquad y=g(\tau ,x,y),
$$
где функция $g(\tau ,x,y)$ --- решение задачи Коши
$$
 \left\{
   \begin{array}{l}
     \pt_\tau g=\xi (\tau ,g(\tau ,x,y)),
     \smallskip
     \\
     g(x,x,y)=y,
   \end{array}
 \right.
  \eqno(2.2)
$$
 называется {\it характеристикой} уравнения
 (2.1).

 {\bf Лемма~2.1.}
$$
    \left( \pt_x
           +\xi \pt_y
    \right)g=0.
    \eqno(2.3)
$$

 {\bf Доказательство.} Предположим, что функции
$g_x,\: g_y$ существуют, и найдем их явный вид, следуя \cite{Поз}.
Продифференцируем соотношение
 (2.2)
 по~$x:$
$$
  \pt_\tau \pt_x g
  =
  (\pt \xi/\pt g)
  (\tau ,g)
  \:
  \pt_x g.
 \eqno(2.4)
$$
Получаем линейное обыкновенное дифференциальное уравнение
относительно $\pt_x g.$ Начальное условие для него получаем
дифференцированием начального условия $(2.2):$
$$
  \pt_x g(x,x,y)=0,
$$
т.е.
$$
  \left.\pt_\tau g(\tau ,x,y)\right| _{\tau =x}+
  \left.\pt_x g(\tau ,x,y)\right| _{\tau =x}=0,
$$
откуда, с учетом (2.2),
$$
  \left.\pt_x g(\tau ,x,y)\right| _{\tau =x}=-\xi (x,y).
  \eqno(2.4^0)
$$
Решая задачу Коши (2.4), $(2.4^0)$, получаем:
$$
  \pt_x g(\tau ,x,y)=-\xi (x,y)\exp
  \left\{
    \int _x^{\tau}
      \frac{\ds \pt \xi}{\ds \pt g}
      (t,g(t,x,y))
    \: dt
  \right\}.
  \eqno(2.5)
$$
Аналогично, решая задачу Коши для функции $\pt_y g$
$$
\left\{
  \begin{array}{l}
    \pt_\tau \pt_y g
    =
      \frac{\ds \pt \xi}{\ds \pt g}(\tau ,g)
      \pt_y
      g,
      \\
    \pt_y g(\tau ,x,y)| _{\tau =x}=1,
  \end{array}
\right.
$$
имеем
$$
  \pt_y g(\tau ,x,y)=\exp
  \left\{
    \int _x^{\tau}
      \frac{\ds \pt \xi}{\ds \pt g}
      (t,g(t,x,y))
    \: dt
  \right\}.
  \eqno(2.6)
$$
Из формул (2.5), (2.6) следует (2.3). $\Box $

Поставим для уравнения (2.1) задачу Коши, задав начальное условие
на оси~$Oy$:
$$
  u(0,y)=u^0(y),
  \eqno(2.1^0)
$$
где $u^0\in C^1({\mathbb{R}}^1).$

Задача (2.1), $(2.1^0)$ корректна \cite{РЯ}. Будем решать ее в
полуплоскости ${x\ge 0}$ (построение решения в полуплоскости
${x\le 0}$ проводится аналогично).

 {\bf Лемма~2.2.}
 Решение задачи (2.1), $(2.1^0)$ задается формулой
$$
  u(x,y)=u^0(g(0,x,y))+\int _0^x f(\tau ,g(\tau ,x,y))d\tau .
  \eqno(2.7)
$$

{\bf Доказательство.} Условие $(2.1^0)$ выполнено:
$$
  u(0,y)=u^o(g(0,0,y))=u^0(y)
$$
в силу начального условия (2.2). При дифференцировании соотношения
(2.7) воспользуемся формулой для дифференцирования интеграла,
зависящего от параметра:
$$
  \pt_x
  \int _0^x
    \varphi(\tau ,x)
  d\tau =
  \int _0^x
    \pt_x\varphi(\tau ,x)
  d\tau
  +\varphi (x,x).
$$
Имеем:
$$
    \left( \pt_x
           +\xi \pt_y
    \right)
    u
    =
    \left.
    u_t^0(t)
    \right|_{t=g(0,x,y)}
    \left( \pt_x
           +\xi \pt_y
    \right)
    g(0,x,y)+
$$
$$
  +\int _0^x
    \frac{\ds \pt f}{\ds \pt g}(\tau ,g(\tau ,x,y))
    \left( \pt_x
           +\xi \pt_y
    \right)
    g(\tau ,x,y)
  d\tau
  +f(x,g(x,x,y)).
$$
 Из леммы~2.1
 и начального условия (2.2) следует, что
$$
  \left(
    \pt_x
    +\xi \pt_y
  \right)
  u(x,y)=f(x,y). \quad \Box
$$

   Процедура получения решения задачи (2.1), $(2.1^0)$ при помощи формулы (2.7)
называется {\it интегрированием уравнения} (2.1) {\it вдоль
характеристики.}


\subsection{Система в римановых инвариантах}

  Развитая выше теория работает и в более общем случае. Пусть имеются
вектор-функции
$$
\begin{array}{lll}
 u=(u_1,\dots ,u_m),\quad & u_i=u_i(x,y)\in C^1({\mathbb{R}}^2),&
 \smallskip
 \\
 \xi =(\xi _1,\dots ,\xi _m),\quad &\xi _i=\xi_i(x,y,u)\in C^1({\mathbb{R}}^2),&
 \smallskip
 \\
 f=(f_1,\dots ,f_m),\quad & f_i=f_i(x,y,u)\in C^1({\mathbb{R}}^2),&
 i=1,\dots ,m.
\end{array}
$$
{\bf Определение.} Система вида
$$
    \left( \pt_x
           +\xi_i(x,y,u)
           \pt_y
    \right)
    u_i(x,y)=f_i(x,y,u), \quad i=1,\dots ,m,
    \eqno(2.8)
$$
называется {\it системой в римановых инвариантах} \cite{РЯ}.

   Поставим для системы (2.8) задачу Коши, задав на оси $Oy$ начальные условия
$$
  u(0,y)=u^o(y),
  \eqno(2.8^0)
$$
где $u^0=(u_1^0,\dots ,u_m^0),\quad
     u_i^0=u_i^0(y)\in C^1(\mathbb{R}),\quad i=1,\dots ,m.$\\
{\bf Определение.} Кривая на плоскости $(x,y),$ задаваемая
уравнениями
$$
  x=\tau, \qquad y=g_i(\tau ,x,y),
$$
где функция $g_i(\tau ,x,y)$ --- решение задачи Коши
$$
\left\{
\begin{array}{lr}
    \pt_\tau g_i
    =
    \xi _i
     (\tau,
      g_i(\tau ,x,y),
      u(\tau ,x,y)
     ),
  \smallskip
  \\
    g_i(x,x,y)=y,
\end{array}
\right.
 \eqno (2.9)
$$
называется {\it i-й характеристикой} системы (2.8).\\
 \addtocounter{lemma_counter}{1}
 {\bf Лемма~2.3.}
 $i\mbox{-я}$ компонента решения задачи
 (2.8), $(2.8^0)$
 может быть представлена как результат интегрирования
$i\mbox{-го}$ уравнения системы вдоль $i\mbox{-й}$ характеристики,
${i=1,\dots ,m:}$
$$
  u_i(x,y)=u_i^0(g_i(0,x,y))+
  {\ds \int _0^x}
    f_i(\tau ,g_i(\tau ,x,y),
      u(\tau ,g_i(\tau ,x,y)))
   \: d\tau .
 \eqno(2.10)
$$

{\bf Доказательство.} Предположим, что решение $u(x,y)$ задачи
(2.8), $(2.8^0)$ известно. Подставляя его в (2.8), убеждаемся в
том, что каждое уравнение системы имеет вид (2.1), так как можно
считать, что
$$
  \xi _i(x,y,u(x,y))=\tilde \xi (x,y),\quad
  f_i(x,y,u(x,y))=\tilde f(x,y),\quad i=1,\dots ,m.
$$
Поэтому для уравнений
$$
    (
      \pt_x
           +\tilde \xi _i (x,y) \pt_y
    )
    \: u_i(x,y)=\tilde f_i(x,y),\quad i=1,\dots ,m,
$$
справедливы леммы 2.1 и 2.2, при этом равенства, аналогичные (2.5)
и (2.6), записываются в виде
$$
  \!\!\!
{\ds
\begin{array}{l}
  \!\!\!
  \pt_x
  g_i(\tau ,x,y)
  \!
  =
  \!\!
  -\xi _i (x,y,u)
  \!
  \exp
  \!
  \left\{ \!
    {\ds\int _x^{\tau}}
    \!\!
    \left[
        \xi _{iy} \! + \!\!\!
        {\ds\sum}_{j=1}^m
        \frac{\ds \pt \xi_i}{\ds \pt u_j}\frac{\ds \pt u_j}{\ds \pt y}
    \right]
    \!\!
    (t,g_i(t,x,y))dt \!
 \right\}
 \!
 ,
 \\
 \\
  \!\!\!\pt_y g_i(\tau ,x,y)=
  \exp
  \left\{
    {\ds\int _x^{\tau}}
    \left[
        \xi _{iy}+
        {\ds\sum}_{j=1}^m
            \frac{\ds \pt \xi_i}{\ds \pt u_j}
            \frac{\ds \pt u_j}{\ds \pt y}
    \right]
    (t,g_i(t,x,y))dt
 \right\}
 \!
 ,
\end{array}
}
  \!\!\!\!\!
 \eqno(2.11)
$$
а равенство (2.7) имеет вид (2.10). $\Box$

Формула (2.10) используется при изучении свойств решения. Она
рекурсивна и поэтому не дает решение задачи (2.8), $(2.8^0)$ в
явном виде. Построение решения можно проводить методом
последовательных приближений.



\section{Системы в римановых инвариантах
         и уравнения Монжа--Ампера гиперболического типа}

 Здесь мы следуем статье Туницкого \cite{Туницкий-деп},
 вошедшей в \cite{Туницкий-канд-дис}.
 Безусловно, о том, что гиперболическое уравнение Монжа--Ампера
 можно свести к системе пяти уравнений первого порядка,
 известно давно (\cite{Ку}, прилож.~1 к гл.~V, \S~2),
 и что эта система сводима к системе в
 римановых инвариантах, также было известно
 (в \cite{Ку} система не в инвариантах).
 У Туницкого система в инвариантах выписана в явном виде,
 и ей можно пользоваться.


\subsection{Задача Коши для уравнения Монжа--Ампера}

На плоскости ${\mathbb{R}^2 = (x,y)}$ \ будем рассматривать
уравнение Монжа--Ампера относительно неизвестной функции
$z=z(x,y)$
   $$
   A+Bz_{xx}+Cz_{xy}+Dz_{yy}+E\,\mbox{\rm hess}\, z=0,
   $$
где $\mbox{\rm hess}\, z=z_{xx}z_{yy}-z^2_{xy},\quad A,B,C,D,E$
--- функции, зависящие от  $x,y,z$,\ $z_x,z_y$; \;\; ${E\ne 0}$. Так
как делением на $E$ всегда можно добиться того, что ${E\equiv 1}$,
то без ограничения общности будем рассматривать уравнение
$$
   A+Bz_{xx}+Cz_{xy}+Dz_{yy}+\mbox{\rm hess}\, z = 0. \eqno (3.1)
$$
 В дальнейшем будем считать, что
 ${A,B,C,D \in C^2(\mathbb{R}^5)}$.

 Предположим, что ${z(x,y)}$ --- дважды непрерывно
 дифференцируемое решение уравнения (3.1). Будем говорить, что
 уравнение (3.1) гиперболично
 на решении ${z(x,y)}$
 (\cite{Ку}, прил.~1 к гл.~5, ${{\S}~2}$),
 если
 $$
 \Delta^2(x,y,z(x,y),z_x(x,y),z_y(x,y))
 =
 C^2 - 4BD + 4AE > 0.
 \eqno(3.2)
 $$
 У нас здесь всегда будет ${E\equiv 1}$.
 Поскольку ${\Delta}^2 > 0$, положим также
 ${\Delta > 0}$.
 Для уравнения (3.1) на оси $Oy$
 поставим задачу Коши
$$
 z(0,y) = z_0(y), \qquad z_x(0,y)=p_0(y);
 \eqno(3.3)
$$
 здесь ${z_0 \in C^3(\mathbb{R}^1)}$, а
 ${p_0 \in C^2(\mathbb{R}^1)}$.
 Предположим, что для начальных значений $z_0$ и $p_0$
 выполнены два условия. Во-первых, ось $Oy$ свободна, то есть
$$
 z_0^{''}(y)
 +
 B
 (0,y,z_0(y),p_0(y),z_0^{'}(y))
 \neq
 0.
 \eqno(3.4)
$$
 Во-вторых, на оси ординат выполнено условие гиперболичности (3.2)
$$
 \Delta^2
 (0,y,z_0(y),p_0(y),z_0^{'}(y))
 >
 0.
 \eqno(3.5)
$$
 Полная формулировка задачи Коши для уравнения Монжа--Ампера будет
 следующей. Требуется найти трижды непрерывно дифференцируемую функцию
 ${z(x,y)}$, которая принимает начальные значения (3.3),
 удовлетворяет уравнению (3.1), и на которой это уравнение
 гиперболично.

 Поясним, откуда берется ${C^3\mbox{-дифференцируемость}}$,
 тогда как классическое решение уравнения (3.1) подразумевается
 ${C^2\mbox{-гладким}}$. Для
 вывода системы в римановых инвариантах будут использованы условия
 интегрируемости
 ${p_{xy} = p_{yx}}$, ${q_{xy} = q_{yx}}$,
 где $p,\ q$ имеют смысл $z_x,\ z_y$.

\subsection{Система в полных дифференциалах}

 Пусть ${z(x,y)}$ --- $C^3\mbox{-решение}$ уравнения (3.1) в некоторой
 области $T$. Потребуем, чтобы уравнение (3.1) было гиперболично
 на решении $z$ и
$$
 z_{yy}(x,y)
 +
 B
 (x,y,z(x,y),z_x(x,y),z_y(x,y))
 \neq
 0
 \eqno(3.6)
$$
 для всех точек ${(x,y)}$ множества $T$. Неравенство (3.6)
 означает, что вертикальные прямые ${x=\mbox{const}}$ свободны.

 Положим
$$
 u_1
 =
 \frac
 {\ds C + \Delta - 2z_{xy}}
 {2(z_{yy} + B)}
 ,
 \qquad
 u_2
 =
 \frac
 {\ds C - \Delta - 2z_{xy}}
 {2(z_{yy} + B)}
 .
 \eqno(3.7)
$$
 Значения функций $u_1$ и $u_2$ совпадают с тангенсами углов
 наклона характеристик уравнения (3.1)
 (\cite{Ку}, прил.~1 к гл.~V, \S~2). Из условия гиперболичности
 (3.2) вытекает неравенство ${u_1 \neq u_2}$, которое позволяет
 однозначно разрешить соотношения (3.7) относительно $z_{xy}$ и
 $z_{yy}$:
$$
 z_{yy}
 =
 \frac
 {\ds \Delta}
 {u_1 - u_2}
 -
 B
 ,
 \qquad
 z_{xy}
 =
 \frac{\ds \Delta}{\ds 2}
 \frac
 {u_1 + u_2}
 {u_2 - u_1}
 +
 \frac{\ds C}{\ds 2}
 .
 \eqno(3.8)
$$
 Подставим выражения (3.8) в уравнение (3.1). Получим линейное
 уравнение относительно $z_{xx}$. Решая его, находим
$$
 z_{xx}
 =
 \Delta
 \frac
 {u_1 u_2}
 {u_1 - u_2}
 -
 D
 .
 \eqno(3.9)
$$
 Следуя Монжу, положим
$$
 z_x = p, \qquad z_y = q.
 \eqno(3.10)
$$
 Из равенств (3.8)--(3.9) ясно, что три функции $z$, $p$ и $q$
 переменных $x$ и $y$ удовлетворяют в области $T$ системе
 уравнений в полных дифференциалах
$$
\begin{array}{l}
 p_x
 =
 \Delta
 (x,y,z,p,q)
 \frac
 {\ds u_1 u_2}
 {\ds u_1 - u_2}
 (x,y)
 -
 D
 (x,y,z,p,q)
 ,
 \medskip
 \\
 p_y
 =
 \frac{\ds \Delta}{\ds 2}
 (x,y,z,p,q)
 \frac
 {\ds u_1 + u_2}
 {\ds u_2 - u_1}
 (x,y)
 +
 \frac{\ds C}{\ds 2}
 (x,y,z,p,q)
 ,
 \medskip
 \\
 q_x
 =
 \frac{\ds \Delta}{\ds 2}
 (x,y,z,p,q)
 \frac
 {\ds u_1 + u_2}
 {\ds u_2 - u_1}
 (x,y)
 +
 \frac{\ds C}{\ds 2}
 (x,y,z,p,q)
 ,
 \medskip
 \\
 q_y
 =
 \frac
 {\ds \Delta (x,y,z,p,q)}
 {\ds (u_1 - u_2)(x,y)}
 -
 B
 (x,y,z,p,q)
 .
\end{array}
 \eqno(3.11)
$$
 Правые части уравнений (3.11) --- непрерывно дифференцируемые
 функции от $x,y,z,p$ и $q$.

 Таким образом, нами установлено следующее. Если ${z(x,y)}$ ---
 трижды непрерывно дифференцируемое решение уравнения (3.1) в
 области $T$, для которого выполнены неравенства (3.2) и (3.6), то
 определены три функции $z,p$ и $q$, являющиеся решением системы
 уравнений (3.10)--(3.11).

 Обратно, пусть в области $T$ существует пара ${(u_1,\ u_2)}$
 непрерывно дифференцируемых функций таких, что ${u_1 \neq u_2}$,
 и в этой области существует $C^1\mbox{-решение}$ системы уравнений
 (3.10)--(3.11) такое, что
$$
 \Delta
 (x,y,z(x,y),z_x(x,y),z_y(x,y))
 >
 0
 .
 \eqno(3.12)
$$
 Понятно, что в этом случае ${z_{xx} = p_x}$, ${z_{xy} = p_y =
 q_x}$, ${z_{yy} = q_y}$. Следовательно, ${z \in C^3(T)}$.
 Подставляя значения вторых производных функции $z$ в левую часть
 соотношения (3.1), убеждаемся, что ${z(x,y)}$ --- решение
 уравнения (3.1). В силу оценки (3.12) это уравнение гиперболично
 на $z$, и справедливо неравенство (3.6).

 Резюмируем полученные результаты в виде леммы.

 {\bf Лемма~3.1.}
 В области $T$ тогда и только тогда определена функция ${z(x,y)}$ класса
 ${C^3(T)}$, удовлетворяющая соотношениям (3.1), (3.2) и (3.6),
 когда в этой области существует $C^1\mbox{-решение}$ ${(z,p,q)}$ системы
 уравнений (3.10)--(3.11), удовлетворяющее условию (3.12). $\Box$

\subsection{Система в римановых инвариантах}

 Система дифференциальных уравнений в полных дифференциалах
 является переопределенной и, вообще говоря, решения не имеет. Для
 существования решения необходимо и достаточно, чтобы выполнялись
 условия интегрируемости. Условия
 интегрируемости уравнений (3.10)--(3.11) сводятся к
 дифференциальным соотношениям между функциями $u_1, u_2, z, p$ и
 $q$.

 Действительно, так как правые части системы (3.10)--(3.11)
 непрерывно дифференцируемы, то функции $z,\ p$ и $q$ дважды
 непрерывно дифференцируемы. Следовательно, на множестве $T$ имеют
 место равенства
$$
 z_{xy} = z_{yx},
 \qquad
 p_{xy} = p_{yx},
 \qquad
 q_{xy} = q_{yx}.
 \eqno(3.13)
$$
 Первое из равенств (3.13) выполняется тождественно в силу системы
 (3.10)--(3.11). Второе приводится к виду
$$
\begin{array}{ll}
 &
 \Delta_y
 \frac
 {\ds u_1 u_2}
 {\ds u_1 - u_2}
 +
 \Delta
 \frac
 {\ds u_1^2 u_{2y} - u_2^2 u_{1y}}
 {\ds (u_1 - u_2)^2}
 -
 D_y
 +
 \left(
    \frac
    {\ds u_1 u_2}
    {\ds u_1 - u_2}
    \Delta_z
    -
    D_z
 \right)
 q
 +
 \medskip
 \\
 &
 +
 \left(
    \frac
    {\ds u_1 u_2}
    {\ds u_1 - u_2}
    \Delta_p
    -
    D_p
 \right)
 \left(
    \frac{\ds \Delta}{\ds 2}
    \frac
    {\ds u_1 + u_2}
    {\ds u_2 - u_1}
    +
    \frac{\ds C}{\ds 2}
 \right)
 +
 \medskip
 \\
 &
 +
 \left(
    \frac
    {\ds u_1 u_2}
    {\ds u_1 - u_2}
    \Delta_q
    -
    D_q
 \right)
 \left(
    \frac
    {\ds \Delta}
    {\ds u_1 - u_2}
    -
    B
 \right)
 =
 \medskip
 \\
 =
 &
 \frac{\ds \Delta_x}{\ds 2}
 \frac
 {\ds u_1 + u_2}
 {\ds u_2 - u_1}
 +
 \Delta
 \frac
 {\ds u_2 u_{1x} - u_1 u_{2x}}
 {\ds (u_2 - u_1)^2}
 +
 \frac{\ds C_x}{\ds 2}
 +
 \left(
    \frac
    {\ds u_1 + u_2}
    {\ds u_2 - u_1}
    \frac{\ds \Delta_z}{\ds 2}
    +
    \frac{\ds C_z}{\ds 2}
 \right)
 p
 +
 \medskip
 \\
 &
 +
 \left(
    \frac
    {\ds u_1 + u_2}
    {\ds u_2 - u_1}
    \frac{\ds \Delta_p}{\ds 2}
    +
    \frac{\ds C_p}{\ds 2}
 \right)
 \left(
    \Delta
    \frac
    {\ds u_1 u_2}
    {\ds u_1 - u_2}
    -
    D
 \right)
 +
 \medskip
 \\
 &
 +
 \left(
    \frac
    {\ds u_1 + u_2}
    {\ds u_2 - u_1}
    \frac{\ds \Delta_q}{\ds 2}
    +
    \frac{\ds C_q}{\ds 2}
 \right)
 \left(
    \frac{\ds \Delta}{\ds 2}
    \frac
    {\ds u_1 + u_2}
    {\ds u_2 - u_1}
    +
    \frac{\ds C}{\ds 2}
 \right)
 .
\end{array}
 \eqno(3.14)
$$
 Аналогичным образом третье из равенств (3.13) сводится к
 соотношению
$$
\begin{array}{ll}
 &
 \frac{\ds \Delta_y}{\ds 2}
 \frac
 {\ds u_1 + u_2}
 {\ds u_2 - u_1}
 +
 \Delta
 \frac
 {\ds u_2 u_{1y} - u_1 u_{2y}}
 {\ds (u_2 - u_1)^2}
 +
 \frac{\ds C_y}{\ds 2}
 +
 \left(
    \frac
    {\ds u_1 + u_2}
    {\ds u_2 - u_1}
    \frac{\ds \Delta_z}{\ds 2}
    +
    \frac{\ds C_z}{\ds 2}
 \right)
 q
 +
 \medskip
 \\
 &
 +
 \left(
    \frac
    {\ds u_1 + u_2}
    {\ds u_2 - u_1}
    \frac{\ds \Delta_p}{\ds 2}
    +
    \frac{\ds C_p}{\ds 2}
 \right)
 \left(
    \frac{\ds \Delta}{\ds 2}
    \frac
    {\ds u_1 + u_2}
    {\ds u_2 - u_1}
    +
    \frac{\ds C}{\ds 2}
 \right)
 +
 \medskip
 \\
 &
 +
 \left(
    \frac
    {\ds u_1 + u_2}
    {\ds u_2 - u_1}
    \frac{\ds \Delta_q}{\ds 2}
    +
    \frac{\ds C_q}{\ds 2}
 \right)
 \left(
    \frac
    {\ds \Delta}
    {\ds u_1 - u_2}
    -
    B
 \right)
 =
 \medskip
 \\
 =
 &
 \frac
 {\ds \Delta_x}
 {\ds u_1 - u_2}
 +
 \Delta
 \frac
 {\ds u_{2x} - u_{1x}}
 {\ds (u_1 - u_2)^2}
 -
 B_x
 +
 \left(
    \frac
    {\ds \Delta_z}
    {\ds u_1 - u_2}
    -
    B_z
 \right)
 p
 +
 \medskip
 \\
 &
 +
 \left(
    \frac
    {\ds \Delta_p}
    {\ds u_1 - u_2}
    -
    B_p
 \right)
 \left(
    \Delta
    \frac
    {\ds u_1 u_2}
    {\ds u_1 - u_2}
    -
    D
 \right)
 +
 \medskip
 \\
 &
 +
 \left(
    \frac
    {\ds \Delta_q}
    {\ds u_1 - u_2}
    -
    B_q
 \right)
 \left(
    \frac{\ds \Delta}{\ds 2}
    \frac
    {\ds u_1 + u_2}
    {\ds u_2 - u_1}
    +
    \frac{\ds C}{\ds 2}
 \right)
 .
\end{array}
 \eqno(3.15)
$$

 Соответствия (3.14)--(3.15) представляют собой линейную систему
 двух алгебраических уравнений с двумя неизвестными
 ${u_{1x} + u_2 u_{1y}}$
 и
 ${u_{2x} + u_1 u_{2y}}$.
 Ее определитель равен
 ${\Delta^2/(u_1 - u_2)^3}$,
 а так как ${u_1 \neq u_2}$,
 то эта система однозначно разрешима:
$$
\left\{
\begin{array}{l}
    u_{1x} + u_2 u_{1y}
    =
    E_o + E_1 u_1 + E_2 u_2 + E_3 u_1^2 +
    E_4 u_1 u_2 + E_5 u_1^2 u_2
    ,
 \medskip
 \\
    u_{2x} + u_1 u_{2y}
    =
    I_o + I_1 u_1 + I_2 u_2 + I_3 u_2^2 +
    I_4 u_1 u_2 + I_5 u_1 u_2^2
    .
\end{array}
\right.
 \eqno(3.16)
$$
 Коэффициенты $E_j,\ I_j,\ {0 \le j \le 5}$,
 определяются через функции $B,C,D,\Delta$
 и их первые производные следующим образом:
$$
\begin{array}{l}
 E_0 = I_0 = D_q,
 \qquad
 E_5 = I_5 = -B_p,
 \medskip
 \\
 E_1
 =
 \phantom{-}
 \alpha_1 + \alpha_2
 +
 \frac{\ds 1}{\ds 4\Delta}
 (
    \phantom{-}
    CC_q - 3\Delta C_q + C\Delta_q + \Delta\Delta_q - 2\Delta D_p
 )
 ,
 \medskip
 \\
 E_2
 =
 -\alpha_1 - \alpha_2
 +
 \frac{\ds 1}{\ds 4\Delta}
 (
   -CC_q - \phantom{3}\Delta C_q - C\Delta_q - \Delta\Delta_q - 2\Delta D_p
 )
 ,
 \medskip
 \\
 E_3
 =
 -\beta_1 + \beta_2
 +
 \frac{\ds 1}{\ds 4\Delta}
 (
   -CC_p + \phantom{3}\Delta C_p + C\Delta_p - \Delta\Delta_p + 2\Delta B_q
 )
 ,
 \medskip
 \\
 E_4
 =
 \phantom{-}
 \beta_1 - \beta_2
 +
 \frac{\ds 1}{\ds 4\Delta}
 (
   \phantom{-}
   CC_p + 3\Delta C_p - C\Delta_p + \Delta\Delta_p + 2\Delta B_q
 )
 ,
 \medskip
 \\
 I_1
 =
 \phantom{-}
 \alpha_1 - \alpha_2
 +
 \frac{\ds 1}{\ds 4\Delta}
 (
    \phantom{-}
    CC_q - \phantom{3}\Delta C_q - C\Delta_q + \Delta\Delta_q - 2\Delta D_p
 )
 ,
 \medskip
 \\
 I_2
 =
 -\alpha_1 + \alpha_2
 +
 \frac{\ds 1}{\ds 4\Delta}
 (
   -CC_q - 3\Delta C_q + C\Delta_q - \Delta\Delta_q - 2\Delta D_p
 )
 ,
 \medskip
 \\
 I_3
 =
 \phantom{-}
 \beta_1 + \beta_2
 +
 \frac{\ds 1}{\ds 4\Delta}
 (
   \phantom{-}
   CC_p + \phantom{3}\Delta C_p + C\Delta_p + \Delta\Delta_p + 2\Delta B_q
 )
 ,
 \medskip
 \\
 I_4
 =
 -\beta_1 - \beta_2
 +
 \frac{\ds 1}{\ds 4\Delta}
 (
   -CC_p + 3\Delta C_p - C\Delta_p - \Delta\Delta_p + 2\Delta B_q
 )
 ,
\end{array}
 \eqno(3.17)
$$
 где
$$
\begin{array}{l}
 \alpha_1
 =
 \frac{\ds 1}{\ds 2\Delta}
 (
    2D_y + C_x + 2D_z q + C_z p + C D_p - D C_p - 2B D_q
 )
 ,
 \medskip
 \\
 \alpha_2
 =
 \frac{\ds 1}{\ds 2\Delta}
 (
    \Delta_x + \Delta_z p - D \Delta_p
 )
 ,
 \medskip
 \\
 \beta_1
 =
 \frac{\ds 1}{\ds 2\Delta}
 (
    2B_x + C_y + 2B_z p + C_z q + C B_q - B C_q - 2D B_p
 )
 ,
 \medskip
 \\
 \beta_2
 =
 \frac{\ds 1}{\ds 2\Delta}
 (
    \Delta_y + \Delta_z q - B \Delta_q
 )
 .
\end{array}
$$

 Пусть функция ${z(x,y)}$ --- $C^3\mbox{-решение}$ уравнения (3.1) в
 области $T$, для которого имеют место неравенства (3.2) и (3.6).
 Согласно лемме~3.1 функции $z,\ p=z_x$ и $q=z_y$ удовлетворяют
 системе уравнений (3.10)--(3.11), где $u_1$ и $u_2$ определены
 равенствами (3.7).

 Умножим второе из равенств (3.10) на $u_1$ и сложим с первым. В
 результате получим
$$
 z_x + u_1 z_y = p + u_1 q.
 \eqno(3.18)
$$
 Аналогичным образом поступим с равенствами (3.11). Первое из них
 прибавим ко второму, умноженному на $u_2$, а третье сложим с
 четвертым, умноженным на $u_1$. Получаем
$$
 p_x + u_2 p_y = \frac{\ds C + \Delta}{\ds 2} u_2 - D
 \eqno(3.19)
$$
 и
$$
 q_x + u_1 q_y = \frac{\ds C + \Delta}{\ds 2} - B u_1
 \eqno(3.20)
$$
 соответственно.

 Таким образом, доказан следующий факт.
 \\
 \textbf{Теорема~3.1.}
 Пусть ${z(x,y)}$ --- $C^3\mbox{-решение}$ уравнения (3.1) в области $T$
 такое, что справедливы условия (3.2) и (3.6). Тогда набор функций
 ${(u_1,\ u_2,\ z,\ p,\ q)}$, где $u_1,\ u_2$ найдены из выражений
 (3.7), ${p=z_x}$ и ${q=z_y}$, является $C^1\mbox{-решением}$ системы
 пяти уравнений (3.16)--(3.20) в этой области. $\Box$

\subsection{Сведение уравнения Монжа--Ампера к системе в инвариантах}

 Допустим, что область $T$ имеет непустое пересечение с осью $Oy$.
 Следующее утверждение является в определенном смысле обратным по
 отношению к теореме~3.1.

 \addtocounter{theorem_counter}{1}
 \textbf{Теорема~3.2.}
 Пусть
 $(u_1,\ u_2,\ z,\ p,\ q)$
 --- $C^1\mbox{-решение}$ системы (3.16)--(3.20) в области $T$,
 принимающее начальные значения
$$
\begin{array}{l}
 z(0,y) = z_0(y),
 \medskip
 \\
 p(0,y) = p_0(y),
 \medskip
 \\
 q(0,y) = z_0^{'}(y),
 \medskip
 \\
 u_1(0,y)
 =
 {\ds
    \frac
    {(C + \Delta) ((0,y,z_0(y),p_0(y),z_0^{'}(y))) - 2p_0^{'}(y)}
    {2(z_0^{''} + B(0,y,z_0(y),p_0(y),z_0^{'}(y)))}
 }
 ,
 \medskip
 \\
 u_2(0,y)
 =
 {\ds
    \frac
    {(C - \Delta) ((0,y,z_0(y),p_0(y),z_0^{'}(y))) - 2p_0^{'}(y)}
    {2(z_0^{''} + B(0,y,z_0(y),p_0(y),z_0^{'}(y)))}
 }
 ,
\end{array}
 \eqno(3.21)
$$
 и пусть $T$ является областью определенности этого решения.
 Тогда, если выполнена оценка (3.12), то
 $z$ является $C^3\mbox{-решением}$ задачи (3.1)--(3.3) на множестве $T$,
 причем
 ${z_x = p},\ {z_y = q},$ и имеет место неравенство (3.6).
 \\
 \textbf{Доказательство.}
 В первую очередь покажем, что ${u_1 \neq u_2}$ в области $T$.
 Согласно начальным значениям (3.21) и оценке (3.5) имеем
$$
 u_1(0,y) - u_2(0,y)
 =
 {\ds
    \frac
    {\Delta
        (0,y,z_0(y),p_0(y),z_0^{'}(y))
    }
    {z_0^{''} (y) + B
        (0,y,z_0(y),p_0(y),z_0^{'}(y))
    }
 }
 \neq 0
 .
 \eqno(3.22)
$$
 Рассмотрим уравнения (3.16). Вычтем из первого уравнения второе.
 Опираясь на представление (3.17), получим равенство
$$
\begin{array}{l}
  (u_1 - u_2)_x + u_2 (u_1 - u_2)_y
  =
 \medskip
 \\
 \phantom{12345}
  =
  (u_1 - u_2)
  \left(
    2\alpha_2 +
    {\ds
      \frac
      {C \Delta_q - \Delta C_q}
      {2\Delta}
    }
   \right.
    +
 \medskip
 \\
 \phantom{12345}
    +
    (u_1 + u_2)
    \left(
      \beta_2
      +
      \frac{\ds B_q}{\ds 2}
      +
      {\ds
        \frac
        {C \Delta_p + \Delta C_p}
        {4 \Delta}
      }
    \right)
      +
 \medskip
 \\
 \phantom{12345}
    +
  \left.
      (u_2 - u_1)
      \left(
        \beta_1
        +
        {\ds
          \frac
          {C C_p + \Delta \Delta_p}
          {4 \Delta}
        }
      \right)
      -
      B_p u_1 u_2 + u_{2y}
  \right)
  .
\end{array}
 \eqno(3.23)
$$
 Так как $T$ --- область определенности решения системы
 (3.16)--(3.20), то на отрезке ${[0,x]}$ существует решение задачи
 Коши
$$
\left\{
\begin{array}{l}
 {\ds
   \frac{d g}{d\tau}
 }
   =
   u_2(\tau, g)
   ,
 \medskip
 \\
   g(x,x,y) = y
   ,
\end{array}
\right.
$$
 и точка ${(\tau,\ g(\tau,x,y))}$
 содержится в $T$ при ${0 \le \tau \le x}$,
 если точка ${(x,y)}$ содержится в $T$.
 Проинтегрируем уравнение (3.23) вдоль характеристики
 ${\eta = g(\tau,x,y)}$ от $0$ до $x$. Получим
$$
\begin{array}{l}
  (u_1 - u_2)(x,y)
  =
 \medskip
 \\
 \phantom{12345}
  =
  (u_1 - u_2)
  (0, g(0,x,y))
  \times
 \medskip
 \\
 \phantom{12345}
  \times
  \exp
  \{
  {\ds \int_0^x}
  \left(
    2\alpha_2 +
    {\ds
      \frac
      {C \Delta_q - \Delta C_q}
      {2\Delta}
    }
   \right.
    +
 \medskip
 \\
 \phantom{12345}
    +
    (u_1 + u_2)
    \left(
      \beta_2
      +
      \frac{\ds B_q}{\ds 2}
      +
      {\ds
        \frac
        {C \Delta_p + \Delta C_p}
        {4 \Delta}
      }
    \right)
      +
 \medskip
 \\
 \phantom{12345}
    +
  \left.
      (u_2 - u_1)
      \left(
        \beta_1
        +
        {\ds
          \frac
          {C C_p + \Delta \Delta_p}
          {4 \Delta}
        }
      \right)
      -
      B_p u_1 u_2 + u_{2y}
  \right)
 \medskip
 \\
 \phantom{12345}
  (\tau, g(\tau, x,y))
  \: d\tau
  \}
  .
\end{array}
 \eqno(3.24)
$$
 Учитывая неравенство (3.22), из выражения (3.24) заключаем, что
 ${u_1 \neq u_2}$ в области $T$.

 Далее, положим
$$
\begin{array}{l}
  r(x,y,z,p,q)
  =
  \Delta
  {\ds
    \frac
    {u_1 u_2}
    {u_1 - u_2}
  }
  -
  D
  ,
 \medskip
 \\
  s(x,y,z,p,q)
  =
  {\ds
    \frac
    {\Delta}
    {2}
  }
  {\ds
    \frac
    {u_1 + u_2}
    {u_2 - u_1}
  }
  +
  {\ds
    \frac
    {C}
    {2}
  }
  ,
 \medskip
 \\
  t(x,y,z,p,q)
  =
  {\ds
    \frac
    {\Delta}
    {u_1 - u_2}
  }
  -
  B
  .
\end{array}
 \eqno(3.25)
$$
 Введенные обозначения позволяют записать уравнения (3.19)--(3.20)
 в виде
$$
  p_x + u_2 p_y = r + u_2 s
  ,
  \eqno(3.19')
$$
$$
  q_x + u_1 q_y = s + u_1 t
  .
  \eqno(3.20')
$$
 Так как ${u_1 \neq u_2}$ в области $T$, то уравнения (3.16)
 эквивалентны уравнениям (3.14)--(3.15). Последние можно записать
 следующим образом:
$$
  r_y + r_z q + r_p s + r_q t = s_x + s_z p + s_p r + s_q s
  ,
  \eqno(3.14')
$$
$$
  s_y + s_z q + s_p s + s_q t = t_x + t_z p + t_p r + t_q s
  .
  \eqno(3.15')
$$
 Характеристики системы уравнений
 (3.18), $(3.19')$, $(3.20')$
 находим, решая задачу Коши для обыкновенного дифференциального
 уравнения
$$
 \left\{
 \begin{array}{l}
   {\ds
    \frac{d g_i}{d \tau}
   }
   =
   u_{3-i} (\tau, g_i)
   ,
 \medskip
 \\
 g_i (x,x,y) = y \qquad (i = 1,2).
 \end{array}
 \right.
 \eqno(3.26)
$$
 Проинтегрируем уравнения
 (3.18), $(3.19')$, $(3.20')$
 вдоль соответствующих характеристик от $0$ до $x$. В результате
 получим
$$
 \begin{array}{l}
   z(x,y)
   =
   z_0(g_2(0,x,y))
   +
   {\ds \int_0^x}
   \{
     p + u_1 q
   \}
   (\tau, g_2(\tau, x,y))
   \: d\tau
   ,
 \medskip
 \\
   p(x,y)
   =
   p_0(g_1(0;x,y))
   +
   {\ds \int_0^x}
   \{
     r + u_2 s
   \}
   (\tau, g_1(\tau, x,y))
   \: d\tau
   ,
 \medskip
 \\
   q(x,y)
   =
   z_0^{'}(g_2(0,x,y))
   +
   {\ds \int_0^x}
   \{
     s + u_1 t
   \}
   (\tau, g_2(\tau, x,y))
   \: d\tau
   .
 \end{array}
 \eqno(3.27)
$$
 Правые части уравнений (3.26) непрерывно дифференцируемы. Поэтому
 функции
 ${g_i(\tau, x,y)}\ {(i=1,2)}$
 непрерывно дифференцируемы и имеют непрерывные смешанные
 производные второго порядка по $\tau,\ x$ и по $\tau,\ y$, а
 между производными функций $g_i$ по $x$ и по $y$ имеется
 зависимость
$$
   g_{ix}
 (\tau,x,y)
 +
 u_{3-i}
 (x,y)
   g_{iy}
 (\tau,x,y)
 =
 0
 .
 \eqno(3.28)
$$

 Вычислим первые производные функций $z,\ p$ и $q$. Для этого
 продифференцируем равенства (3.27). Принимая во внимание
 соотношения $(3.14')-(3.15')$, (3.28), интегрируя по частям и
 учитывая начальные значения (3.21), получаем
$$
 \begin{array}{l}
   z_y(x,y)
   =
   q(x,y)
   +
   {\ds \int_0^x}
   \{
     (p_y - s) + (s - q_x)
   \}
   (\tau, g_2(\tau, x,y))
    g_{2y}
    (\tau,x,y)
   \: d\tau
   ,
 \medskip
 \\
   z_x(x,y)
   =
   p(x,y)
   +
   {\ds \int_0^x}
   \{
     (p_y - s) + (s - q_x)
   \}
   (\tau, g_2(\tau, x,y))
    g_{2x}
    (\tau,x,y)
   \: d\tau
   ,
 \medskip
 \\
   p_y(x,y)
   =
   s(x,y,z(x,y),p(x,y),q(x,y))
   +
 \medskip
 \\
   \phantom{123}
   +
   {\ds \int_0^x}
   \{
     s_z(p - z_x) + r_z(z_y - q) + s_p(r - p_x) + r_p(p_y - s)
     +
 \medskip
 \\
   \phantom{123}
   +
     s_q(s - q_x) + r_q(q_y - t)
   \}
   (\tau, g_1(\tau, x,y))
     g_{1y}
    (\tau,x,y)
   \: d\tau
   ,
 \medskip
 \\
   p_x(x,y)
   =
   r(x,y,z(x,y),p(x,y),q(x,y))
   +
 \medskip
 \\
   \phantom{123}
   +
   {\ds \int_0^x}
   \{
     s_z(p - z_x) + r_z(z_y - q) + s_p(r - p_x) + r_p(p_y - s)
     +
 \medskip
 \\
   \phantom{123}
   +
     s_q(s - q_x) + r_q(q_y - t)
   \}
   (\tau, g_1(\tau, x,y))
     g_{1x}
    (\tau,x,y)
   \: d\tau
   ,
 \medskip
 \\
   q_y(x,y)
   =
   t(x,y,z(x,y),p(x,y),q(x,y))
   +
 \medskip
 \\
   \phantom{123}
   +
   {\ds \int_0^x}
   \{
     t_z(p - z_x) + s_z(z_y - q) + t_p(r - p_x) + s_p(p_y - s)
     +
 \medskip
 \\
   \phantom{123}
   +
     t_q(s - q_x) + s_q(q_y - t)
   \}
   (\tau, g_2(\tau, x,y))
    g_{2y}
    (\tau,x,y)
   \: d\tau
   ,
 \medskip
 \\
   q_x(x,y)
   =
   s(x,y,z(x,y),p(x,y),q(x,y))
   +
 \medskip
 \\
   \phantom{123}
   +
   {\ds \int_0^x}
   \{
     t_z(p - z_x) + s_z(z_y - q) + t_p(r - p_x) + s_p(p_y - s)
     +
 \medskip
 \\
   \phantom{123}
   +
     t_q(s - q_x) + s_q(q_y - t)
   \}
   (\tau, g_2(\tau, x,y))
    g_{2x}
    (\tau,x,y)
   \: d\tau
   .
 \end{array}
 \eqno(3.29)
$$
 Выражения (3.29) позволяют сделать вывод о том, что
$$
 \begin{array}{l}
   z_y(x,y)
   =
   q(x,y)
   ,
 \medskip
 \\
   z_x(x,y)
   =
   p(x,y)
   ,
 \medskip
 \\
   p_y(x,y)
   =
   q_x(x,y)
   =
   s(x,y,z(x,y),p(x,y),q(x,y))
   ,
 \medskip
 \\
   p_x(x,y)
   =
   r(x,y,z(x,y),p(x,y),q(x,y))
   ,
 \medskip
 \\
   q_y(x,y)
   =
   t(x,y,z(x,y),p(x,y),q(x,y))
   .
 \end{array}
$$
 С учетом соответствий (3.25) это означает, что три функции $z,\
 p$ и $q$ --- $C^1\mbox{-решение}$ системы (3.10)--(3.11). Отсюда
 согласно лемме~3.1 вытекает, что $z$ --- $C^3\mbox{-решение}$ задачи
 (3.1)--(3.3) в области $T$.
 $\Box$

 \textbf{Замечание.}
 Система (3.16)--(3.20) состоит из пяти уравнений и содержит пять
 неизвестных функций $u_1,\ u_2,\ z,\ p$ и $q$. В случае, если
 коэффициенты $A,\ B,\ C$ и $D$ уравнения (3.1) не зависят от
 значения $z$, уравнения (3.16), (3.19) и (3.20) образуют
 замкнутую систему четырех уравнений с четырьмя неизвестными
 функциями $u_1,\ u_2,\ p$ и $q$. Зная $u_1,\ u_2,\ p$ и $q$, мы
 можем найти $z$ из уравнения (3.18) или первого равенства (3.27).
 Определенный интерес представляет ситуация, при которой уравнения
 (3.16) образуют замкнутую систему с двумя неизвестными $u_1$ и
 $u_2$. Очевидно, что указанная ситуация имеет место тогда и
 только тогда, когда
$$
 {\ds
   \frac{\pt E_j}{\pt z}
   =
   \frac{\pt E_j}{\pt p}
   =
   \frac{\pt E_j}{\pt q}
   =
   \frac{\pt I_j}{\pt z}
   =
   \frac{\pt I_j}{\pt p}
   =
   \frac{\pt I_j}{\pt q}
   =
   0
 }
 \eqno(3.30)
$$
 ${(j = 0,\dots, 5)}$, где $E_j$ и $I_j$ вычислены согласно
 выражениям (3.17). В этом случае после нахождения решения
 ${u_1,\ u_2}$ значения $z,\ p$ и $q$ могут быть найдены из
 уравнений (3.18)--(3.20) или (3.11).



\subsection{Окончательный вид системы в инвариантах}

Система (3.16), полученная Д.~В.~Туницким \cite{Туницкий-деп},
имеет в правой части многочлены с образующими $u_1,\ u_2,$ а
функции
 $p,\ q,$ тоже искомые, содержатся в коэффициентах. Это может вызвать
трудности при попытке наложить условия на коэффициенты.
 В работах автора предложена другая запись
системы (3.16): искомые функции $u_1, u_2, p, q$ являются
образующими, а коэффициенты состоят из известных функций.

 Введем переобозначение
$$
   {r = u_1, \quad s = u_2.}
 \eqno(3.31)
$$
 Подчеркнем, что здесь $r,\ s$ являются характеристическими
 переменными (\cite{Ку}, прил.~1 к гл.~V, \S~2)
 (они же римановы инварианты и собственные значения
 системы (\cite{РЯ}, гл.~1)), а не обозначениями Монжа
 для вторых производных из
 (3.25).
 И те, и другие обозначения являются традиционными.

 Система имеет вид
$$
    \left(
      \partial_x
           +\xi(\omega) \partial_y
    \right)
    \omega = f_\omega(x,y,r,s,p,q,z),
 \eqno(3.32)
$$
где индекс $\omega$ принимает значения $r,s,p,q,z$,
 функция ${\xi(\omega)}$ имеет вид
$$
 \xi(r)=s, \quad
 \xi(s)=r, \quad
 \xi(p)=s, \quad
 \xi(q)=r, \quad
 \xi(z)=r \, ,
$$
$$
\begin{array}{l}
f_r=\rho_0+\rho_1r+\rho_2s+\rho_3pr+\rho_4qr+\rho_5ps+\rho_6qs+\rho_7r^2
+\rho_8rs+\rho_9pr^2+
\smallskip \\
\phantom{f^r}+\rho_{10}qr^2+\rho_{11}prs+\rho_{12}qrs+\rho_{13}r^2s
=f_r(\rho ,r,s,p,q),
 \smallskip \\
 f_s(\sigma ,r,s,p,q)=f_r(\sigma ,s,r,p,q),
 \smallskip \\
 f_p=\pi_0+\pi_1 s,\quad
 f_q=\kappa_0+\kappa_1 r,\quad
 f_z=p+qr,
\end{array}
$$
вектор-функции $\rho ,\sigma ,\pi ,\kappa$ зависят от $x,y,z,p,q$,
$$
\begin{array}{ll}
\rho_{0\phantom{0}}=D_q\: , & \sigma_{0\phantom{0}}=D_q\: ,
\smallskip \\
 \rho_{3\phantom{0}}=\phantom{-}\frac{1}{2\Dt}(C+\Dt)_z\:
,&
            \sigma_{3\phantom{0}}=-\frac{1}{2\Dt}(C-\Dt)_z\: ,
\smallskip \\
\rho_{4\phantom{0}}=\phantom{-}\frac{1}{\Dt}D_z\: ,&
            \sigma_{4\phantom{0}}=-\frac{1}{\Dt}D_z\: ,
\smallskip \\
\rho_{5\phantom{0}}=-\frac{1}{2\Dt}(C+\Dt)_z\: ,&
            \sigma_{5\phantom{0}}=\phantom{-}\frac{1}{2\Dt}(C-\Dt)_z\: ,
\smallskip \\
\rho_{6\phantom{0}}=-\frac{1}{\Dt}D_z\: ,&
            \sigma_{6\phantom{0}}=\phantom{-}\frac{1}{\Dt}D_z\: ,
\smallskip \\
\rho_{9\phantom{0}}=-\frac{1}{\Dt}B_z\: ,&
            \sigma_{9\phantom{0}}=\phantom{-}\frac{1}{\Dt}B_z\: ,
\smallskip \\
\rho_{10}=-\frac{1}{2\Dt}(C-\Dt)_z\: ,\phantom{1234567}&
                      \sigma_{10}=\phantom{-}\frac{1}{2\Dt}(C+\Dt)_z\: ,
\smallskip \\
\rho_{11}=\phantom{-}\frac{1}{\Dt}B_z\: ,&
                      \sigma_{11}=-\frac{1}{\Dt}B_z\: ,
\smallskip \\
\rho_{12}=\phantom{-}\frac{1}{2\Dt}(C-\Dt)_z\: ,&
                      \sigma_{12}=-\frac{1}{2\Dt}(C+\Dt)_z\: ,
\smallskip \\
\rho_{13}=-B_p\: ,&\sigma_{13}=-B_p\: ,
\end{array}
$$
$$
\begin{array}{l}
\rho_1=\phantom{-}\frac{1}{2\Dt} ( \phantom{-}
\Dt_x+C_x+2D_y+CD_p-D\Dt_p-DC_p-2BD_q+\frac{1}{2}CC_q-\frac{3}{2}\Dt C_q+
\smallskip \\
\phantom{\rho_1=} +\frac{1}{2}C\Dt_q+\frac{1}{2}\Dt\Dt_q-\Dt D_p
)\: ,
\smallskip \\
\rho_2=-\frac{1}{2\Dt} ( \phantom{-}
\Dt_x+C_x+2D_y+CD_p-D\Dt_p-DC_p-2BD_q+\frac{1}{2}CC_q+\frac{1}{2}\Dt C_q+
\smallskip \\
\phantom{\rho_1=} +\frac{1}{2}C\Dt_q+\frac{1}{2}\Dt\Dt_q+\Dt D_p
)\: ,
\smallskip \\
\rho_7=-\frac{1}{2\Dt} (
-\Dt_y+C_y+2B_x+CB_q+B\Dt_q-BC_q-2DB_p+\frac{1}{2}CC_p-\frac{1}{2}\Dt C_p-
\smallskip \\
\phantom{\rho_1=} -\frac{1}{2}C\Dt_p+\frac{1}{2}\Dt\Dt_p-\Dt B_q
)\: ,
\smallskip \\
\rho_8=\phantom{-}\frac{1}{2\Dt} (
-\Dt_y+C_y+2B_x+CB_q+B\Dt_q-BC_q-2DB_p+\frac{1}{2}CC_p+\frac{3}{2}\Dt C_p-
\smallskip \\
\phantom{\rho_1=} -\frac{1}{2}C\Dt_p+\frac{1}{2}\Dt\Dt_p+\Dt B_q
)\: ,\smallskip \\
\sigma_1=-\frac{1}{2\Dt} (
-\Dt_x+C_x+2D_y+CD_p+D\Dt_p-DC_p-2BD_q+\frac{1}{2}CC_q+\frac{3}{2}\Dt C_q-
\smallskip \\
\phantom{\rho_1=} -\frac{1}{2}C\Dt_q+\frac{1}{2}\Dt\Dt_q+\Dt D_p
)\: ,
\smallskip \\
\sigma_2=\phantom{-}\frac{1}{2\Dt} (
-\Dt_x+C_x+2D_y+CD_p+D\Dt_p-DC_p-2BD_q+\frac{1}{2}CC_q-\frac{1}{2}\Dt C_q-
\smallskip \\
\phantom{\rho_1=} -\frac{1}{2}C\Dt_q+\frac{1}{2}\Dt\Dt_q-\Dt D_p
)\: ,
\end{array}
$$
$$
\begin{array}{l}
\sigma_7=\phantom{-}\frac{1}{2\Dt} ( \phantom{-}
\Dt_y+C_y+2B_x+CB_q-B\Dt_q-BC_q-2DB_p+\frac{1}{2}CC_p+\frac{1}{2}\Dt C_p+
\smallskip \\
\phantom{\rho_1=} +\frac{1}{2}C\Dt_p+\frac{1}{2}\Dt\Dt_p+\Dt B_q
)\: ,
\smallskip \\
\sigma_8=-\frac{1}{2\Dt} ( \phantom{-}
\Dt_y+C_y+2B_x+CB_q-B\Dt_q-BC_q-2DB_p+\frac{1}{2}CC_p-\frac{3}{2}\Dt
C_p+
\smallskip \\
\phantom{\rho_1=}
+\frac{1}{2}C\Dt_p+\frac{1}{2}\Dt\Dt_p-\Dt B_q )\: ,
\end{array}
$$
$$
\begin{array}{ll}
\pi_0=-D\: ,                  & \kappa_0=\frac{1}{2}(C+\Dt) \: ,
\\
\pi_1=\frac{1}{2}(C+\Dt)\: ,\phantom{1234567} & \kappa_1=-B \: .
\end{array}
$$

 Начальные условия (3.21) для системы (3.32)
 с учетом переобозначений (3.31)
 имеют вид
$$
\begin{array}{l}
  r(0,y)=r^0(y)=\frac{\ds (C+\Dt)(0,y,z^0(y),p^0(y),z^0_y(y))-2p^0_y(y)}
           {\ds 2(z^0_{yy}(y)+B(0,y,z^0(y),p^0(y),z^0_y(y)))}\; ,
  \medskip \\
  s(0,y)=s^0(y)=\frac{\ds (C-\Dt)(0,y,z^0(y),p^0(y),z^0_y(y))-2p^0_y(y)}
           {\ds 2(z^0_{yy}(y)+B(0,y,z^0(y),p^0(y),z^0_y(y)))}\; ,\medskip \\
  p(0,y)=p^0(y)\; ,\quad
  q(0,y)=q^0(y)=z^0_y(y)\; ,\quad
  z(0,y)=z^0(y)\; .
\end{array}
 \eqno(3.32^0)
$$

Запишем систему $(3.32^0)$ в виде
$$
 \left\{
 \begin{array}{l}
 p_y^0 = - \frac{\ds \Delta}{\ds 2}
 \frac{\ds r^0 + s^0}{\ds r^0 - s^0} +
 \frac{\ds C}{\ds 2} \, ,
 \medskip
 \\
 q_{y}^0 =
 \frac{\ds \Delta}{\ds r^0 - s^0} - B \, ,
 \medskip
 \\
 z_y^0 = q^0 \, .
 \end{array}
 \right.
 \eqno (3.32^{00})
$$
Удобно задавать в качестве начальных данных функции $r^0,\ s^0$, а
функции $z^0,\ p^0$ находить из системы $(3.32^{00})$. Если
коэффициенты $A,B,C,D$ зависят только от переменных $x,\ y$,
система является линейной. В общем случае она нелинейна.

 Уравнение с коэффициентами, зависящими только от $x,\ y$, сводится к
 системе
$$
\begin{array}{l}
    ( \partial_x + s \, \partial_y )
    r
    =
    (r - s)(a_{r 1}+a_{r 2}\, r)\, ,
  \medskip \\
    ( \partial_x + r \, \partial_y )
    s
    =
    (r - s)(a_{s 1}+a_{s 2}\, s)\, ,
\end{array}
 \eqno (3.33)
$$
где
$$
\begin{array}{ll}
    a_{r1}=\frac{\ds 1}{\ds 2\Dt}(2D_y+C_x+{\Dt}_x)\, , &
    a_{r2}=\frac{\ds 1}{\ds 2\Dt}(-2B_x-C_y+{\Dt}_y)\, ,
    \medskip
    \\
    a_{s1}=\frac{\ds 1}{\ds 2\Dt}(2D_y+C_x-{\Dt}_x)\, , &
    a_{s2}=\frac{\ds 1}{\ds 2\Dt}(-2B_x-C_y-{\Dt}_y) \, .
\end{array}
$$
 После нахождения функций $r,\ s$ решаются
 линейные уравнения (3.32) для функций $p, \ q, \ z$.
 Начальные данные $p^0,\ q^0,\ z^0$ в этом случае выражаются через
 $r^0,\ s^0$ при помощи $(3.32^{00})$.



\section{Последовательность приближенных решений}
\subsection{Итерационная последовательность}

   Пусть известны функции ${\up{n}{\omega}}(x,y),\;
\omega=r,s,p,q,z.$ Функции ${\up{n+1}{\omega}}(x,y),\linebreak
\omega=r,s,p,q,z,$ определим как решение нелинейной задачи Коши
$$
\left\{
  \begin{array}{l}
    (
        \pt_x
        +
        {\up{n+1}{\xi}}(\omega)
        \pt_y
    )
    \up{n+1}{\omega}
    =
    f_\omega(x,y,\up{n}{r},\up{n}{s},\up{n}{p},\up{n}{q},\up{n}{z}),
    \smallskip \\
    {\up{n+1}{\omega}}
    (0,y)=\omega^0(y),
    \quad \omega=r,s,p,q,z.
  \end{array}
\right.
 \eqno(4.1)
$$

 Задача Коши (4.1) распадается на четыре
независимые задачи: нелинейная система для
 ${\up{n+1}{r}}$, ${\up{n+1}{s}}$
 и три независимых линейных уравнения для оставшихся функций.
 Доказательству разрешимости нелинейной системы мы предпошлем
 установление априорных оценок.

   Начальное приближение определим так:
$$
  {{\up{0}{\omega}}(x,y)}
  =
 \omega^0(y), \quad \omega=r,s,p,q,z.
 \eqno(4.2)
$$

   $C^1\mbox{-решение}$ задачи (3.32), $(3.32^0)$ существует, если существуют
непрерывные функции $\omega,\ \omega_y, \; {\omega=r,s,p,q,z}$,
так как согласно уравнениям (3.32) вектор-функция $\omega_x$
выражается через $\omega,\ \omega_y$. Доказательство существования
непрерывных вектор-функций $\omega,\ \omega_y$ есть доказательство
равномерной сходимости последовательностей
 $\{{\up{n}{\omega}}\},\ \{{\up{n}{\omega}_y}\}, \ {\omega=r,s,p,q,z,}$
 и в этом
доказательстве решающую роль играет равномерная ограниченность
указанных последовательностей.


\subsection{Равномерная ограниченность}

   Пусть вектор-функции
 $\rho, \sigma, \pi, \kappa$
 и функции
  $r^0,\: s^0,\: p^0,\: q^0$
 ${C^1\mbox{-гладки}}$ и ограничены,
 ${z^0 \in C^1(\mathbb{R})}$.

 В выражениях для функций ${f_r,\, f_s}$
 из правых частей (3.32)
 мономы разбиваются на две категории:
 содержащие ${p,\ q}$ и не содержащие.
 Соответственно требуется различать
 коэффициенты
 ${\rho, \, \sigma}$
 из разных категорий.
 Для различения введем множества индексов:
 ${J_{rs} = \{ 0,1,2,7,8,13 \} }$
 для мономов, не содержащих ${p,\ q,}$
 и
 ${J_{pq} = \{ 3,4,5,6,9,10,11,12 \} }$
 --- для содержащих.

   Введем обозначения
$$
\begin{array}{rcl}
    U_0&=&\max\lms_{\omega=r,s}\:\sup\lms_{y\in \mathbb{R}}
    |\omega^0(y)|
    = \hbox{const},
 \medskip\\
    \alpha_1(x)&=&
      \sup\lms_
      {
        {
            {
            \atop\scriptstyle (y,z,p,q)\in \mathbb{R}^4
            }
            \atop
            {
            \atop\scriptstyle j \, \in \, J_{rs}
            }
        }
        \atop
        {
          \atop\scriptstyle a \, \in \, \{ \rho, \sigma \}
        }
      }
      |a_j(x,y,z,p,q)|,
 \medskip\\
    \alpha_2(x)&=&
      \sup\lms_
      {
        {
            {
            \atop\scriptstyle (y,z,p,q)\in \mathbb{R}^4
            }
            \atop
            {
            \atop\scriptstyle j \, \in \, J_{pq}
            }
        }
        \atop
        {
          \atop\scriptstyle a \, \in \, \{ \rho, \sigma \}
        }
      }
      |a_j(x,y,z,p,q)|,
 \medskip\\
    \alpha_3 &=&
      \sup\lms_
      {
        {
            {
            \atop\scriptstyle (x,y,z,p,q)\in \mathbb{R}^5
            }
            \atop
            {
            \atop\scriptstyle j=0,1
            }
        }
        \atop
        {
          \atop\scriptstyle a \, \in \, \{ \pi, \kappa \}
        }
      }
      |a_j(x,y,z,p,q)|
      = \hbox{const}.
\end{array}
 \eqno(4.3)
$$

{\bf Лемма~4.1.} Пусть
$$
\begin{array}{l}
  |p^0|\le 1
   \; , \quad
   |q^0|\le 1
   \; ,
  \bigskip \\
  U_0
  +
  6{\ds\int_{-\infty}^{+\infty}}
    \alpha_1(x)
  \: dx
  +
  8{\ds\int_{-\infty}^{+\infty}}
    (1 + 2\alpha_3 |x|)\alpha_2(x)
  \: dx
  \le 1
  \; .
\end{array}
  \eqno(4.4)
$$
Предположим, что для некоторого номера ${n \ge 0}$ имеют место
оценки
$$
\begin{array}{ll}
 |\up{n}{\omega}(x,y)|
 \le
 1
  \: ,
 &
 \omega=r,s
 \; ,
 \smallskip
 \\
 |\up{n}{\omega}(x,y)|
 \le
 1 + 2\alpha_3 |x|
 \: ,
 &
 \omega=p,q
 \; ,
 \smallskip
 \\
 |\up{n}{z}(x,y)|
 \le
   \max\lms_{t \in [y - |x|,\ y + |x|]}
   |z^0(t)|
   + 2|x| + 2\alpha_3 x^2
\end{array}
  \eqno(4.5)
$$
для ${\forall \: (x,y)\in [0,+\infty )\times \mathbb{R}}$.
 Тогда для номера ${n+1}$
 при условии существования функций
 ${\up{n+1}{\omega}},\ \omega = r,s,p,q,z$,
 справедливы те же оценки (4.5).

{\bf Доказательство.} Из соотношений (4.2), (4.4) следует, что для
${n=0}$ неравенства (4.5) верны. Предположим, что неравенства
(4.5) выполнены для некоторого номера $n$, и докажем их для номера
${n+1}$.

 Воспользуемся представлением правых частей $f_\omega$ в системе (3.32).
Из формулы (2.10) получаем: для $\omega=r,s$
$$
\begin{array}{rl}
  |{\up{n+1}{\omega}}(x,y)|
  &
  \le
  |\omega^0|
  +
  {\ds \int_0^x}
  |f_\omega|
  \: d\tau
  \le
 \medskip
 \\
  &
 \le
 U_0
 +
 6
 {\ds \int_0^{+\infty}}
 \!\!\!
   \alpha_1(\tau)
 \: d\tau
 +
 4
 {\ds \int_0^{+\infty}}
 \!\!\!
 \alpha_2(\tau)
 \: |{\up{n}{p}}|
 \: d\tau
 +
 4
 {\ds \int_0^{+\infty}}
 \!\!\!
 \alpha_2(\tau)
 \: |{\up{n}{q}}|
 \: d\tau
 \le
 \medskip
 \\
  &
 \le
 U_0
 +
 6
 {\ds \int_0^{+\infty}}
 \!\!\!
   \alpha_1(\tau)
 \: d\tau
 +
 8
 {\ds \int_0^{+\infty}}
 \!\!\!
 \alpha_2 (\tau)
 (1 + 2\alpha_3 \tau)
 \: d\tau
 \le
 1
 \: .
\end{array}
$$
Для $\omega=p,q,\; a=\pi,\kappa$
$$
\begin{array}{rl}
  |{\up{n+1}{\omega}}(x,y)|
  &
  \le
  |\omega^0|
  +
  {\ds \int_0^x}
  |f_\omega|
  \, d\tau
  \le
  \phantom{12}
    \langle\mbox{см. (4.3) для }\alpha_3\rangle
 \medskip\\
  &
 \le
 1
 +
 {\ds \int_0^x}
 (|a_0|+|a_1|)
 \, d\tau
 \le
 1
 +
 2\alpha_3
 {\ds \int_0^x}
 d\tau
 \le
 1
 +
 2\alpha_3
 \:
 |x|
 \: .
\end{array}
$$
 Для ${\omega=z,}$ так как ${|{\up{n}{r}}|\le 1,}$ т.е.
 ${{\up{n+1}{g_z}}(\tau,x,y)\in [y - x + \tau, \ \ y + x - \tau]:}$
$$
\begin{array}{l}
  |{\up{n+1}{z}}(x,y)|
  \le
  |z^0|
  +
  {\ds \int_0^x}
  |\up{n}{p}+\up{n}{q}\up{n}{r}|
  \, d\tau
  \le
 \max\lms_{t\in [y - |x|,\ y + |x|]} |z^0(t)|
 +
 \medskip\\
 +
 {\ds \int_0^x}
     2(1 + 2\alpha_3 \tau)
 \, d\tau
 =
 \max\lms_{t\in [y - |x|,\ y + |x|]} |z^0(t)|
 +
 2 |x|
 +
 2\alpha_3 x^2
 \: .
 \ \Box
\end{array}
$$

{\bf Следствие.} При выполнении условий (4.4) семейства
${\{\up{n}{\omega}\}}\: ,\; {\omega=r,s,p,q,z}$, равномерно
ограничены на компакте
$$
  G(\bar x,\bar y)=\{(x,y)|\; x\in[0,\bar x],\;
   y\in[\bar y - \bar x + x,\ \ \bar y + \bar x - x]\}
   \eqno (4.6)
$$
 для
 ${\forall \: (\bar x,\bar y)\in [0,+\infty)\times\mathbb{R}}$.

{\bf Доказательство.} Функции $\up{n}{\omega},\;
{\omega=r,s,p,q,z}$, определены на компакте ${G(\bar x,\bar y)}$,
так как характеристики, выходящие из любой точки ${(x,y)\in G(\bar
x,\bar y)}$, лежат в компакте ${G(\bar x,\bar y)}$ в силу
${|{\up{n}{r}}|\le 1,}\; {|{\up{n}{s}}|\le 1}$. Равномерная
ограниченность следует из неравенств (4.5). $\Box$


\subsection{Гиперболичность в узком смысле}

Потребуем, чтобы начальные данные $r^0, s^0$ были разделены
некоторой постоянной. Более точно: $\exists \, \delta > 0$:
$$
  \inf \lms_{\scriptstyle y\in \mathbb{R}} r^0(y)
  -
  \sup \lms_{\scriptstyle y\in \mathbb{R}} s^0(y)
  \geq \delta > 0 \, .
  \eqno (4.7)
$$

{\bf Лемма~4.2.} Пусть $\exists \, \varepsilon \in (0,\delta]$:
$$
  6
  {\ds\int_{-\infty}^{+\infty}}
      \alpha_1(x)
  \: dx
  +
  8
  {\ds\int_{-\infty}^{+\infty}}
    (1 + 2\alpha_3 |x|)
    \alpha_2(x)
  \: dx
  \leq
  \frac{\delta - \varepsilon}{2}
  \, .
  \eqno (4.8)
$$
Тогда для $n=0,1,2,\dots$
$$
\inf \lms_{\scriptstyle (x,y)\in \mathbb{R}^2} \up{n}{r}(x,y) -
  \sup \lms_{\scriptstyle (x,y)\in \mathbb{R}^2} \up{n}{s}(x,y)
  \geq \varepsilon > 0 .
  \eqno (4.9)
$$

{\bf Доказательство.} Согласно формуле (2.10), учитывая оценки
(4.3), (4.5), (4.8), имеем
$$
\begin{array}{rl}
 {\up{n+1}{r}}(x,y)
 &
 \ge
 \inf \lms_{\scriptstyle y\in \mathbb{R}} r^0(y)
 -
 \left(
  6
  {\ds\int_{-\infty}^{+\infty}}
    \alpha_1(x)
  \: dx
  +
  8
  {\ds\int_{-\infty}^{+\infty}}
    (1 + 2\alpha_3 |x|)
    \alpha_2(x)
  \: dx
 \right)
 \ge
 \medskip
 \\
 &
 \ge
 \inf \lms_{\scriptstyle y\in \mathbb{R}} r^0(y)
 - \frac{\ds \delta - \varepsilon}{\ds 2} \, ,
 \medskip
  \\
 {\up{n+1}{s}}(x,y)
 &
 \le
 \sup \lms_{\scriptstyle y\in \mathbb{R}} s^0(y)
 +
 \left(
  6{\ds\int_{-\infty}^{+\infty}}\alpha_1(x)dx
  +
  8
  {\ds\int_{-\infty}^{+\infty}}
  (1 + 2\alpha_3 |x|)
  \alpha_2(x)
  \: dx
 \right)
 \leq
 \medskip
 \\
 &
 \leq
 \sup \lms_{\scriptstyle y\in \mathbb{R}}
 s^0(y)
 +
 \frac{\ds \delta - \varepsilon}{\ds 2}
 \, .
\end{array}
$$
Тогда согласно условию (4.7)
$$
 \inf \lms_{\scriptstyle (x,y)\in \mathbb{R}^2}
 {\up{n+1}{r}}(x,y)
 -
 \sup \lms_{\scriptstyle (x,y)\in \mathbb{R}^2}
 {\up{n+1}{s}}(x,y)
 \geq
 \inf \lms_{\scriptstyle y\in \mathbb{R}} r^0(y)
 -
 \sup \lms_{\scriptstyle y\in \mathbb{R}} s^0(y)
 - (\delta - \varepsilon) \geq \varepsilon
 \; .
 \;\square
$$


\subsection{Разрешимость итерационной системы}

\textbf{Лемма~4.3.} ${C^1\mbox{-решение}}$ задачи (4.1) существует
на всей полуплоскости.
\\
 \textbf{Доказательство.} Рассмотрим
 подсистему двух уравнений относительно
 $\up{n+1}{r}$, $\up{n+1}{s}$. Она является слабо-нелинейной
 (\cite{РЯ}, гл.~1, ${\S}~10,$ п.~3).
 (Для системы (2.8) слабая нелинейность --- это
  ${\ds {\pt \xi_i}/{\pt u_i} = 0}$ для всех $i$.)
 Решение этой системы ограничено на
 всей полуплоскости
 согласно оценкам (4.5). Система является
 гиперболической в узком смысле согласно оценкам (4.9).
 Согласно теореме Рождественского --- Сидоренко
 (\cite{РЯ}, гл.~1, ${\S}~10,$ п.~3)
 первые производные функций
 $\up{n+1}{r}$, $\up{n+1}{s}$
 не обращаются в бесконечность при конечном значении $x$.
 Следовательно, согласно следствию из этой
 теоремы
 (\cite{РЯ}, гл.~1, ${\S}~10,$ п.~3),
 подсистема первых двух уравнений
 задачи (4.1) разрешима на всей
 (полу)плоскости, т.е. имеет глобальное $C^1\mbox{-гладкое}$ решение.

 Действительно, задача (4.1) разрешима локально, т.е. в
 некоей окрестности прямой ${x = 0}$. Этот хорошо известный факт
 следует, например, из теоремы существования
 (\cite{РЯ}, гл.~1, ${\S}~8,$ п.~2)
 с учетом следствия из леммы~4.1.
 При продолжении локального
 решения сильный разрыв, т.е. уход первой производной в бесконечность при
 конечном $x,$ невозможен согласно
 теореме Рождественского --- Сидоренко.
 В доказательстве этой теоремы строится мажоранта, ограничивающая
 классическое решение (по модулю), и, таким образом, перед уходом
 в бесконечность решение, оставаясь классическим, должно выйти за
 мажоранту, чего быть не может.

  Далее, как следует из изучения
 слабых разрывов (разрывов первого рода, т.е. конечных скачков первой производной),
 \textit{слабый разрыв
 системы квазилинейных уравнений гиперболического типа, распространяясь
 вдоль характеристики, не может ни возникнуть, ни исчезнуть,
 если только решение и его первые производные остаются
 ограниченными}
 (\cite{РЯ}, гл.~1, ${\S}~10,$ п.~1). Поскольку у нас начальные
 данные $C^1\mbox{-гладкие},$ слабых разрывов нет. Таким образом,
 $C^1\mbox{-гладкое}$ решение существует на всей полуплоскости.

 Оставшиеся три уравнения системы (4.1)
 являются линейными и,
 снабженные начальными данными (4.1),
 имеют решения
 на всей полуплоскости. $\Box$

 \textbf{Лемма~4.4.}
 Пусть
 ${{\up{n+1}{\omega}}(x,y)}$ --- решение
 задачи (4.1).
 Тогда
 характеристики задачи (4.1)
 существуют для любого значения параметра
 $\tau \in [0,x]$.

 \textbf{Доказательство.}
 Характеристики определяются
 как решения задачи (2.9).
 Классическая теорема Коши существования и
 единственности для систем вида
$$
 \frac{d}{d \tau}
 \ u_i
 =
 f_i (\tau , u_1, \dots , u_n),
 \quad i=1, \dots , n
$$
имеет следующий вид:

 {\it
 Пусть в замкнутой области
$$
 \bar G =
 \{
   (\tau , u_1, \dots , u_n)
   | \quad
   |\tau - \tau^0| \leq a, \quad
   |u_i - u_i^0| \leq b, \quad
   i=1, \dots , n
 \}
$$
 выполнены условия:\\
 1) функции $f_i$ непрерывны,\\
 2) ${|f_i| \leq A}$;\\
 3) функции $f_i$ липшицевы по переменным ${u_1, \dots , u_n}$.\\
 Тогда система c начальными условиями ${u_i(\tau^0) = u^0}$ имеет
 единственное решение для ${|\tau - \tau^o| \leq \min (a, b/A)}$.
 }

 В нашем случае ${A = 1}$ согласно
 равномерной оценке (4.5).
 Далее, ${a = x}$, и мы можем выбрать сколь угодно
 большое значение $b$. Липшицевость обеспечивается ограниченностью
 в области $\bar G$ частных производных по ${u_1, \dots , u_n}$: согласно
 лемме~4.3 правые части принадлежат к классу
 ${C^1([0, +\infty) \times \mathbb{R})}$
 как классическое решение задачи (4.1).
 Следовательно, обе характеристики существуют для
 ${\forall \: \tau \in [0, x]}$.
 Таким образом, мы можем
 интегрировать вдоль характеристик.
 $\ \Box$



\section{Производные приближенных решений}

  Распространим теорему Рождественского -- Сидоренко
 (\cite{РЯ}, гл.~1, ${\S}~10,$ п.~3), доказанную для
точного решения в случае системы двух уравнений, на метод
последовательных приближений для системы пяти уравнений с двумя
различными характеристиками.

 Запись
${{(\up{n}{u})}_y}$ означает, что к $u$ сначала применяется
оператор взятия ${n\mbox{-го}}$ приближения, и лишь затем оператор
дифференцирования по $y,$ а не наоборот.

Пусть ${G(\bar x, \bar y)}$ --- компакт, определенный условием
(4.6),
 ${\forall (\bar x,\bar y)\in [0,+\infty )\times \mathbb{R}}$.

 {\bf Лемма~5.1.} Существует функция $\Phi (x) \in C^0(\mathbb{R})$
такая, что
$$
|{(\up{n}{\omega})_y} (x,y)| \leq \Phi (x)
 \eqno (5.1)
$$
для $\omega = r,s,p,q,z$, $\forall (\bar x,\bar y)\in G(\bar x,
\bar y)$, $n=0,1,2,\dots$

 {\bf Доказательство.}
  Введем обозначение ${\up{n}{g}_{\omega}}
 (x,y_0) = {\up{n}{g}_{\omega}}
 (x,0,y_0)$. Тогда ${\up{n}{g}_{\omega}} (x,y_0)$ есть решение
 задачи Коши
 $$
\left\{
 \begin{array}{l}
    \pt_x
    {\up{n}{g}}_{\omega}
    =
    \xi ({\up{n}{\omega}})
    (x, \up{n}{g}_{\omega}),
 \smallskip \\
    {\up{n}{g}_{\omega}}(0,y_0) = y_0, \quad \omega = r,s,p,q,z,
 \end{array}
 \right.
 \eqno (5.2)
 $$
 т.е. кривая $(x, {\up{n}{g}_{\omega}}(x,y_0))$ --- характеристика,
 проходящая через точку $(0, y_0)$. Формула (2.10) принимает вид
$$
{\up{n+1}{\omega}} (x, {\up{n+1}{g_\omega}} (x,y_0)) =
 \omega^0 (y_0) + \int_0^x f_{\omega} (\tau, {\up{n+1}{g_\omega}} (\tau,y_0),
  {\up{n}{\Omega}} (\tau, {\up{n+1}{g_\omega}} (\tau,y_0)) )\; d\tau
  \, ,
$$
где $\Omega = (\omega) = (r,s,p,q,z)$.
 Введя обозначение
$$
{\up{n+1}{\bar\omega}} (x,y_0) = {\up{n+1}{\omega}} (x,
{\up{n+1}{g_\omega}} (x,y_0)) \, ,
 \eqno (5.3)
$$
получим для ${\omega = r,s,p,q,z}$:
$$
 {\up{n+1}{\bar\omega}} (x, {\up{n+1}{g_\omega}} (x,y_0)) =
 \omega^0 (y_0) + {\ds \int}_0^x f_{\omega} (\tau, {\up{n+1}{g_\omega}} (\tau,y_0),
  {\up{n}{\Omega}} (\tau, {\up{n+1}{g_\omega}} (\tau,y_0)) )\; d\tau
  \, .
 \eqno (5.4)
$$
Отсюда следует, что ${\up{n+1}{\bar\omega}} (x,y_0)$ есть решение
задачи Коши
$$
\left\{
\begin{array}{l}
 \pt_x {\up{n+1}{\bar\omega}} =
  f_{\omega} (x, {\up{n+1}{g_\omega}} (x,y_0),
  {\up{n}{\Omega}} (x, {\up{n+1}{g_\omega}} (x,y_0)) )\, ,
  \smallskip \\
 {\up{n+1}{\bar\omega}}(0,y_0) = \omega^0 (y_0), \quad
 \omega = r,s,p,q,z.
\end{array}
\right.
 \eqno (5.5)
$$
Дифференцированием по $y_0$ равенства (5.3) получим
$$
 \pt_{y_0} {\up{n+1}{\bar\omega}} (x, y_0) =
 \left.
 \pt_y {\up{n+1}{\bar\omega}} (x, y)
 \right|_{y = {\up{n+1}{g_\omega}} (x,y_0)}
 \pt_{y_0} {\up{n+1}{g_\omega}} (x,y_0) \, ,
$$
откуда
$$
 \pt_y {\up{n+1}{\omega}} (x, y) =
 \pt_{y_0}{\up{n+1}{\bar\omega}} (x, y_0)
 \pt_{y_0}
 {\up{n+1}{g_\omega}} (x,y_0) \, .
 \eqno (5.6)
$$
Дифференцируя по $y_0$ формулы (5.2), получим уравнение
$$
\pt_x
 \left(\pt_{y_0}
 {\up{n}{g}_{\omega}} (x,y_0)
 \right) =
 \left.
 \pt_y \xi ({\up{n}{\omega}}) (x,y)
 \right|_{y = {\up{n}{g}_{\omega}} (x,y_0)}
 \pt_{y_0}
 {\up{n}{g}_{\omega}} (x,y_0) \, ,
$$
или
$$
 \pt_x
 \ln
 \left(
 \pt_{y_0}
 {\up{n}{g}_{\omega}} (x,y_0)
 \right)
 =
 \left.
 \pt_y \xi ({\up{n}{\omega}}) (x,y)
 \right|_{y = {\up{n}{g}_{\omega}} (x,y_0)}
 \, ,
\eqno (5.7)
$$
и начальное условие
$$
 (
 \pt_{y_0}
 {\up{n}{g}_{\omega}}
 )
  (0,y_0) = 1 \, .
 \eqno (5.7^0)
$$
Пусть $v(x,y)\in C^1$ --- произвольная функция. Введем обозначение
$$
 \left(
     \frac{d}{dx}
     \: v
 \right)_{\omega}
 =
 (
    \pt_x
    +
    {\up{n+1}{\xi}} (\omega)
    \pt_y
 )
 \: v
 \, .
$$
Вычитая из равенства
$$
 (
    \pt_x
    +
    {\up{n+1}{r}} \pt_y
 )
 {\up{n+1}{s}}
 =
 f_s(x,y,{\up{n}{\Omega}})
$$
равенство
$$
 (
    \pt_x
    +
    {\up{n+1}{s}} \pt_y
 )
 {\up{n+1}{s}}
 =
 \left(
     \frac{d}{dx} {\up{n+1}{s}}
 \right)_r \, ,
$$
получим
$$
 \pt_y
     {\up{n+1}{s}} =
 \frac{f_s(x,y,{\up{n}{\Omega}})
         -
         (\frac{\ds d}{\ds dx}{\up{n+1}{s}})_r
      }
      {\ds {\up{n+1}{r}} - {\up{n+1}{s}}} \, .
$$
Преобразуем это равенство:
$$
\begin{array}{rl}
 \pt_y
 {\ds\up{n+1}{s}}
 &
 =
 \
 \frac{\ds f_s(x,y,{\up{n}{\Omega}})
           -
           f_r(x,y,{\up{n}{\Omega}})
      }
      {\ds {\up{n+1}{r}} - {\up{n+1}{s}}}
 +
 \frac{{\ds f_r(x,y,{\up{n}{\Omega}})} \, - \,
 {\ds (}
 \frac{\ds d}{\ds dx}{\ds\up{n+1}{s}}
        {\ds )}_r
      }
      {\ds {\up{n+1}{r}} - {\up{n+1}{s}}}
 =
  \smallskip \\
 &
 =
 \
 \frac{\ds f_s(x,y,{\up{n}{\Omega}})
           -
           f_r(x,y,{\up{n}{\Omega}})
      }{\ds {\up{n+1}{r}} - {\up{n+1}{s}}}
 +
 \frac{
        {\ds (}
           \frac{\ds d}{\ds dx}
           {\ds \up{n+1}{r}}
        {\ds )}_r
   \, - \,
 {\ds (}
      \frac{\ds d}{\ds dx}
      {\ds \up{n+1}{s}}
 {\ds )}_r
 }
 {\ds {\up{n+1}{r}} - {\up{n+1}{s}}}
 =
  \smallskip \\
 &
 =
 \
 \frac{\ds f_s(x,y,{\up{n}{\Omega}}) -
           f_r(x,y,{\up{n}{\Omega}})
      }
      {\ds {\up{n+1}{r}} - {\up{n+1}{s}}}
 +
 \frac{
        \left(
            \frac{\ds d}{\ds dx}
            \left(
                {\ds {\up{n+1}{r}} - {\up{n+1}{s}} }
            \right)
        \right)_r
      }
      {\ds {\up{n+1}{r}} - {\up{n+1}{s}}}
 =
  \smallskip \\
 &
 =
 \
 \frac{\ds f_s(x,y,{\up{n}{\Omega}}) -
           f_r(x,y,{\up{n}{\Omega}})
      }
      {\ds {\up{n+1}{r}} - {\up{n+1}{s}}}
 +
   \left(
     \frac{\ds d}{\ds dx}
     \ln
     \left(
        {\ds {\up{n+1}{r}} - {\up{n+1}{s}} }
     \right)
   \right)_r
   \, .
 \end{array}
 \eqno (5.8)
$$
Используем равенство (5.7):
$$
 \left.
    \left(
 \!\!
        \left(
            \pt_x
            +
            {\up{n+1}{s}}
            \pt_y
        \right)
        \ln
        \frac
        {
            \pt_{y_0}
            {\up{n+1}{g_r}} (x,y_0)
        }
        {
            {\ds {\up{n+1}{r}} - {\up{n+1}{s}} }
        }
    \right)
 \right|
        _{y = {\up{n+1}{g_r}} (x,y_0)}
 \!\!\!\!\!\!\!\!
 \!\!\!\!\!\!\!\!
 \!\!\!
 =
 \
    \left.
        \frac{\ds f_s(x,y,{\up{n}{\Omega}}) -
                f_r(x,y,{\up{n}{\Omega}})
            }
            {\ds {\up{n+1}{r}} - {\up{n+1}{s}}}
    \:
    \right|
           _{y = {\up{n+1}{g_r}} (x,y_0)}
   \, .
$$
Интегрируем вдоль характеристики согласно формуле (2.10), учитывая
начальное условие $(5.7^0)$:
$$
            \pt_{y_0}
            {\up{n+1}{g_r}} (x,y_0)
 =
        \frac
        {
            {\ds {\up{n+1}{r}} - {\up{n+1}{s}} }
        }
        {\ds r^0 - s^0}
    \exp
    \left\{
        {\ds \int_0^x}
        \frac
             {\ds {\up{n}{f_s}} - {\up{n}{f_r}}}
             {\ds {\up{n+1}{r}} - {\up{n+1}{s}}}
        \: d\tau
    \right\}
   \, .
   \eqno(5.9)
$$
Так как функции $f_{\omega},\ \omega = r,s,p,q,z$, непрерывно
дифференцируемы по $r,s,p,q,z,x,y$, а семейство
 $\{ {\up{n}{\Omega}} \}$
равномерно ограничено, то существует константа $a$ такая, что для
 ${\omega = r,s,p,q,z,}$
 ${\mu = r,s,p,q,z,x,y}$
$$
|f_{\omega}|\leq a,
 \quad
 \left|
   \frac{\ds \pt f_\omega}{\ds \pt \mu}
 \right|
  \leq a.
\eqno (5.10)
$$
 Используя оценки (4.5), (4.9), (5.10), получаем из (5.9) оценку
$$
\frac{1}{\psi (x)} \leq
   \pt_{y_0} {\ds\up{n+1}{g_r}} (x,y_)
   \leq \psi (x) \, ,
   \eqno (5.11)
$$
где
$$
\psi (x) =
 (2/\varepsilon)
 \exp
 \left\{
     2ax/\varepsilon
 \right\}
 \, .
$$

Оценки для $\pt_{y_0}{\up{n+1}{g_s}}(x,y_0)$ получаются
аналогично.
 Напомним, что всего характеристик две, т.е.
 для ${\omega = r,s,p,q,z}$
 ${{\ds\up{n+1}{g_\omega}} (x,y_0)}$
 совпадает либо с
 ${{\ds\up{n+1}{g_r}} (x,y_0)}$,
 либо с
 ${{\ds\up{n+1}{g_s}} (x,y_0)}$.
 Тем самым для
 $\omega = r,s,p,q,z$, $n=1,2,3,\dots$ справедливы оценки
$$
\frac{1}{\psi (x)} \leq
   \pt_{y_0} {\ds\up{n+1}{g_\omega}} (x,y_0)
   \leq \psi (x) \, .
   \eqno (5.12)
$$
Оценка для ${n = 0}$ нам не потребуется.

Дифференцируя формулы (5.5) по параметру $y_0$, получим
$$
\left\{
 \begin{array}{l}
   \pt_x
   (
     \pt_{y_0}
     {\up{n+1}{\bar \omega}}
   )
     =
     \sum \lms_{\mu = r,s,p,q,z}
        \frac{\ds\pt f_\omega}
          {\ds\pt {\up{n}{\mu}}}
        \frac{\ds\pt {\up{n}{\mu}}}
          {\ds\pt {\up{n+1}{g_\omega}}}
        \frac{\ds\pt {\up{n+1}{g_\omega}}}
          {\ds\pt y_0}
        +
        \frac{\ds\pt f_\omega}
          {\ds\pt {\up{n+1}{g_\omega}}}
        \frac{\ds\pt {\up{n+1}{g_\omega}}}
          {\ds\pt y_0}
   \, ,
   \medskip \\
       \pt_{y_0}
     {\up{n+1}{\bar \omega}}
     (0, y_0)
     =
     \omega_y^0(y_0)
     \, , \qquad \omega = r,s,p,q,z
     \, .
 \end{array}
 \right.
 \eqno (5.13)
$$
Вводя обозначение
$$
V_0 = \max\lms_{\omega = r,s,p,q,z}
 \sup\lms_{y_0 \in
                  \left.
                      {G(\bar x, \bar y)}
                  \right|_{x = 0}
          }
 |\omega_y^0 (y_0)|
 \, ,
$$
рассмотрим ассоциированную с задачей (5.13) мажорантную задачу
$$
        \frac{\ds d}{\ds dx}
        V
        =
        5a\psi^2 (x)\; V + a\psi (x) \, ,
  \qquad
        V(0)=V_0 \, .
 \eqno (5.14)
$$
Задача Коши (5.14) в силу линейности имеет решение на всем
компакте ${{G(\bar x, \bar y)}}$.

Начальное приближение $\up{0}{\omega}$ удовлетворяет неравенству
$$
|\pt_y {\up{0}{\omega}}|
 \leq
 \psi (x)\; V(x) \, , \quad
 \omega = r,s,p,q,z,
$$
так как ${\psi \geq 1}$ и правые части задачи (5.14)
неотрицательны. Предположим, что
$$
|\pt_y {\up{n}{\omega}}|
 \leq
 \psi (x)\; V(x) \, , \quad
 \omega = r,s,p,q,z.
 \eqno (5.15)
$$
Учитывая оценки (5.10), (5.12), (5.15), получим
$$
 \begin{array}{l}
 \left|
 \pt_{y_0} {\up{n+1}{\bar\omega}}
 \right|
 \up{(5.13),\ (2.10)}{\le}
 |\omega^0_y| +
 \medskip \\
 + {\ds\int_0^x}
 \left(
   \sum\lms_{\mu = r,s,p,q,z}
     \left|
        \frac{\ds\pt f_\omega}
          {\ds\pt {\up{n}{\mu}}}
     \right|
     \left|
        \frac{\ds\pt {\up{n}{\mu}}}
          {\ds\pt {\up{n+1}{g_\omega}}}
     \right|
     \left|
        \frac{\ds\pt {\up{n+1}{g_\omega}}}
          {\ds\pt y_0}
     \right|
     +
     \left|
        \frac{\ds\pt f_\omega}
          {\ds\pt {\up{n+1}{g_\omega}}}
     \right|
     \left|
        \frac{\ds\pt {\up{n+1}{g_\omega}}}
          {\ds\pt y_0}
     \right|
 \right)
 \; d\tau
 \le
 \medskip \\
 \le
 V_0
 +
 {\ds\int_0^x}
 (5a\psi(\tau)\psi(\tau)\; V(\tau) + a\psi(\tau))
 \; d\tau
 \up{(5.14)}{=}
 V(x)
 \, .
 \end{array}
\eqno(5.16)
$$
Из формулы (5.6), пользуясь оценками (5.12), (5.16), получаем
$$
 \left|
   \pt_y
   {\up{n+1}{\omega}} (x,y)
 \right|
 \leq
 \psi (x) \; V(x) \, , \quad \omega = r,s,p,q,z.
$$
Итак,
$$
\Phi (x) = \psi (x) \; V(x) \, . \ \square
$$



\section{Существование и единственность решения}

 Единственность решения задачи Коши (3.32), $(3.32^0)$ следует из
 теоремы единственности (\cite{РЯ}, гл.~1, {\S}~8, п.~2).

   В этом разделе, как и в разделе 2,
будем рассматривать задачу Коши для системы общего вида (2.8),
$(2.8^0)$ относительно неизвестной вектор-функции $u=(u_1,\dots
,u_m)$. Будем предполагать, что последовательность приближенных
решений $\{{\up{n}{u}}(x,y)\}$ и последовательность их производных
$\{{(\up{n}{u})}_y(x,y)\}$ равномерно ограничены на компакте (4.6)
для произвольной точки
 ${(\bar x, \bar y)\in[0,+\infty)\times \mathbb{R}}$.
Все рассмотрения будем проводить на этом компакте.
В~доказательствах будем следовать стандартной схеме, приведенной в
\cite{РЯ}. Также использована \cite{Туницкий-канд-дис}.


   Пусть $\varphi (u)\in C^1$. Рассмотрим функцию
$\tilde \varphi (\lambda )=\varphi (\bar u+\lambda (u-\bar u)),
\lambda \in \mathbb{R}$. По формуле Ньютона-Лейбница имеем
\cite{Туницкий-канд-дис}
$$
  \tilde \varphi (1)-\tilde \varphi (0)
  \! = \!
  \varphi (u)-\varphi (\bar u)
  \! =
  \!\!
  {\ds \int_0^1}
  \!\!\!
  \tilde \varphi _{\lambda}d\lambda =
  \!\!
  {\ds \int_0^1
  \!\!
  \left(
     \sum\lms_{j=1}^{m}\frac{\pt \varphi}{\pt u_j}
     (\bar u+\lambda (u-\bar u))(u_j-\bar u_j)
     \!\!
  \right)
  \!\!
   d\lambda.
}
$$
Доказанная формула
$$
  \varphi (u)-\varphi (\bar u)=
{\ds
   \sum\lms_{j=1}^m(u_j-\bar u_j)
   \int_0^1 \frac{\pt \varphi}{\pt u_j} (\bar u+\lambda (u-\bar u))
   \; d\lambda
}
 \eqno(6.1)
$$
называется {\it формулой конечных приращений.}


\subsection{Непрерывность решения}

{\bf Лемма~6.1.}\ Вектор-функция
$\lim\lms_{n\to\infty}{\up{n}{u}}(x,y)$ непрерывна.
\\
 {\bf Доказательство.}\ Классическая теорема анализа дает достаточное
условие сходимости функциональной последовательности к непрерывной
функции: элементы последовательности должны быть непрерывными
функциями, и последовательность должна равномерно сходиться.

   Равномерная сходимость последовательности $\{\up{n}{u}\}$ будет
следовать из равномерной сходимости ряда
${\ds\sum_{n=0}^{\infty}}(\up{n+1}{u}-\up{n}{u})(x,y)$.

 Пусть имеются два последовательных приближения
 ${\up{n+1}{u}}$,
 ${\up{n}{u}}$,
 т.е., согласно (4.1),
$$
\begin{array}{l}
    ( \pt_x
           +
           \xi_i(x,y,{\up{n+1}{u}})
           \pt_y
    )
    {\up{n+1}{u_i}}
    =
    f_i(x,y,{\up{n}{u}})
     \; ,
    \smallskip \\
    ( \pt_x
           +
           \xi_i(x,y,{\up{n+1}{u}})
           \pt_y
    )
    {\up{n}{u}}_i
    +
    \xi_i(x,y,{\up{n}{u}})
    {\up{n}{u}}_{iy}
    =
    f_i(x,y,{\up{n-1}{u}})
    +
    \xi_i(x,y,{\up{n+1}{u}})
    {\up{n}{u}}_{iy}
      \; ,
    \smallskip \\
    {\up{n+1}{u}}
    (0,y)
    =
    {\up{n}{u}}
    (0,y)
    =
    u^0(y)
    \, .
\end{array}
$$
 Вычитанием получаем задачу Коши для
 ${{\up{n+1}{u}} - {\up{n}{u}}}$:
$$
\left\{
\begin{array}{ll}
    (
           \pt_x
           +
           \xi_i(x,y,{\up{n+1}{u}})
           \pt_y
    )
    (
      {{\up{n+1}{u_i}} - {\up{n}{u}}_i}
    )
    \;
    =
    &
    f_i(x,y, {\up{n}{u}})
    -
    f_i(x,y, {\up{n-1}{u}})
    -
       \smallskip \\
     &
    -
    {\up{n}{u}}_{iy}
    (
        \xi_i(x,y,{\up{n+1}{u}})
        -
        \xi_i(x,y,{\up{n}{u}})
    )
    \; ,
    \smallskip \\
    ({\up{n+1}{u}} - {\up{n}{u}})
    (0,y)
    =0
    \, .
\end{array}
\right.
$$
   Для ${\up{n+1}{r}} = {\up{n+1}{u}} - {\up{n}{u}}$
 с помощью формулы конечных приращений (6.1) получим систему
$$
\left\{
\begin{array}{l}
    ( \pt_x
           +
           \xi_i(x,y,\up{n+1}{u})
           \pt_y
    )
    \up{n+1}{r_i}
    =
    \smallskip
    \\
{\ds
   =
   \sum\lms_{j=1}^m
   {\up{n}{r}}_j
   \int_0^1
\left(
   \frac{\pt f_i}{\pt u_j}
   (
     x,y,
     \up{n-1}{u}
     +
     \lambda
     \up{n}{r}
   ) -
   {\up{n}{u}}_{iy}
   \frac{\pt \xi_i}{\pt u_j}
   (
       x,y,
       {\up{n}{u}}
       +
       \lambda
       {\up{n+1}{r}}
   )
\right)
   d\lambda
}
     \; ,
     \smallskip \\
     {\up{n+1}{r}}
     (0,y)
     =
     0
     \, .
\end{array}
\right.
$$
Воспользовавшись формулой (2.10), получим
$$
    |\up{n+1}{r_i}|\le
{\ds
   \int_0^x
   \max\lms_k|{\up{n}{r}}_k|
   \sum\lms_{j=1}^m
   \int_0^1
\left(
   \left|\frac{\pt f_i}{\pt u_j}\right| +
   |{\up{n}{u}}_{iy}|
   \left|\frac{\pt \xi_i}{\pt u_j}\right|
\right)
   d\lambda \; d\tau.
}
 \eqno(6.2)
$$
Из формул (3.32) видно, что из {$C^1\mbox{-ограниченности}$}
вектор-функций $\rho,\sigma,\pi,\kappa$ и равномерной
ограниченности последовательности $\{\up{n}{u}\}$ следует
равномерная ограниченность функций ${\ds \left|\frac{\pt f_i}{\pt
u_j}\right|}$. Величина ${\ds \left|\frac{\pt \xi_i}{\pt
u_j}\right|}$ есть либо~0, либо~1, последовательность
$\{{\up{n}{u}}_{yi}\}$ равномерно ограничена, и, значит,
$$
{\ds
   \sum\lms_{j=1}^m
   \int_0^1
\left(
   \left|\frac{\pt f_i}{\pt u_j}\right| +
   |{\up{n}{u}}_{iy}|
   \left|\frac{\pt \xi_i}{\pt u_j}\right|
\right)
   d\lambda
   \le C\; ,
}
$$
где C --- некоторая константа, т.е. левая часть неравенства
равномерно ограничена. Введем обозначение
$$
  R_n(x)=\max\lms_i\sup\lms_{(\tau ,y)\in [0,x]\times \mathbb{R}}
  |{\up{n}{r}}_i(\tau ,y)|.
$$
Тогда формула (6.2) примет вид
$$
  R_{n+1}(x)\le C\int_0^xR_n(\tau )\; d\tau .
   \eqno(6.3)
$$
Последовательно применяя формулу (6.3), получаем:
$$
  R_{n+1}(x)
  \le
  \max\lms_{\tau\in [0,x]}R_1(\tau )\frac{(Cx)^n}{n!}
  \le
  \max\lms_{\tau\in [0,\bar x]}R_1(\tau )\frac{(C \bar x)^n}{n!}\; ,
$$
откуда согласно признаку Вейерштрасса следует равномерная
сходимость ряда
${\ds\sum_{n=0}^{\infty}}({\up{n+1}{u}}-{\up{n}{u}})(x,y)$
 в полосе
 ${[0,\bar x]\times \mathbb{R}}$. Поэтому предельная
 вектор-функция
 $\lim\lms_{n\to\infty}{\up{n}{u}}(x,y)$
 непрерывна
в полосе ${[0,\bar x]\times \mathbb{R}.}\ \Box$


\subsection{Непрерывная дифференцируемость решения}

{\bf Лемма~6.2.}\ Вектор-функция
$\lim\lms_{n\to\infty}{\up{n}{u}}(x,y)$ непрерывно
дифференцируема.
\\
{\bf Доказательство.}\ Для доказательства непрерывности
вектор-функции $\lim\lms_{n\to\infty}{{(\up{n}{u})}_y}$
воспользуемся теоремой Арцела. Равномерная ограниченность
семейства $\{{\up{n}{u}}_y\}$ уже доказана. Докажем равностепенную
непрерывность семейства $\{{\up{n}{u}}_y\}$.

   Сначала докажем равностепенную непрерывность семейств
$\{\up{n}{u}\},\{\up{n}{g}\}$. По формуле конечных приращений
(6.1) имеем:
$$
 \!
 \begin{array}{ll}
  {\up{n}{u}}(x_1,y_1)-{\up{n}{u}}(x_2,y_2)&
  \!\!\!\!
  =(x_1-x_2){\ds\int_0^1}
  {\up{n}{u}}_x(x_2+\lambda (x_1-x_2),y_2+\lambda (y_1-y_2))
  \; d\lambda +
 \medskip   \\
  &
  \!\!
  +(y_1-y_2){\ds\int_0^1}
  {\up{n}{u}}_y(x_2+\lambda (x_1-x_2),y_2+\lambda (y_1-y_2))
  \; d\lambda.
\end{array}
$$
 Равномерная ограниченность семейства $\{{\up{n}{u}}_x\}$ следует
 из формул
 (4.1) и равномерной ограниченности семейств
 $\{\up{n}{u}\},\, \{{\up{n}{u}}_y\}$. Отсюда следует
 равностепенная непрерывность последовательности $\{\up{n}{u}\}$.

   Для вектор-функции ${\up{n}{g}}(\tau ,x,y)$ по формуле конечных приращений
(6.1) имеем:
$$
\begin{array}{l}
  {\up{n}{g}}(\tau ,x_1,y_1)-{\up{n}{g}}(\tau ,x_2,y_2)=
                                \medskip   \\
  =(x_1-x_2){\ds\int_0^1}{\up{n}{g}}_x(
  \tau ,x_2+\lambda (x_1-x_2),y_2+\lambda (y_1-y_2))
  \; d\lambda +
                                \medskip   \\
  +\:(y_1-y_2){\ds\int_0^1}{\up{n}{g}}_y(
  \tau ,x_2+\lambda (x_1-x_2),y_2+\lambda (y_1-y_2))
  \; d\lambda.
\end{array}
$$
Равномерная ограниченность семейств
$\{{\up{n}{g}}_x\},\{{\up{n}{g}}_y\}$ следует из формул (2.11) и
из равномерной ограниченности семейств
$\{\up{n}{u}\},\{{\up{n}{u}}_y\}$ (в нашем случае ${\xi_{yi}=0,}$
 ${\pt\xi_i/\pt{\up{n}{u}}_j}$ либо 0, либо 1). Поэтому
$$
  |{\up{n}{g}}_i(\tau ,x_1,y_1)-{\up{n}{g}}_i(\tau ,x_2,y_2)|
  \le \mbox{\rm const}\,(|x_1-x_2|+|y_1-y_2|),
$$
т.е. семейство $\{\up{n}{g}\}$ равностепенно непрерывно.

   Рассмотрим теперь на компакте $G(\bar x, \bar y)$ функцию $\up{n+1}{u_{iy}}$ :
$$
\begin{array}{l}
  {\up{n+1}{u_{iy}}}(x,y)
  =
  u^0_{iy}
  ({\up{n+1}{g_i}}(0,x,y))
  +
 \medskip   \\
  \phantom{123}
  +
  {\ds\int_0^x}
  \left\{
    -
    {\up{n+1}{u_{iy}}}
    \left(
        \xi_{iy}
        +
        \sum\lms_{j}
            \frac{\ds \pt \xi_{i}}
                 {\ds \pt {\up{n+1}{u_j}}}
            \cdot
            {\up{n+1}{u_{jy}}}
    \right)
    +
    \left(
        f_{iy}
        +
        \sum\lms_{j}
            \frac{\ds \pt f_{i}}
                 {\ds \pt {\up{n}{u}_j}}
            \cdot
            {\up{n}{u}_{jy}}
    \right)
  \right\}
 \medskip
 \\
  \phantom{123}
  \left(
     \tau ,
     {\up{n+1}{g_i}}
     (\tau ,x,y),
     {\up{n}{u}}
     (\tau ,{\up{n+1}{g_i}}(\tau ,x,y))
   \right)
   d\tau
   \, .
\end{array}
 \eqno(6.4)
$$
В нашем случае ${\xi_{iy} = 0}$,
 ${
            \pt \xi_i/
                 \pt {\up{n+1}{u_j}}
 }$
 либо $0$, либо $1$.
   Поскольку равностепенно непрерывны семейства
$\{\up{n}{u}\},\{\up{n}{g}\}$, то равностепенно непрерывны и
семейства
$$
\begin{array}{l}
\{u^0_{iy}({\up{n+1}{g_i}}(0,x,y))\}
                          \: ,     \medskip   \\
\left\{
    f_{iy}
  \left(
     \tau ,{\up{n+1}{g_i}}(\tau ,x,y),
     {\up{n}{u}}(\tau ,{\up{n+1}{g_i}}(\tau ,x,y))
   \right)
\right\}
                          \: ,     \medskip   \\
\left\{
            \frac{\ds \pt f_{i}}
                 {\ds \pt {\up{n}{u}_j}}
  \left(
     \tau ,{\up{n+1}{g_i}}(\tau ,x,y),
     {\up{n}{u}}(\tau ,{\up{n+1}{g_i}}(\tau ,x,y))
   \right)
\right\}
 \: .
\end{array}
$$
Воспользуемся теоремой Кантора: \textit{непрерывная на компакте в
$\mathbb{R}^s$ функция равномерно непрерывна на этом компакте}.
   Учитывая равномерную ограниченность семейства
$\{{\up{n}{u}}_y\}$ и теорему Кантора, получаем равностепенную
непрерывность правых частей равенства (6.4), т.е. семейства
$\{{\up{n}{u}}_y\}$, на компакте $G(\bar x,\bar y)$.

   Воспользуемся теоремой Арцела: \textit{если функциональная последовательность
равномерно ограничена и равностепенно непрерывна на компакте, то
из нее можно выделить подпоследовательность, равномерно сходящуюся
на этом компакте.} Поэтому существует равномерно сходящаяся
подпоследовательность $\{{\up{n_k}{u_y}}\}$. По теореме о
почленном дифференцировании функциональной последовательности
вектор-функция $\lim\lms_{k\to\infty}{\up{n_k}{u}}(x,y)$
непрерывно дифференцируема по $y$, и
$$
  \pt_y \lim\lms_{k\to\infty}{\up{n_k}{u}}=
  \lim\lms_{k\to\infty}{\up{n_k}{u_y}}.                     \eqno(6.5)
$$
Из доказанной ранее равномерной сходимости последовательности
$\{\up{n}{u}\}$ следует, что
$$
  \lim\lms_{k\to\infty}{\up{n_k}{u}}=
  \lim\lms_{n\to\infty}{\up{n}{u}}.                          \eqno(6.6)
$$
Воспользовавшись равенством (6.6), получаем из равенства (6.5)
правило дифференцирования по $y$ для вектор-функции
$\lim\lms_{n\to\infty}{\up{n}{u}}$:
$$
  \pt_y \lim\lms_{n\to\infty}{\up{n}{u}}=
  \lim\lms_{k\to\infty}{\up{n_k}{u_y}}.                      \eqno(6.7)
$$
Вектор-функция $\pt_y \lim\lms_{n\to\infty}{\up{n}{u}}$
непрерывна, так как непрерывны вектор-функции $\up{n_k}{u_y}$.

   Правило дифференцирования по $x$ для вектор-функции
$\lim\lms_{n\to\infty}{\up{n}{u}}(x,y)$ получаем, переходя в
равенстве (4.1) к пределу:
$$
  \pt_x \lim\lms_{n\to\infty}{\up{n}{u}}_i(x,y)=
  f_i(x,y,\lim\lms_{n\to\infty}{\up{n}{u}})-
  \xi_i(x,y,\lim\lms_{n\to\infty}{\up{n}{u}})
  (\pt_y \lim\lms_{n\to\infty}{\up{n}{u}}_i).      \eqno(6.8)
$$
Непрерывность производной по $x$ обеспечивается непрерывностью
вектор-функций $f,\; \xi ,\;\lim\lms_{n\to\infty}{\up{n}{u}},\;
\pt_y \lim\lms_{n\to\infty}{\up{n}{u}}$.

   Итак, доказано, что
 ${\lim\lms_{n\to\infty}{\up{n}{u}}(x,y)
 ~\in~C^1(G(\bar x,\bar y))}$.
Подставив построенную функцию в задачу Коши (2.8), $(2.8^0)$ и
воспользовавшись равенством (6.8), видим, что она является
решением задачи Коши.

   Остается сказать, что для любой точки
${(x,y)\in [0,+\infty)\times \mathbb{R}}$ можно построить компакт
вида ${G(\bar x,\bar y)}$, содержащий эту точку. Поскольку в
пересечении двух таких компактов начальные приближения совпадают и
равны начальному приближению, заданному в полуплоскости
${[0,+\infty)\times \mathbb{R}}$, то совпадают и предельные
функции, являющиеся, следовательно, ограничениями предельной
функции, определенной в полуплоскости
 ${[0,+\infty)\times \mathbb{R}}$,
 на соответствующую область.

 Отметим, что теорема Арцела, являющаяся чистой теоремой существования,
 использутся в доказательстве, но в
 построении решения не используется. $\Box$



\section{Основной результат}
\subsection{Теорема существования и единственности}

Выше было установлено, что при выполнении неравенств
 (4.4), (4.7), (4.8)
 в полуплоскости ${x \geq 0}$ существует единственное
${C^1\mbox{-решение}}$ задачи (3.32), $(3.32^0)$. Сформулируем
условия на коэффициенты уравнения (3.1) и начальные данные,
обеспечивающие выполнение неравенств
 (4.4), (4.7), (4.8).

Пусть $M_1, M_2$ --- произвольные положительные постоянные,
$\varepsilon , \delta$ --- постоянные, связанные условиями
$$
0 < \varepsilon \leq \delta ,
 \quad (\delta - \varepsilon)/2 < 1 \, ,
 \eqno (7.1)
$$
$\eta (x) \in C^0(\mathbb{R})$ --- произвольная неотрицательная
функция.

 Введем обозначения
$$
\begin{array}{l}
 N_1 = \max \{
               M_1, \frac{1}{2} M_2 (4M_1 + 9 M_1^2)
            \},
 \smallskip
 \\
    N_2 = M_1 M_2 \, ,
 \medskip
 \\
    \tilde{\eta}(x) = \frac{\ds 1}{\ds 1 + 2M_1 |x|} \; \eta (x) \, .
\end{array}
\eqno (7.2)
$$
Будем предполагать выполненными следующие условия:
$$
 \begin{array}{lll}
   1) & \mbox{коэффициенты}\ A,B,C,D \in C^2(\mathbb{R}^5)\, ;
      & \phantom{1234567890}
   \smallskip
   \\
   2)
      & \mbox{начальные функции} \ z^0 \in C^3(\mathbb{R}),
        \ p^0 \in C^2(\mathbb{R});
      &
   \medskip
   \\
   3) & |a|\leq M_1, \ \ 1/\Delta \leq M_2 \, ,
      &
   \medskip
   \\
      & \left|
          \frac{\ds \pt a}{\ds \pt \omega}
        \right|
        \leq M_1 \eta (x) \, ,
        \qquad
        \left|
          \frac{\ds \pt a}{\ds \pt z}
        \right|
        \leq M_1 \tilde{\eta} (x) \, ,
      &
   \smallskip
   \\
      & \mbox{где} \ \ a=B,C,D,\Delta, \quad \omega = x,y,p,q;
      &
   \smallskip
   \\
   4) & (6 N_1 + 8 N_2)
        {\ds\int_{-\infty}^{+\infty}}
          \eta (x)
        \; dx
        \leq
        (\delta - \varepsilon)/2
        \, ;
      &
   \medskip
   \\
   5) & |r^0|, |s^0|
        \le
        1 -
        (\delta - \varepsilon)/2
        \, , \quad
        |z_y^0|, |p^0| \le 1
        \, ;
      &
   \medskip
   \\
   6) & \inf\lms_{y \in \mathbb{R}} r^0 (y) -
        \sup\lms_{y \in \mathbb{R}} s^0 (y) \geq \delta > 0 \, .
\end{array}
\eqno (7.3)
$$

{\bf Теорема~7.1.} При выполнении условий (7.3)
$C^3\mbox{-гладкое}$ решение задачи Коши (3.1), (3.3) существует и
единственно.

 {\bf Доказательство.}
 Докажем, что условия (7.3) обеспечивают выполнение условий
 (4.4), (4.7), (4.8)
 и требований, которые предшествуют им в начале раздела 4.2.

 Условия (4.4) сформулированы в предположении
 ${C^1\mbox{-гладкости}}$
 и ограниченности
 начальных данных $r^0, s^0$.
 Сначала покажем, что условия (7.3)
 обеспечивают условие ${r^0, s^0 \in C^1}$. Из первых двух
 уравнений системы $(3.32^0)$ получаем
$$
 {\ds
    r^0 (y) - s^0 (y)
    =
    \frac{\Delta
        (0,y,z^0(y),p^0(y),z_y^0(y))
    }
    {z_{yy}^0 (y) + B
        (0,y,z^0(y),p^0(y),z_y^0(y))
    }
 }
 \, ,
$$
откуда
$$
 z_{yy}^0 (y) + B
 (0,y,z^0(y),p^0(y),z_y^0(y))
  =
 \frac{\Delta
     (0,y,z^0(y),p^0(y),z_y^0(y))
 }
 {r^0 (y) - s^0 (y)}
 \, .
$$
Таким образом, неравенство
$$
 z_{yy}^0 (y) + B
 (0,y,z^0(y),p^0(y),z_y^0(y))
 \neq
 0
$$
обеспечивается условием (7.3), пункт 6), и условием
  ${1/\Delta \leq M_2}$ из условия (7.3), пункт 3).

 Тогда условие ${r^0, s^0 \in C^1}$ обеспечивают
 формулы $(3.32^0)$
 и условия
 (7.3), пункты 1), 2).
 Ограниченность
 функций $r^0, s^0$ обеспечивается
 условием (7.3), пункт 5).

 Условия (4.4) сформулированы также в предположении
 ${C^1\mbox{-гладкости}}$
 и ограниченности
 начальных данных $p^0, q^0$.
 Согласно $(3.32^0)$\ ${q^0 = z_y^0}$.
 $C^1\mbox{-гладкость}$ функций $p^0, q^0$
 обеспечивается
 условием (7.3), пункт 2).
 Ограниченность
 функций $p^0, q^0$
 обеспечивается
 условием (7.3), пункт 5).

 Далее, условия (4.4) сформулированы в предположении
 ${C^1\mbox{-гладкости}}$
 и ограниченности
 вектор-функций $\rho, \sigma, \pi, \kappa$.
 Согласно (3.32) в выражения для этих вектор-функций входят первые
 производные функций
 ${B, C, D, \Delta}$
 и выражение ${1/\Delta}$.
 Согласно условиям (7.3), пункт 3) функция
 $1/\Delta$ отделена от нуля положительной константой.
 Тогда ${C^1\mbox{-гладкость}}$
 вектор-функций $\rho, \sigma, \pi, \kappa$
 обеспечивается $C^2\mbox{-гладкостью}$ функций
 ${A, B, C, D}$
 согласно условиям (7.3), пункт 1).
 Ограниченность
 вектор-функций $\rho, \sigma, \pi, \kappa$
 обеспечивается
 условиями (7.3), пункт 3).

 Используя условие (7.3), пункт 3), получаем оценки
$$
 \begin{array}{l}
   |\rho_i|, |\sigma_i|\leq M_1 \eta \, , \quad i=0,13;
 \smallskip
 \\
   |\rho_i|, |\sigma_i|\leq \frac{1}{2} M_2 (4M_1 + 9 M_1^2) \eta \, , \quad i=1,2,7,8;
 \smallskip
 \\
   |\rho_i|, |\sigma_i|\leq M_1 M_2 \tilde{\eta} \, , \quad i=3,4,5,6,9,10,11,12;
 \smallskip
 \\
   |\pi_i|, |\kappa_i|\leq M_1 \, , \quad i=0,1;
 \end{array}
$$
Отсюда согласно формулам (4.3), (7.2) имеем оценки
$$
 \begin{array}{l}
   \alpha_1 (x) \leq \max
     \{
       M_1, \frac{1}{2} M_2 (4M_1 + 9 M_1^2)
     \}
     \eta (x) = N_1 \eta (x) \, ,
 \smallskip
 \\
   \alpha_2 (x) \leq M_1 M_2 \tilde{\eta} (x) =
      N_2 \tilde{\eta} (x) \, ,
 \smallskip
 \\
   \alpha_3 \leq M_1 \, .
 \end{array}
 \eqno (7.4)
$$
Тогда пункты 4), 5) условия (7.3) обеспечивают выполнение второго
неравенства (4.4):
$$
 \begin{array}{l}
   U_0 +
        6
        {\ds\int_{-\infty}^{+\infty}}
            \alpha_1 (x)
        \; dx
        +
        8
        {\ds\int_{-\infty}^{+\infty}}
        (1 + 2\alpha_3 |x|)
        \; \alpha_2 (x)
        \; dx
        \le
 \medskip
 \\
   \le
        U_0 +
        (6 N_1 + 8 N_2)
        {\ds\int_{-\infty}^{+\infty}}
            \eta (x)
        \; dx
        \le
        (
            1
            -
            (\delta - \varepsilon)/2
        )
        +
            (\delta - \varepsilon)/2
        =
        1
        \, .
 \end{array}
$$

Первое неравенство (4.4) следует из условия (7.3), пункт 5).
Неравенство (4.7) есть условие (7.3), пункт 6). Неравенство (4.8)
с учетом обозначений (7.2) и оценок (7.4) есть условие (7.3),
пункт 4).

Таким образом, условия (7.3) обеспечивают выполнение неравенств
 (4.4), (4.7), (4.8),
  которые в свою очередь обеспечивают существование единственного
$C^1\mbox{-решения}$ задачи (3.32), $(3.32^0)$ в полуплоскости ${x
\ge 0}$.

Решение в полуплоскости ${x \leq 0}$ строится аналогично.
$C^1\mbox{-гладкость}$ решения
 на всей плоскости
 обеспечивается одинаковыми
начальными данными.

$C^1\mbox{-решению}$ задачи (3.32), $(3.32^0)$ соответствует
{$C^3\mbox{-решение}$} задачи (3.1), (3.3), поскольку
 согласно теореме~3.2 для этого достаточно
 выполнения неравенства
    $
    {\Delta > 0},
    $
 справедливого в силу условия ${1/\Delta \leq M_2}$
 (условие (7.3), пункт 3)).
 $\square$

 {\bf Пример~7.1.} Покажем, что условия (7.3) задают непустое
множество входных данных. Пусть $
 {A=1/16,} \quad
 {B=1/2,} \quad
 {C=0,} \quad {D=0,} \quad
 {z^0 = 1,} \quad {p^0 = 1.}
$
 Тогда согласно (3.2) ${\Delta^2 = 4A,}$ откуда
 ${\Delta = 1/2.}$
 Выражения $(3.32^0)$ для начальных данных принимают вид
$
 {r^0, s^0 = (C \pm \Delta)/2B,}
$
 откуда
$
 {r^0, s^0 = \pm 1/8.}
$
 Таким образом, существуют постоянные
$
 {M_1 = 1/2,} \quad
 {M_2 = 2}\, , \quad
 {\delta = 1/4,} \quad
 {\varepsilon = 1/4}
$
 и функции
 ${\eta(x) = 0}, \quad {\tilde \eta(x) = 0}$,
 для которых выполнены соотношения (7.3).
 Постоянные
 $N_1,\ N_2$ определяются выражениями (7.2). $\Box$

{\bf Пример~7.2.} Уравнение с постоянными коэффициентами. Система
(3.32) сводится к независимой подсистеме (3.33), которая принимает
вид
$$
 \begin{array}{l}
 ( \partial_x
   +
   s \, \partial_y
 )
 \: r = 0 \, ,
 \medskip
 \\
 ( \partial_x +
    r \, \partial_y
 )
 \: s = 0 \, .
 \end{array}
$$
Для этой системы от условий (7.3) остаются только условия на
начальные функции $r^0,\ s^0$.
 После нахождения функций $r,\ s$ решаются
 линейные уравнения (3.32) для функций $p, \ q, \ z$.
 Начальные данные $p^0,\ q^0,\ z^0$ в этом случае выражаются через
 $r^0,\ s^0$ при помощи $(3.32^{00})$.
 $\Box$


\subsection{Начальные данные}

Конкретизируем достаточные условия, обеспечивающие выполнение
условий (7.3) в пунктах, касающихся начальных данных.

Пусть
$$
 {m_1 \geq 0}, \quad {m_2 \geq 0}, \quad
 {L_1 > 0}, \quad {L_2 >0}, \quad {L_3 > 0}
$$
 --- произвольные постоянные. Введем постоянные
$$
\begin{array}{l}
 m_3 = \frac{\ds m_1}{\ds L_1 - m_1} \: (1 + L_3 + \frac{\ds m_2}{\ds 2L_1}) +
       \frac{\ds m_2}{\ds 2L_1} \, ,
 \medskip
 \\
 m_4 = \frac{\ds m_1 (2 + L_3)}{\ds L_1 - m_1} +
       \frac{\ds m_2}{\ds 2(L_1 - m_1)} \, ,
 \medskip
 \\
 L_4 = \frac{\ds 1}{\ds 2 M_2 L_2} \, ,
\end{array}
\eqno(7.5)
$$
где $M_2$ --- постоянная из условий (7.3), пункт 3).

Будем предполагать выполненными следующие условия:
$$
 \begin{array}{rlrl}
   1)
      & z^0 \in C^3(\mathbb{R}),
        \ p^0 \in C^2(\mathbb{R}) \, ;
        \phantom{1234567}
 &  6) & |C| \leq m_2 \, ;
   \smallskip
   \\
   2) & |z_y^0| \leq 1, \ \ |p^0| \leq 1 \, ;
 &   7) & \frac{\ds |\Delta|}{\ds 2L_1} \le L_3 \, ;
   \medskip
   \\
   3) & |z_{yy}^0| \leq m_1, \quad |p_y^0| \leq m_1 \, ;
 &   8) & L_3 + m_3 \le 1 - \frac{\ds \delta - \varepsilon}{\ds 2} \, ;
   \medskip
   \\
   4) & 0 < L_1 \le B \le L_2 \, ;
 &   9) & L_4 - m_4 \ge \delta/2 > 0 \, .
   \medskip
   \\
   5) & L_1 - m_1 > 0 \, ;
\end{array}
\eqno (7.6)
$$

{\bf Теорема~7.2.} Пусть коэффициенты $A, B, C, D$ уравнения (3.1)
удовлетворяют условиям (7.3), и пусть сверх того выполнены условия
(7.6). Тогда начальные функции $z^0,\ p^0$ и определяемые ими с
помощью формул $(3.32^0)$ функции $r^0,\ s^0$ удовлетворяют
условиям (7.3).

{\bf Доказательство.}
 Пункт 2) условий (7.3) следует из пункта 1) условий (7.6).
 Пункт 5) условий (7.3) для функций $z_y^0, \ p^0$
 есть пункт 2) условий (7.6).

 Покажем, что для функций $r^0, \ s^0$ пункт 5) условий (7.3)
 выполнен. Так как
$$
\begin{array}{l}
\ds
 \frac{C \pm \Delta - 2 p_y^0}{2(z_{yy}^0 + B)} -
 \frac{C \pm \Delta}{2 B}
 =
 \frac{ (C \pm \Delta)B - 2p_y^0 B - (C \pm \Delta)(z_{yy}^0 + B)}
 {2(z_{yy}^0 + B)B}
 =
 \bigskip
 \\
\ds
 =
 \frac{-2 p_y^0 B - (C \pm \Delta)z_{yy}^0}
 {2(z_{yy}^0 + B)B}
 =
 - \frac{p_y^0}{z_{yy}^0 + B} -
 \frac{(C \pm \Delta)z_{yy}^0}{2(z_{yy}^0 + B)B} \, ,
\end{array}
$$
то
$$
 \frac{C \pm \Delta - 2 p_y^0}{2(z_{yy}^0 + B)}
 =
 \frac{C \pm \Delta}{2 B} -
 \frac{p_y^0}{z_{yy}^0 + B} -
 \frac{(C \pm \Delta)z_{yy}^0}{2(z_{yy}^0 + B)B} \, .
\eqno(7.7)
$$
Вводя индекс ${\omega = r,s,}$ из равенств $(3.32^0)$, (7.7)
получаем оценки
$$
\begin{array}{l}
\ds
 |\omega^0|
 =
 \left|
 \frac{C \pm \Delta - 2 p_y^0}{2(z_{yy}^0 + B)}
 \right|
 \le
 \frac{|C| + |\Delta|}{2|B|} +
 \frac{|p_y^0|}{|z_{yy}^0 + B|} +
 \frac{(|C| + |\Delta|)|z_{yy}^0|}{2|z_{yy}^0 + B|\cdot |B|}
 \le
 \bigskip
 \\
\ds
 \up{(7.6),\ \mbox{{\scriptsize п.}}\, 3)-6)}{\le}
 \frac{|\Delta| + m_2}{2L_1} +
 \frac{1}{L_1 - m_1} \:
 \left(
 m_1 +
 \frac{m_1(|\Delta| + m_2)}{2L_1}
 \right)
 =
 \bigskip
 \\
\ds
 =
 \frac{|\Delta|}{2L_1} +
 \frac{m_2}{2L_1} +
 \frac{m_1}{L_1 - m_1} \:
 \left(
 1 +
 \frac{|\Delta|}{2L_1} +
 \frac{m_2}{2L_1}
 \right)
 \up{(7.6),\ \mbox{{\scriptsize п.}}\, 7)}{\le}
 \bigskip
 \\
\ds
 \le
 L_3 +
 \frac{m_2}{2L_1} +
 \frac{m_1}{L_1 - m_1} \:
 \left(
 1 +
 L_3 +
 \frac{m_2}{2L_1}
 \right)
 \up{(7.5)}{=}
\ds
 L_3 + m_3
 \up{(7.6),\ \mbox{{\scriptsize п.}}\, 8)}{\le}
 1 - \frac{\delta - \varepsilon}{2} \, .
\end{array}
$$
Покажем, что выполнен пункт 6) условий (7.3). Так как
$$
\begin{array}{l}
\ds
 \frac{C \pm \Delta - 2 p_y^0}{2(z_{yy}^0 + B)} -
 \frac{\pm \Delta}{2 B}
 =
 \frac{(C \pm \Delta)B - 2p_y^0 B - (\pm \Delta)(z_{yy}^0 + B)}
 {2(z_{yy}^0 + B)B}
 =
 \bigskip
 \\
\ds
 =
 \frac{CB - 2 p_y^0 B - (\pm \Delta)z_{yy}^0}
 {2(z_{yy}^0 + B)B}
 =
 \frac{C - 2 p_y^0}{2(z_{yy}^0 + B)} -
 \frac{(\pm \Delta)}{2B} \cdot
 \frac{z_{yy}^0}{z_{yy}^0 + B}
 \, ,
\end{array}
$$
то
$$
 \frac{C \pm \Delta - 2 p_y^0}{2(z_{yy}^0 + B)}
 =
 \frac{\pm \Delta}{2 B}
 -
 \frac{(\pm \Delta)}{2B} \cdot
 \frac{z_{yy}^0}{z_{yy}^0 + B}
 +
 \frac{C - 2 p_y^0}{2(z_{yy}^0 + B)}
 \, .
\eqno(7.8)
$$
Вводя индекс ${\omega = r,s}$, из равенств $(3.32^0)$, (7.8)
получаем оценки
$$
\begin{array}{l}
\ds
 |\omega^0|
 =
 \left|
 \frac{C \pm \Delta - 2 p_y^0}{2(z_{yy}^0 + B)}
 \right|
 \ge
 \bigskip
 \\
 \ds
 \ge
 \inf
 \left|
   \frac{\pm \Delta}{2B}
 \right|
 -
 \sup
 \left|
   \frac{\pm \Delta}{2B} \cdot
   \frac{z_{yy}^0}{z_{yy}^0 + B}
 \right|
 -
 \sup
 \left|
   \frac{C - 2 p_y^0}{2(z_{yy}^0 + B)}
 \right|
 \up{(7.6),\ \mbox{{\scriptsize п.}}\, 4)}{\ge}
 \bigskip
 \\
 \ds
 \up{(7.3),\ \mbox{{\scriptsize п.}}\, 3)}{\ge}
 \frac{1}{M_2} \cdot \frac{1}{2L_2}
 -
 \sup
 \left|
   \frac{\pm \Delta}{2B} \cdot
   \frac{z_{yy}^0}{z_{yy}^0 + B}
 \right|
 -
 \sup
 \left|
   \frac{C - 2 p_y^0}{2(z_{yy}^0 + B)}
 \right|
 \ge
 \bigskip
 \\
 \ds
 \up{(7.6),\ \mbox{{\scriptsize п.}}\, 3),5),7)}{\ge}
 \frac{1}{M_2} \cdot \frac{1}{2L_2}
 -
 L_3 \cdot \frac{m_1}{L_1 - m_1}
 -
 \frac{m_2 + 2m_1}{2(L_1 - m_1)}
 \up{(7.5)}{\ge}
 \bigskip
 \\
 \ds
 \ge
 L_4 - m_4
 \up{(7.6),\ \mbox{{\scriptsize п.}}\, 9)}{\ge}
 \delta/2 > 0 \, .
\end{array}
$$
Поскольку знак функции $\omega^0$ определяется ее главной частью
${\pm \Delta/2B,}$ то функции $r^0$ и $s^0$ имеют разные знаки,
причем из условия ${B > 0}$ ((7.6), пункт 4)) и условия ${\Delta >
0}$ ((7.3), пункт 3)) следует, что ${r^0 > 0}, \quad {s^0 < 0}$.
Тогда из доказанных соотношений
 ${|\omega^0| \ge \delta/2}$
получаем $
 {\inf\lms_{y \in \mathbb{R}} r^0 (y)
 -
 \sup\lms_{y \in \mathbb{R}} s^0 (y) \geq \delta > 0}
$.
 $\Box$

{\bf Пример~7.3.} Условиям (7.6) удовлетворяют коэффициенты и
начальные функции примера~7.1. В этом случае ${M_2 = 2}$,
 ${m_i = 0}$, ${L_i = 1/2}$, ${i = 1,2,3,4}$.
Задавая малые возмущения входных данных из примера~7.1, получаем
пример уравнения с переменными коэффициентами, зависящими от
$x,y,z,z_x,z_y$, и с непостоянными начальными функциями. $\Box$


\section{Приложение. Контактный подход}
\subsection{Контактное преобразование может переводить
 классическое решение в решение, особое в каждой точке.
 }



 Рассмотрим нижеследующий пример.
 Автор благодарен В.~В.~Лычагину
 \cite{Лычагин-75}, \cite{Лычагин-79} и Л.~В.~Зильберглейту,
 введшим его в курс дела.

 {\bf Пример~8.1.} Рассмотрим пространство
 джетов первого порядка ${J^1 (\mathbb{R}^2)}$ с координатами
 $x,y,z,p,q$, где $x,\, y$ имеют смысл независимых переменных, $z$
 --- смысл искомой функции $z(x,y)$, и $p,\, q$ --- смысл первых производных
 $z_x,\, z_y$ соответственно.
 На ${J^1 (\mathbb{R}^2)}$ имеется пространство внешних форм
 2-го порядка ${\Lambda^2 J^1 (\mathbb{R}^2)}$.

Имея в виду, что
 ${p = f_x}$,
 ${q = f_y}$,
 где $f(x,y)$
 --- некоторая
 функция на $\mathbb{R}^2$,
  получаем равенство
$$
\begin{array}{l}
 dx \wedge dq + dy \wedge dp =
 dx \wedge d (f_y) +
 dy \wedge d (f_x) =
 dx \wedge
 \left(
 f_{xy} \: dx +
 f_{yy} \: dy
 \right) +
 \medskip\\
 +
 dy \wedge
 \left(
 f_{xx} \: dx +
 f_{xy} \: dy
 \right) =
 f_{xy} \:
 dx \wedge dx +
  f_{yy} \:
 dx \wedge dy +
  f_{xx} \:
 dy \wedge dx +
 \medskip\\
 +
  f_{xy} \:
 dy \wedge dy =
  f_{yy} \:
 dx \wedge dy -
  f_{xx} \:
 dx \wedge dy =
 \left(
  f_{yy} -
  f_{xx}
 \right)
 dx \wedge dy
\end{array}
$$
и равенство
$$
\begin{array}{l}
 dp \wedge dq + dx \wedge dy =
 d
  (f_x)
  \wedge d
  (f_y)
  + dx \wedge dy =
 \left(
  f_{xx}
  \: dx +
  f_{xy}
  \: dy
 \right)
 \wedge
 (
 f_{xy}
  \: dx +
 \medskip\\
 +
  f_{yy}
  \: dy
 )
 +
 dx \wedge dy
 =
 f_{xx}
 f_{xy}
 \: dx \wedge dx
 +
  f_{xx}
  f_{yy}
 \: dx \wedge dy
 +
 f_{xy}
 f_{xy}
 \: dy \wedge dx
 +
 \medskip\\
 +
 f_{xy}
 f_{yy}
 \: dy \wedge dy
 +
 dx \wedge dy
 =
 \left(
  f_{xx}
  f_{yy}
  -
    \left(
      f_{xy}
    \right)
    ^2
 \right)
 \: dx \wedge dy
 +
 dx \wedge dy
 =
 \medskip\\
 =
 \left(
 \mbox{hess}\: f + 1
 \right)
 dx \wedge dy \, .
\end{array}
$$
 Здесь использованы следующие тождества алгебры внешних форм
 \cite{Бурбаки-Алгебра-1}:
$$
\begin{array}{l}
  \omega_1 \wedge (\omega_2 + \omega_3) =
  \omega_1 \wedge \omega_2 +
  \omega_1 \wedge \omega_3 \, ,
  \smallskip \\
  a \omega_1 \wedge \omega_2 = \omega_1 \wedge a \omega_2
  = a (\omega_1 \wedge \omega_2) \, ,
  \smallskip \\
  \omega \wedge \omega = 0 \, ,
\end{array}
$$
где $\omega,\, \omega_1,\, \omega_2,\, \omega_3$ --- внешние
формы, $a$ --- константа.

 Итак, форме
$$
 dx \wedge dq + dy \wedge dp
$$
 из
 ${\Lambda^2 J^1 (\mathbb{R}^2)}$
 мы поставили в соответствие форму
$$
 \left(
  f_{yy} -
  f_{xx}
 \right)
 dx \wedge dy
$$
 из
 ${\Lambda^2 (\mathbb{R}^2)}$,
 т.{\,}е. линейное волновое уравнение,
 а форме
$$
dp \wedge dq + dx \wedge dy
$$
 из
 ${\Lambda^2 J^1 (\mathbb{R}^2)}$
 --- форму
$$
 \left(
 \mbox{hess}\: f + 1
 \right)
 dx \wedge dy
$$
 из
 ${\Lambda^2 (\mathbb{R}^2)}$,
 т.{\,}е. простейшее уравнение Монжа-Ампера.

 Рассмотрим преобразование Ампера
$$
 \bar x = - p , \quad
 \bar y = y , \quad
 \bar z = z - p \: x , \quad
 \bar p = x , \quad
 \bar q = q .
$$
Оно является контактным, т.{\,}е. сохраняет форму
 ${dz - p\: dx - q \: dy:}$
$$
\begin{array}{l}
   d \bar z - \bar p \: d \bar x - \bar q \: d \bar y =
   d(z - p \: x) - x \: d(-p) - q \: dy =
  \smallskip \\
   =
   dz - dp \: x - p \: dx + x \: dp - q \: dy =
   dz - p \: dx - q \: dy \, .
\end{array}
$$
Преобразование Ампера переводит уравнение Монжа-Ампера
${\mbox{hess} \: z = -1}$ в линейное волновое уравнение
 ${z_{xx} - z_{yy} = 0}$, так как
$$
 d \bar p \wedge d \bar q
 +
 d \bar x \wedge d \bar y
 =
 dx \wedge dq
 +
 d(-p) \wedge dy
 =
 dx \wedge dq
 -
 dp \wedge dy
 =
 dx \wedge dq
 +
 dy \wedge dp
 \, .
$$
Преобразование Ампера переводит классическое решение ${z=xy}$
уравнения ${\mbox{hess} \: z = -1}$, т.{\,}е. двумерное
интегральное многообразие ${(u,v,uv,v,u)}$, в интегральное
многообразие ${(-v,v,0,u,u)}$, т.{\,}е. многозначное решение
волнового уравнения ${z_{xx} - z_{yy} = 0}$. Это последнее
интегральное многообразие при проектировании на плоскость
${\mathbb{R}^2 = (x,y)}$ дает прямую ${(-v,v)}$, а не двумерную
область, т.{\,}е. оно ни в одной своей точке не может
использоваться в качестве классического решения волнового
уравнения. $\square$


\label{references}